\documentclass[12pt]{article}

\usepackage{amsmath, amssymb, eucal, amscd, color}
\begin{document}

\newtheorem{theorem}{Theorem}
\renewcommand{\thetheorem}{}
\newtheorem{proposition}{Proposition}
\renewcommand{\theproposition}{}
\newtheorem{lemma}{Lemma}
\renewcommand{\thelemma}{}
\newtheorem{corollary}{Corollary}
\newtheorem{remark}{Remark}
\newcommand{\Def}{{\bf Definition.}\quad}
\newcommand{\Thm}{{\bf Theorem.}\quad}
\newcommand{\Prop}{{\bf Proposition.}\quad}
\newcommand{\Lem}{{\bf Lemma.}\quad}
\newcommand{\Ex}{{\bf Example.}\quad}
\newcommand{\Rem}{{\bf Remark.}\quad}
\newtheorem{example}{Example}
\newcommand{\ca }{\mathcal}
\newcommand{\CC}{{\mathbb{C}}}
\newcommand{\QQ}{{\mathbb{Q}}}
\newcommand{\RR}{{\mathbb{R}}}
\newcommand{\ZZ}{{\mathbb{Z}}}
\newcommand\Hom{\operatorname{Hom}}
\newcommand\End{\operatorname{End}}
\newcommand\Spec{\operatorname{Spec}}
\newcommand\CH{\operatorname{CH}}
\newcommand\Cone{\operatorname{Cone}}
\newcommand{\Gr}{\operatorname{Gr}}
\newcommand{\cA}{{\cal A}} \newcommand{\cB}{{\cal B}}\newcommand{\cC}{{\cal C}}
\newcommand{\cD}{{\cal D}}\newcommand{\cE}{{\cal E}}\newcommand{\cF}{{\cal F}}
\newcommand{\cG}{{\cal G}}\newcommand{\cH}{{\cal H}}\newcommand{\cI}{{\cal I}}
\newcommand{\cJ}{{\cal J}} \newcommand{\cK}{{\cal K}} \newcommand{\cL}{{\cal L}}\newcommand{\cM}{{\cal M}} \newcommand{\cN}{{\cal N}} \newcommand{\cO}{{\cal O}}\newcommand{\cP}{{\cal P}} \newcommand{\cQ}{{\cal Q}} \newcommand{\cR}{{\cal R}}\newcommand{\cS}{{\cal S}} \newcommand{\cT}{{\cal T}} 
\newcommand{\cZ}{{\cal Z}} \newcommand{\Z}{{\cal Z}}
\newcommand{\cU}{{\cal U}} \newcommand{\cV}{{\cal V}}
\newcommand{\bF}{{\mathbf{F}}}
\newcommand{\bU}{ {\mathbf{U}} }
\newcommand{\bK}{ {\mathbf{K}} }
\newcommand{\bL}{ {\mathbf{L}} }

\newcommand{\chZ}{ {\cal Z}} 

\newcommand{\al}{\alpha} \newcommand{\be}{\beta} \newcommand{\ga}{\gamma} 
\newcommand{\de}{\delta} \newcommand{\eps}{\epsilon} \newcommand{\la}{\lambda}
\newcommand{\De}{\Delta}

\newcommand{\BA}{{\mathbb{A}}} \newcommand{\BB}{{\mathbb{B}}}
\newcommand{\BC}{{\mathbb{C}}} \newcommand{\BD}{{\mathbb{D}}}
\newcommand{\BE}{{\mathbb{E}}} \newcommand{\BF}{{\mathbb{F}}}
\newcommand{\BG}{{\mathbb{G}}} \newcommand{\BI}{{\mathbb{I}}} 
\newcommand{\BIc}{{{I}}}
\newcommand{\bbI}{{\mathbb{I}}}
\newcommand{\BJ}{{\mathbb{J}}} \newcommand{\BK}{{\mathbb{K}}} 
\newcommand{\BL}{{\mathbb{L}}} \newcommand{\BM}{{\mathbb{M}}} 
\newcommand{\BN}{{\mathbb{N}}} \newcommand{\BO}{{\mathbb{O}}} 
\newcommand{\BP}{{\mathbb{P}}} \newcommand{\BQ}{{\mathbb{Q}}}
\newcommand{\BR}{{\mathbb{R}}} \newcommand{\BS}{{\mathbb{S}}}
\newcommand{\BT}{{\mathbb{T}}} \newcommand{\BU}{{\mathbb{U}}} 
\newcommand{\BV}{{\mathbb{V}}} 
\newcommand{\BW}{{\mathbb{W}}} 
\newcommand{\Bf}{{\bold{f}}}  \newcommand{\Bg}{{\bold{g}}}
\newcommand{\Bu}{{\bold{u}}}\newcommand{\Bv}{{\bold{v}}}
\newcommand{\Bw}{{\bold{w}}}

\newcommand{\sq}{\square} 
\newcommand{\ul}{\underline}
\newcommand{\ulm}{ {\ul{m}\,} } \newcommand{\uln}{ {\ul{n}} } 
\newcommand{\ot}{\otimes}
\newcommand{\scirc}{\circ} 
\newcommand{\hts}{\hat{\otimes}}
\newcommand{\htimes}{\hat{\times}}
\newcommand{\tildetimes}{{\tilde\times}}
\newcommand{\tbullet}{^{\bullet\,\bullet\,\bullet} }
\newcommand{\dbullet}{^{\bullet\,\bullet} }

\newcommand{\ts}{\otimes}
\newcommand{\vphi}{{\varphi}}  
\newcommand{\ddd}{\cdots}
\newcommand{\hlim}{\operatornamewithlimits{holim}} %
\newcommand{\bd}{\partial}
\newcommand{\sAb}{s Ab}       %
\newcommand{\stm}{\times}
\newcommand{\Tot}{\operatorname{Tot}}
\newcommand{\Ker}{\operatorname{Ker}}
\newcommand{\Cok}{\operatorname{Cok}}
\newcommand{\bop}{\mathop{\textstyle\bigoplus}}
\newcommand{\bts}{\mathop{\textstyle\bigotimes}}
\newcommand{\bhts}{\mathop{\textstyle \widehat{\bigotimes}}}
\newcommand{\Cyc}{\Z} \newcommand{\C}{\Z}
\newcommand{\injto}{\hookrightarrow}
\newcommand{\isoto}{\overset{\sim}{\to}}
\newcommand{\J}{\cJ}
\newcommand{\ti}{\tilde}
\newcommand{\CHM}{{C\!H\!\cal M}}
\newcommand{\dotsigma}{\dot{\sigma} }
\newcommand{\dotdelta}{\dot{\delta} }
\newcommand{\dotvphi}{\dot{\vphi} }
\newcommand{\maruphi}{\varphi^{in}}

\newcommand{\term}{\operatorname{tm}}
\newcommand{\tm}{\operatorname{tm}}
\newcommand{\ctop}{\overset{\,\,\circ}}
\newcommand{\init}{\operatorname{in}}

\newcommand{\Symb}{\mathop{Symb}}
\newcommand{\bound}{\partial}
\newcommand{\proj}{\pi}
\newcommand{\surjto}{\twoheadrightarrow}
\newcommand{\sgn}{\operatorname{sgn}}
\newcommand{\tbar}{\,\rceil \!\! \lceil}
\newcommand{\dbar}{|} 
\newcommand{\tp}{\mbox{top}}
\newcommand{\bsigma}{\mbox{\boldmath $\sigma$}}
\newcommand{\smallbsigma}{\mbox{\boldmath $\scriptstyle{\sigma}$}}
\newcommand{\bbrace}[1]{\mbox{\boldmath \{ }#1
\mbox{\boldmath \} }}
\newcommand{\shift}{\operatorname{shift}}

\newcommand{\bvphi}{\mbox{\boldmath $\vphi$}}
\newcommand{\bpartial}{\mbox{\boldmath $\partial$}}
\newcommand{\incl}{\operatorname{incl}}
\newcommand{\bDelta}{\mbox{\boldmath $\Delta$}}
\newcommand{\diag}{\operatorname{diag}}

\newcommand{\bskip}{\bigskip}
\newcommand{\sskip}{\smallskip}
\newcommand{\SmProjS}{(Smooth/$k$, Proj/$S$)}

\newcommand{\Seq}{{\BS eq}}

\newcommand{\mapright}[1]{%
  \smash{\mathop{%
    \hbox to
     1cm{\rightarrowfill}}\limits^{#1} } } 
\newcommand{\mapr}{\mapright}
\newcommand{\smapr}[1]{%
  \smash{\mathop{%
    \hbox to 0.5cm{\rightarrowfill}}\limits^{#1} } } 
\newcommand{\maprb}[1]{%
  \smash{\mathop{%
    \hbox to 1cm{\rightarrowfill}}\limits_{#1} } } 
\newcommand{\mapleft}[1]{%
  \smash{\mathop{%
    \hbox to 1cm{\leftarrowfill}}\limits^{#1} } }
\newcommand{\mapl}{\mapleft}
\newcommand{\maplb}[1]{%
  \smash{\mathop{%
    \hbox to 1cm{\leftarrowfill}}\limits_{#1} } }
\newcommand{\mapdown}[1]{\Big\downarrow
  \llap{$\vcenter{\hbox{$\scriptstyle#1\, $}}$ }}
\newcommand{\mapd}{\mapdown}
\newcommand{\mapdownr}[1]{\Big\downarrow
  \rlap{$\vcenter{\hbox{$\scriptstyle#1\, $}}$ }}
\newcommand{\mapdr}{\mapdownr}
\newcommand{\mapup}[1]{\Big\uparrow
  \llap{$\vcenter{\hbox{$\scriptstyle#1\, $}}$ }}
 \newcommand{\mapu}{\mapup}
\newcommand{\mapupr}[1]{\Big\uparrow
  \rlap{$\vcenter{\hbox{$\scriptstyle#1\, $}}$ }}
 \newcommand{\mapur}{\mapupr}

\newcommand{\RefBlone}{[1]}
\newcommand{\RefBltwo}{{[2]}}
\newcommand{\RefBlthree}{{[3]}}
\newcommand{\RefCH}{{[4]}}
\newcommand{\RefFu}{{[5]}}
\newcommand{\RefHaoneone}{{[6]}}
\newcommand{\RefHaonetwo}{{[7]}}
\newcommand{\RefHaonethree}{{[8]}}
\newcommand{\RefHatwo}{{[9]}}
\newcommand{\RefHathree}{{[10]}}

\newcommand{\sss}[1]{\noindent(#1)}

\newcommand{\Multicpx}{0.1}
\newcommand{\Finiteorderedset}{0.3}
\newcommand{\Tensorproductcpx}{0.2}
\newcommand{\Phicpx}{0.4}

\newcommand{\Icov}{1.1}
\newcommand{\ZMU}{1.2}
\newcommand{\ZMUpush}{1.3}
\newcommand{\restrictedts}{1.4}
\newcommand{\Seqfiberings}{1.5}
\newcommand{\distsubcpx}{1.6}
\newcommand{\ZMtensor}{1.7}
\newcommand{\prophtsZM}{1.8}
\newcommand{\htsZMU}{1.9}
\newcommand{\prophtsZMU}{1.10}

\newcommand{\secfuncpx}{2}
\newcommand{\ZXU}{\secfuncpx.1}
\newcommand{\ZXUgeneral}{\secfuncpx.2}
\newcommand{\cFIcJ}{\secfuncpx.3}
\newcommand{\cFIcJSigma}{\secfuncpx.4}
\newcommand{\dimcyclcpx}{\secfuncpx.5}
\newcommand{\cFI}{\secfuncpx.6}
\newcommand{\propcFI}{\secfuncpx.6.1}
\newcommand{\lemmaFI}{\secfuncpx.7}
\newcommand{\cFISigma}{\secfuncpx.8}
\newcommand{\barcpx}{\secfuncpx.9}
\newcommand{\cpxFI}{\secfuncpx.10}
\newcommand{\phiphicompatible}{\secfuncpx.11}
\newcommand{\acyclicitysigma}{\secfuncpx.12}
\newcommand{\proofphiphicompatible}{\secfuncpx.13}
\newcommand{\proofacyclicitysigma}{\secfuncpx.14}
\newcommand{\symbols}{\secfuncpx.15}

\newcommand{\set}{s}
\newcommand{\Tee}{T}
\newcommand{\secdistcpx}{3}
\newcommand{\distcpxZM}{\secdistcpx.1}
\newcommand{\distcpxZMU}{\secdistcpx.2}
\newcommand{\almstdisj}{\secdistcpx.3}
\newcommand{\distcpxcFI}{\secdistcpx.4}
\newcommand{\distcpxFI}{\secdistcpx.5}

\newcommand{\secdiag}{4}
\newcommand{\diagcycle}{\secdiag.1}
\newcommand{\diagcycleProp}{\secdiag.2}
\newcommand{\deltadiag}{\secdiag.3}
\newcommand{\Deltadiag}{\secdiag.4}
\newcommand{\defdiag}{\secdiag.5}
\newcommand{\diagcomp}{\secdiag.6}

\title{\bf Cycle theory of relative correspondences}
\author{Masaki Hanamura }
\date{}
\maketitle

\begin{abstract} 
We establish a theory of complexes of 
{\it relative correspondences\/}. 
The theory generalizes the known theory of complexes of correspondences 
of smooth projective varieties. 
It will be applied in the sequel of this paper 
to the construction of the triangulated category 
of motives over a base variety.
\end{abstract}
\renewcommand{\thefootnote}{\fnsymbol{footnote}}

There is the theory of algebraic correspondences of 
smooth projective varieties for the Chow group and for the 
higher Chow group. 
We first recall the classical theory of correspondences for the 
Chow group. 
For smooth projective varieties $X$, $Y$ over a field $k$, let 
$\CH^r(X\times Y)$ be the Chow group of codimension $r$ cycles of 
$X\times Y$. An element of this group is said to be a  correspondence
from $X$ to $Y$. 
 Let $Z$ be another smooth projective variety. For  $u\in \CH^r(X\times Y)$ and $v\in \CH^s(Y\times Z)$, the composition $u\scirc v\in 
\CH^{r+s-\dim Y}(X\times Z)$ is defined by 
$$u\scirc v=p_{13\, *}(p_{12}^*u\cdot p_{23}^*v)$$
where for example $p_{12}$ is the projection from $X\times Y\times Z$ to $X\times Y$. One has associativity for composition: $(u\scirc 
v)\scirc w=u\scirc (v\scirc w)$.
The theory of {\it motives}
(to be precise, Chow motives) over $k$ is based on this
theory of correspondences. The basic idea is to consider the additive 
category where objects are smooth projective varieties, morphisms are given 
by correspondences, and composition given by composition of correspondences. 

Instead of the Chow group one can take the higher Chow group, 
and still has the theory of correspondences. 
Recall for a variety $X$ the cycle complex $(\Z^r(X, \bullet), 
\partial)$ is a chain complex where $\Z^r(X, n)$ is 
the free abelian group on the set of non-degenerate
irreducible  subvarieties $V$ of $X\times\sq^n$ meeting faces properly
(see \S 0 for details). The boundary map $\partial$
is given by restricting cycles  to 
 codimension one faces and taking an alternating sum. 
The homology of this complex is the group $\CH^r(X, n)$. 
One has $\CH^r(X, 0)=\CH^r(X)$. 

For elements $u\in  \CH^r(X\times Y, n)$ 
and $v\in \CH^s(Y\times Z, m)$ the composition $u\scirc v\in 
\CH^{r+s-\dim Y}(X\times Z, n+m)$ is defined by the same formula,
and one has associativity.  
 Indeed we can do this at the level of chain 
complexes. 
For $X$ and $Y$ smooth projective, 
$\Z^r(X\times Y, \bullet)$ is the complex of ``higher" correspondences from $X$ to $Y$.  For $u\in \Z^r(X\times Y, n)$ and $v\in \Z^s(Y\times Z, m)$ 
the pull-backs $p_{12}^*u$ and $p_{23}^*v$ may not meet properly
in $X\times Y\times Z\times\sq^{n+m}$. 
But according to a moving lemma the  subcomplex
$$\Z^r(X\times Y,\bullet)\hts \Z^s(Y\times Z, \bullet)$$
of $\Z^r(X\times Y,\bullet)\ts \Z^s(Y\times Z, \bullet)$  
generated by elements $u\ts v$, 
where $u$, $v$ are non-degenerate irreducible 
subvarieties  such that $p_{12}^*u$ and $p_{23}^*v$
meet properly, is a quasi-isomorphic subcomplex. 
For such $u$ and $v$, 
 the composition $u\scirc v\in \Z^{r+s-\dim Y}(X\times Z, \bullet)$ 
 is defined by the same formula,  yielding a map of complexes 
$$\rho:\Z^r(X\times Y, \bullet)\hts \Z^s(Y\times Z, \bullet)\to 
\Z^{r+s-\dim Y}(X\times Z, \bullet)\,\,.$$
If $W$ is a fourth smooth projective variety, the subcomplex
$\Z(X\times Y, \bullet)\hts \Z(Y\times Z, \bullet)\hts \Z(Z\times W, \bullet)$,
generated by $u\ts v\ts w$ such that the triple 
$p_{12}^*u, p_{23}^*v, p_{34}^*w$
 is properly intersecting on the four-fold product, is a 
quasi-isomorphic subcomplex.  For such $u, v, w$, one has $u\scirc v\scirc w
\in \Z(X\times W, \bullet)$ 
defined by $p_{14\, *}(p_{12}^*u\cdot p_{23}^*v\cdot p_{34}^*w)$, 
 and the following holds: $u\scirc v\scirc w=(u\scirc v)\scirc w=
u\scirc (v\scirc w)$. 

Complexes $\Z(X\times Y, \bullet)$ and the partially defined 
composition were used in the construction of 
a theory of the triangulated category of
{\it mixed motives} over $k$, see \RefHaoneone, \RefHaonetwo.  An object of the category is 
{\it a diagram of smooth projective varieties } which  consists of 
a sequence of smooth projective varieties and higher correspondences 
between them, subject to ``cocycle" conditions. 
\bigskip

We shall generalize this to relative correspondences. Let
$S$ be a quasi-projective variety over $k$. 
By a {\it smooth variety $X$ over $S$} we mean a smooth variety over $k$, 
equipped with a projective map to $S$ (the map $X\to S$ need not be 
smooth).  The class of such varieties we denote by
\SmProjS.
Let $X$ and $Y$ be smooth varieties over $S$. 
A natural choice for the complex
of correspondences from $X$ to $Y$ would be  $\Z_a(X\times_S Y, \bullet)$, 
the cycle complex of dimension $a$ cycles of 
the fiber product $X\times_S Y$. 
Since the variety $X\times_S Y$ is not smooth, 
we replace this with another complex of abelian groups $F(X, Y)$.
Concretely $F(X, Y)$ is the cone of the restriction map 
of the cycle complexes $\Z(X\times Y, \bullet)\to \Z(X\times Y-X\times_S Y, 
\cdot)$,  shifted by $-1$. (To be precise one should keep 
track of the dimensions of the cycle complex, which we ignore now.)
With this modification, we still does not have a partially defined 
composition map, but we achieve something close: 
\smallskip 

(1) There is  
 an injective quasi-isomorphism of complexes 
$ \Z(X\times_S Y, \bullet)\to F(X, Y)$.   

(2) If $Z$ is another smooth variety, projective over $S$, there is 
a quasi-isomorphic subcomplex 
$$\iota: F(X, Y)\hts F(Y, Z)\injto F(X, Y)\ts F(Y, Z)\,\,.$$

(3) There is another complex $F(X, Y, Z)$ and a surjective quasi-isomorphism 
$$\sigma: F(X, Y, Z)\to F(X, Y)\hts F(Y, Z)\,\,.$$

(4) There is a map of complexes 
$\vphi: F(X, Y, Z)\to F(X, Z)\,\,.$
\smallskip 

\noindent In the derived category at least, one has an induced map 
$F(X, Y)\ts F(Y, Z)\to F(X, Z)$ obtained by composing $(\iota\sigma)^{-1}$ and $\vphi$.  This map plays the role of composition. 
One should note, in contrast to the case $S=\Spec k$, 
there is no composition map defined on 
$F(X, Y)\hts F(Y, Z)$;  the composition $\vphi$ is defined only on 
$F(X, Y, Z)$. 

 The pattern persists for more than three varieties. 
For the formulation it is convenient to change the notation
as follows.  
In the above situation, write $X_1$, $X_2$ and $X_3$ in place of 
$X, Y, Z$; let 
$$F(X_1, X_2, X_3\tbar \{2\}):= F(X_1, X_2)\ts F(X_2, X_3)$$
and 
$$F(X_1, X_2, X_3| \{2\}):= F(X_1, X_2)\hts F(X_2, X_3)\,\,.$$
Then the maps are of the form 
$\iota_2: F(X_1, X_2, X_3\tbar \{2\})\injto F(X_1, X_2, X_3|\{2\})$, 
$\sigma_2: \! F(X_1, X_2, X_3)$
\newline $\to F(X_1, X_2, X_3\mid \{2\})$, and 
$\vphi_2: F(X_1, X_2, X_3)\to F(X_1, X_3)$. 
The generalization goes as follows. 
\smallskip 

(1) For each sequence of objects in \SmProjS,
 $X_1, \ddd ,X_n$ ($n\ge 2$), 
there corresponds  a complex   $F(X_1, \cdots , X_n)$. 
If $n=2$ there is an injective quasi-isomorphism 
$\Z(X_1\times_S X_2, \bullet)\to F(X_1, X_2)$. 

 For a subset of integers 
 $S=\{i_1, \cdots, i_{a-1}\}\subset (1, n)$, let 
$i_0=1$, $i_a=n$ and 
$$F(X_1, \cdots , X_n\tbar S):=F(X_{i_0}, \cdots, X_{i_1})\ts
F(X_{i_1}, \cdots, X_{i_2})\ts\cdots\ts
F(X_{i_{a-1}}, \cdots, X_{i_a})\,\,.
$$
This is an $a$-tuple complex. 
There is  an $a$-tuple complex  $F(X_1, \cdots , X_n| S)$ and an injective 
quasi-isomorphism 
$$\iota_S:F(X_1, \cdots , X_n| S)\injto F(X_1, \cdots , X_n\tbar  S)\,\,.$$
We assume  $F(X_1, \cdots , X_n\dbar \emptyset)
=F(X_1,\cdots, X_n)$.

(2) For $S\subset S'$ there is a surjective quasi-isomorphism 
of multiple complexes 
$$\sigma_{S\, S'}: F(X_1, \cdots , X_n\dbar S)\to 
F(X_1, \cdots , X_n\dbar S')\,\,.$$
For $S\subset S'\subset S''$, $\sigma_{S\,S''}=
\sigma_{S'\,S''}\sigma_{S\,S'}$.  In particular we have 
$\sigma_S:=\sigma_{\emptyset\, S} : F(X_1, \cdots , X_n)\to 
F(X_1, \cdots , X_n\dbar S)$. 

(3)  For $K=\{k_1, \ddd , k_b\}\subset (1, n)$ disjoint from $S$, a map of multiple complexes
$$\vphi_K: F(X_1, \ddd, X_n | S)\to F(X_1, \ddd,\widehat{X_{k_1}}, 
\ddd , \widehat{X_{k_b}}, \ddd ,  X_n | S)\,\,.$$
If $K$ is the disjoint union of $K'$ and $K''$, one has 
$\vphi_K=\vphi_{K'}\vphi_{K''}$. 

(4) If $K$ and $S'$ are disjoint $\sigma_{S\, S'}$ and $\vphi_K$ commute.  
\smallskip

In \S 1 and 2 of this paper, we define the complexes  $F(X_1, \cdots, X_n)$ as above for
a sequence of smooth quasi-projective varieties $X_1, \cdots, X_n$, each 
equipped with a projective map to the base variety $S$. 
We now explain the ideas for the construction in case $n\le 3$.

In \S 1, 
given a smooth variety $M$ and an open covering $\cU=\{U_i\}_{i\in I}$ (indexed by a finite totally ordered set $I$)
of an open set $U\subset M$, 
we define a complex $\C(M, \cU)$ (called the \v Cech cycle complex) which is
quasi-isomorphic to the cycle complex $\Z(A, \bullet)$ of $A=M-U$. 
If $\cU=\{U\}$, the covering consisting of $U$ only, 
$\Z(M, \cU)$ is the cone of the restriction map $\Z(M)\to \Z(U)$, 
shifted by $-1$.  In general one replaces $\Z(U)$ by $\Z(\cU)$, 
the \v Cech complex with respect to the covering. 

Assume $M'$ is another smooth variety, $\cU'$ a finite ordered
 open covering of $U'\subset M'$; assume also there are smooth maps 
 $q:M\to Y$ and $q':M'\to Y$. 
 Let $M\times_Y M'$ be the fiber product and $p: M\times_Y M'\to M$, 
$p':M\times_Y M'\to M'$ be the projections. 
One has the covering $p^{-1}\cU\amalg {p'}^{-1}\cU'$ of
the open set $p^{-1}U\cup
{p'}^{-1}U'$ of $M\times_Y  M'$. 
For $u\in\Z(M, \cU)$ and $v\in\Z(M', \cU')$ one has the pull-backs
$p^*u\in \Z(M\times_Y M', p^{-1}\cU)$ and 
${p'}^*v\in \Z(M\times_Y M', {p'}^{-1}\cU')$, and if 
they meet properly, their {\it product} is defined as an element 
of $\Z(M\times_Y M', p^{-1}\cU\amalg {p'}^{-1}\cU')$. 
Let $\Z(M,\cU)\hts \Z(M', \cU')\subset \Z(M, \cU)\ts \Z(M', \cU')$
be the subcomplex generated by $u\ts v$ for such $u, v$; 
this is called the {\it restricted tensor product}, and shown to be a quasi-isomorphic
subcomplex. The product gives  a map of complexes
$$\rho:  \Z(M, \cU)\hts \Z(M', \cU')\to \Z(M\times_Y M', 
p^{-1}\cU\amalg {p'}^{-1}\cU')\,\,.$$

If $p: M\to N$ is a projective map, $\cV$ a covering of an open set
of $V\subset N$, then $p^{-1}\cV$ is an open covering of $p^{-1}V
\subset M$, and there is the projection map 
$p_*: \Z(M, p^{-1}\cV)\to \Z(N, \cV)$. 

If we apply this to $A=X\times_SY \subset M=X\times Y$ and 
the covering consisting only of $U_{12}:=M-A$, one obtains 
the complex $\Z(X\times Y, \{U_{12}\})$. If we set $F(X, Y)$ to 
be this complex our problem is partially solved. 
If $Z$ is another variety 
over $S$, one has $F(Y, Z)=\Z(Y\times Z, \{U_{23}\})$
with $U_{23}=Y\times Z-Y\times_S Z$, and 
there is the product map 
$$\rho: \Z(X\times Y, \{U_{12}\})\hts \Z(Y\times Z, \{U_{23}\})
\to \Z(X\times Y\times Z, \{p_{12}^{-1}(U_{12}), p_{23}^{-1}(U_{23})\}\, )
\,\,.$$
The problem remains, since  from the target of $\rho$ there is no 
projection $p_{13\, *}$ to the cycle complex $\Z(X\times Z, \{U_{13}\})$
where $U_{13}=X\times Z-X\times_S Z$. 

One notices here  that there is the restriction  map 
$$r: \Z(X\times Y\times Z, \{U_{123}\}) \to \Z(X\times Y\times Z, 
\{p_{12}^{-1}(U_{12}), p_{23}^{-1}(U_{23})\}\, )\,\,,$$
where $U_{123}:= X\times Y\times Z- X\times_S Y\times_S Z$, since 
$U_{123}$ contains both $p_{12}^{-1}(U_{12})$ and $ p_{23}^{-1}(U_{23})$. 
The map $r$ is a quasi-isomorphism, since both complexes are 
quasi-isomorphic to $\Z(X\times_S Y\times_S Z)$. 
Assume for simplicity $Y$ is projective. One then defines the  projection 
along $p_{13}$ as the composition
$$p_{13\,* }:\Z(X\times Y\times Z, \{U_{123}\})\to
\Z(X\times Y\times Z, \{p_{13}^{-1} U_{13}\}\})\to 
 \Z(X\times Z, \{U_{13}\})=F(X, Z)\,\,.$$
Here the first map is the restriction, which is defined since 
$U_{123}\supset p_{13}^{-1} U_{13}$, and the second map is the projection 
along $p_{13}$. 
Consider now the double complex 
$$\begin{array}{ccc}
& &\Z(X\times Y, \{U_{12}\})\hts \Z(Y\times Z, \{U_{23}\}) \\
& &\mapdr{\rho} \\
\Z(X\times Y\times Z, \{U_{123}\})&\mapr{r} &\Z(X\times Y\times Z, 
\{p_{12}^{-1}(U_{12}), p_{23}^{-1}(U_{23})\})
\end{array}
$$
where the upper right corner and lower left corner are
 placed in degree 0,
and let $F(X, Y, Z)$ be the total complex. 
In other words it is the cone of $r+\rho$ shifted by $-1$. 
We thus have a homotopy commutative diagram of complexes 
$$\begin{array}{ccc}
F(X, Y, Z)&\mapr{} &\Z(X\times Y, \{U_{12}\})\hts \Z(Y\times Z, \{U_{23}\}) \\
\mapd{}& &\mapdr{\rho} \\
\Z(X\times Y\times Z, \{U_{123}\})&\mapr{r} &\Z(X\times Y\times Z, 
\{p_{12}^{-1}(U_{12}), p_{23}^{-1}(U_{23})\})
\end{array}
$$
The required properties are satisfied with this: 
The map $\sigma: F(X, Y, Z)\to F(X, Y)\hts F(Y, Z)$ is given 
by the projection to $\Z(X\times Y, \{U_{12}\})\hts \Z(Y\times Z, \{U_{23}\})$,
and the map $\vphi: F(X, Y, Z)\to F(X, Z)$ is obtained by composing 
the projection to $\Z(X\times Y\times Z, \{U_{123}\})$ with the map 
$p_{13\, *}$. 

The construction of the complexes  
$F(X_1, \cdots, X_n)$ for $n\ge 3$ and the 
maps $\sigma$, $\vphi$ consists of a systematic generalization of the
above.  In \S 1 we discuss the properties of the complexes $\Z(M, \cU)$
and their functorial properties. 
We also study the restricted tensor product of them. 
Using these, we construct in \S 2 the complexes $F(X_1, \cdots, X_n|S)$, the
maps $\iota$, $\sigma$, $\vphi$, and verify the properties 
(1)-(4) above.   
As in case $n=2, 3$ explained above, $F(X_1, \cdots, X_n)$ is 
built up of complexes of the form $\cZ(M, \cU)$ for appropriate
$(M, \cU)$. 

In \S 3 we study further properties of the complexes 
$F(X_1, \cdots, X_n|S)$ and show:
\smallskip 

$\bullet$\quad The multiple complex $F(X_1, \cdots, X_n)$ is a degreewise free
$\ZZ$-module on a given set of generators $\set F(X_1, \cdots, X_n)$.

$\bullet$\quad Let $S=\{i_1, \cdots, i_{a-1}\}\subset (1, n)$ and 
$I_j:=[i_{j-1}, i_j]$ for $j=1, \cdots, a$. 
For an element 
$$(\al_j)\in \prod_{j=1, \cdots, a} \set F(I_j)$$
we have the condition whether or not it is {\it properly intersecting}. (The condition reduces to the properly intersecting
property in relevant cycle complexes.)

$\bullet$\quad
The multiple complex $F(X_1, \cdots, X_n|S)$ has an alternative 
description as the subcomplex generated by $\al_1\ts\cdots\ts \al_{a}$ for $\al_j\in \set F(I_j)$, with  $(\al_j)$ properly intersecting.

$\bullet$\quad 
One can define a class of quasi-isomorphic subcomplexes of $F(X_1, \cdots, X_n\tbar S)$ called {\it distinguished subcomplexes}. 
The careful exposition of the details takes a large part of \S 3.

The complex $F(X_1, \cdots, X_n|S)$ is an example of a distinguished subcomplex. 
To give another typical example, 
let  $n<m$ and  given a sequence of varieties $X_1, \cdots, X_m$, 
a subset $S\subset (1, n)$, and an element $f\in \set F(X_n, \cdots, X_m)$. 
The subcomplex $[F(X_1, \cdots, X_n|S)]'$  of $F(X_1, \cdots, X_n\tbar S)$ generated by 
$\al_1\ts\cdots\ts\al_{a}$ such that $\{\al_1, \cdots,\al_{a}, f\}$
is properly intersecting is a distinguished subcomplex
of $F(X_1, \cdots, X_n\tbar S)$. One has by definition
$[F(X_1, \cdots, X_n|S)]'\subset F(X_1, \cdots, X_n|S)$. 
Then the map 
$$F(X_1, \cdots, X_n|S)\to F(X_1, \cdots, X_n|S)\ts F(X_n, \cdots, X_m)$$
sending $u$ to $u\ts f$ restricts to  give a map 
$[F(X_1, \cdots, X_n|S)]'\to F(X_1, \cdots, X_m|S\cup\{n\})$, as is obvious 
from the definitions. 
\bigskip

In \S 4 we construct the {\it diagonal cycles} and 
the {\it diagonal extension} which play the role of 
the identity. Let  $\Delta_X\in \Z(X\times_S X, 0)$ be the element 
given by the diagonal $X\subset X\times_S X$. Its image under the 
inclusion to $F(X, X)$ is also denoted $\Delta_X$; it has degree 0
and boundary zero. 
One can construct, for $n\ge 2$, an element 
$\bDelta_X(1, \cdots, n)\in F(\overbrace{X, \cdots, X}^{n})$ of 
degree 0 with boundary zero, satisfying the properties 
(compatibility with the maps $\sigma$ and $\vphi$) below. 
For the statement we introduce some notation. 
When $X$ is understood, for any
subset $I=\{j_1, \cdots, j_m\}\subset [1, n]$ set 
$F(I)=F(\overbrace{X, \cdots, X}^{m})$ and 
$\bDelta_X(I)=\bDelta_X(j_1, \cdots, j_m)\in F(I)$. 
For $S\subset (1, n)$ let  $\tau_S: F(X_1, \cdots, X_n)\to 
F(X_1, \cdots, X_n\tbar S)$ be the composition of $\sigma_S$ and $\iota_S$. 
\smallskip 

$\bullet$\quad One has $\bDelta_X(1, 2)=\Delta_X \in F(X, X)$. 

$\bullet$\quad If $S=\{i_1, \cdots, i_{a-1}\}
\subset (1, n)$, and $I_j=[i_{j-1}, i_j]$ for $j=1, \cdots, a$, one has 
$$\tau_S(\bDelta_X(1, \cdots, n)\, )=\bDelta(I_1)\ts\cdots\ts
\bDelta(I_{a})$$
in $F(X, \cdots, X\tbar S)=F(I_1)\ts\cdots\ts F(I_{a})$. 

$\bullet$\quad For $K\subset (1, n)$, 
$\vphi_K(\bDelta(1, \cdots ,n )\,)=\bDelta([1, n]-K)$.
\smallskip   

Let $\lambda: [1, m]\to [1, n]$ be a surjective map, and 
$X$ be a sequence of objects in \SmProjS\,
indexed by $[1, n]$. 
Let $\lambda^*X: i\mapsto X_{\lambda(i)}$ be the induced sequence 
of objects on $[1, m]$. We shall construct a map of complexes
called the diagonal extension
$$\operatorname{\lambda^*}:
F(X_1, \cdots, X_n)\to F(X_{\lambda(1)}, \cdots, X_{\lambda(m)})$$
and show that it is compatible with the maps $\sigma$ and $\vphi$ in 
an appropriate sense. 
\bigskip

The constructions and results in this paper show that the classes of
smooth varieties over $S$,  the complexes $F(X_1, \cdots, X_n)$
and the maps $\sigma$, $\vphi$ form a {\it quasi DG category}, 
a generalization of a DG category, to be introduced in another paper \RefHathree\,. 
A quasi DG category consists of a class of objects, 
complexes $F(X_1, \cdots, X_n|S)$ for a sequence of objects, 
maps $\sigma_{S\, S'}$, $\vphi_K$ and additional structure 
that are subject to a set of axioms. 
The axioms  is an abstraction of the properties 
verified for the relative cycle complexes. 

To be more precise, a {\it symbol} over $S$ is by definition a formal finite sum 
$\bop_{\al}(X_\al/S, r_\al)$
where $X_\al$ is a smooth variety over $S$ and $r_\al\in \ZZ$. 
To a finite sequence of symbols $K_1, \cdots, K_n$ ($n\ge 2$) 
and a subset $S\subset (1, n)$ one 
can associate a complex of abelian groups 
$F(K_1, \cdots, K_n|S)$; if $K_i=(X_i, r_i)$, then 
$F(K_1,\cdots, K_n|S)$ is the complex $F(X_1, \cdots, X_n|S)$,
the integers $r_i$ specifying the dimensions of the cycle complexes
involved. (If $n=2$, $F((X_1/S), (X_2/S)\,)$ is quasi-isomorphic to $\Z_{\dim X_2-r_2+r_1}(X_1\times_S X_2, \bullet)$.)
One has maps $\sigma_{S\, S'}$ and $\vphi_K$ for 
$F(K_1, \cdots, K_n|S)$ as well.  The class of symbols over $S$, 
the complexes $F(K_1, \cdots, K_n|S)$, the maps $\sigma_{S\, S'}$, 
$\vphi_K$, along with additional structure -- generating set for the 
complex, notion of properly intersecting elements, distinguished 
subcomplexes with respect to constraints, diagonal cycles and 
diagonal extension mentioned above -- constitute a quasi DG category. 
\bigskip

\setcounter{section}{-1}

\section{Basic notions.}

Subsections (\Multicpx) and (\Finiteorderedset) are used throughout this paper, 
(\Tensorproductcpx) and (\Phicpx) are needed in \S\S 2 and 3. 
\bigskip

\sss{\Multicpx} {\it Multiple complexes.}\quad
By a complex of abelian groups we mean a graded abelian group $A^\bullet$
with a map $d$ of degree one satisfying $dd=0$. 
If $u: A\to B$ and $v: B\to C$ are maps of complexes, 
we define $u\cdot v: A\to C$ by $(u\cdot v)(x)=v(u(x))$. 
So $u\cdot v$ is $v\scirc u$ in the usual notation. 
As usual we also write $vu$ for $v\scirc u$ (but not for $v\cdot u$).

A double complex $A=(A^{i, j}; d', d'')$ is a bi-graded abelian group with 
differentials $d'$ of degree $(1, 0)$, $d''$ of degree $(0, 1)$, satisfying
$d'd''+ d''d'=0$. Its total complex $\Tot(A)$ is the complex with 
$\Tot(A)^k=\bop_{i+j=k}A^{i, j}$ and the differential $d=d'+d''$. 
In contrast a ``double" complex 
$A=(A^{i, j}; d', d'')$ is a bi-graded abelian group with 
differentials $d'$ of degree $(1, 0)$, $d''$ of degree $(0, 1)$, satisfying
$d'd''=d''d'$. Its total complex $\Tot(A)$ is given by 
$\Tot(A)^k=\bop_{i+j=k}A^{i, j}$ and the differential $d=d'+(-1)^i d''$ on 
$A^{i, j}$.  
(Note that the totalization depends on the ordering of the 
two gradings; if we reverse the order, the corresponding totalization has 
differential $(-1)^j d'+ d''$ on $A^{i, j}$.)
A ``double" complex can be viewed as a double complex by taking the differentials 
to be $(d', (-1)^i d'')$ (or, when we reverse the order, $(\,(-1)^j d', d'')$).

Let $(A, d_A)$ and $(B, d_B)$ be complexes. Then 
$(A^{i, j}=A^{j}\ts B^i; d_A\ts 1, 1\ts d_B)$ is a ``double" complex.
Its total complex has 
differential $d$ given by 
$$d(x\ts y)= (-1)^{\deg y} dx\ts y + x\ts dy\,\,$$
(for the tensor product complex, 
we always take the reverse order of the gradings for the totalization). 
Note this differs from the usual convention. 

More generally for
$n\ge 2$ one has the notion of $n$-tuple complex and ``$n$-tuple" 
complex. An $n$-tuple (resp. ``$n$-tuple") complex is 
a $\ZZ^n$-graded abelian group $A^{i_1, \cdots, i_n}$ with differentials
$d_1, \cdots, d_n$, $d_k$ raising $i_k$ by 1, such that 
for $k\neq \ell$, $d_kd_\ell+ d_\ell d_k=0$ (resp. $d_kd_\ell=d_\ell d_k$).
An ``$n$-tuple" complex $A^{i_1, \cdots, i_n}$ is an $n$-tuple complex
with respect to the differentials $(d_1, (-1)^{\al_1} d_2, \cdots, (-1)^{\al_n} d_n)$
where $\al_k=\sum_{j<i} \deg_j$.  
In this way we turn an ``$n$-tuple" into an $n$-tuple complex. 
(As for double complexes, one may reverse the order of the gradings; in that case 
we will explicitly mention it.)
A single
complex $\Tot(A)$, called the total complex, is defined in either case. 

As a variant one can define partial totalization. 
To explain it, 
let $S_1, \cdots, S_m$ be an ordered set of non-empty
subsets of $[1, n]:=\{1, \cdots, n\}$
 such that $S_i\cap S_j=\emptyset$ for $i\neq j$ and $\cup S_i=[1, n]$.
Such data corresponds to a surjective map $f: [1, n]\to [1, m]$. 
Then, given an $n$-tuple complex (resp. an ``$n$-tuple" complex)
$A^{i_1, \cdots, i_n}$
 one can ``totalize" in degrees in $S_i$, and form an 
 $m$-tuple (resp. ``$m$-tuple") complex denoted 
$\Tot^{S_1, \cdots, S_m}(A)$ or $\Tot^f(A)$. 
Given surjective maps $f: [1, n]\to [1, m]$ and $g: [1, m]\to [1, \ell]$, one has 
$\Tot^g\Tot^f(A)=\Tot^{gf}(A)$.  
For example, if a subset $S=[k, \ell]\subset [1, n]$ is specified, 
one can ``totalize"
in degrees in $S$, so the result $\Tot^S(A)$ is an $m$-tuple 
(resp. ``$m$-tuple") complex, where $m=n-|S|+1$. 

For $n$ complexes $A_1^\bullet, \cdots, A_n^\bullet$, the tensor
product $A_1^\bullet\ts \cdots \ts A_n^\bullet$ is an ``$n$-tuple" 
complex.  For the tensor product, we view it as an $n$-tuple complex 
with respect to the 
reverse order of the gradings.  This $n$-tuple complex and its 
total complex will be still denoted $A_1^\bullet\ts\cdots\ts A_n^\bullet$. 

The only difference between $n$-tuple and ``$n$-tuple" complexes is 
that of signs. 
For the rest of this section, what we say for 
$n$-tuple complexes equally applies to ``$n$-tuple" complexes. 

If $A$ is an $n$-tuple complex and $B$ an $m$-tuple complex, and when 
$S=[k, \ell]\subset [1, n]$ with $m=n-|S|+1$ is specified, one can talk of 
maps of $m$-tuple complexes $\Tot^S(A)\to B$.  When the choice of 
$S$ is obvious from the context, we just say maps of multiple complexes
$A\to B$. For example if $A$ is an $n$-tuple complex and $B$ an 
$(n-1)$-tuple complex, for each set $S=[k, k+1]$ in $[1, n]$ one can 
speak of maps of $(n-1)$-tuple complexes $\Tot^S(A)\to B$; if $n=2$
there is no ambiguity. 
\bigskip 

(\Multicpx.1) {\it Multiple subcomplexes of a tensor product complex.}\quad
Let $A$ and $B$ be complexes. A  double subcomplex $C^{i, j}\subset
A^i\ts B^j$ is a submodule closed under the two differentials. 
If $\Tot(C)\injto \Tot(A\ts B)$ is a quasi-isomorphism, we say 
$C^{\bullet\, \bullet}$ is a quasi-isomorphic subcomplex. 
It is convenient to let $A^\bullet\hts B^\bullet$ denote
such a subcomplex. 
(Note it does not mean the tensor product of subcomplexes of $A$ and $B$.)
Likewise a quasi-isomorphic 
multiple subcomplex of $A_1^\bullet\ts\cdots \ts A_n^\bullet$
is denoted $A_1^\bullet\hts\cdots \hts A_n^\bullet$. 
\bigskip

\sss{\Tensorproductcpx} {\it  Tensor product of ``double'' complexes.}
\quad Let $A\dbullet=(A^{a,p}; d'_A, d''_A)$ be a ``double'' complex
 (so  $d'$ has degree $(1, 0)$, $d''$ has degree $(0,1)$, and $d'd'=0$, 
$d''d''=0$ and $d'd''=d''d'$). The associated total complex $\Tot(A)$
has differential $d_A$ given by $d_A= d'+(-1)^a d''$ on $A^{a,p}$. 
The association $A\mapsto \Tot(A)$ forms a functor.
Let $(B^{b,q}; d'_B, d''_B)$ be another ``double'' complex. 
Then the tensor product of $A$ and $B$  as ``double'' complexes, 
denoted $A\dbullet\times B\dbullet$, 
is by definition the ``double'' complex $(E^{c,r}; d'_E, d''_E)$, where
$$E^{c,r}=\bop_{a+b=c\,,p+q=r}A^{a,p}\ts B^{b,q}$$
and $d'_E=(-1)^b d'_A\ts 1+  1\ts d'_B$, $d''_E=(-1)^q 
d''_A\ts 1+ 1\ts 
d''_B$. 

The tensor product complex $\Tot(A)\ts \Tot(B)$ and 
the total complex of $A\dbullet\times B\dbullet$ are 
related as follows. There is an isomorphism 
of complexes 
$$u: \Tot(A)\ts \Tot (B)\to \Tot(A\dbullet\times B\dbullet)$$
given by 
$u=(-1)^{aq}\cdot id$ on the summand $A^{a,p}\ts B^{b,q}$.

Let $A\dbullet$, $B\dbullet$, $C\dbullet$ be ``double'' complexes. One has
an obvious isomorphism of ``double'' 
complexes
$(A\dbullet\times B\dbullet)\times C\dbullet=A\dbullet\times 
(B\dbullet\times C\dbullet)$; it is denoted 
$A\times B\times C$. 
We will often suppress the double dots for simplicity. 
The following diagram commutes:
$$
\begin{array}{ccc}
\Tot(A)\ts\Tot(B)\ts\Tot(C)&\mapr{u\ts 1} &\Tot(A\times B)\ts\Tot(C) \\
\mapd{1\ts u}& &\mapdr{u} \\
\Tot(A)\ts\Tot(B\times C)&\mapr{u} &\Tot(A\times B\times C)\,\,.
\end{array}
$$
The composition defines an isomorphism 
 $u: \Tot(A)\ts \Tot(B)\ts \Tot(C)\isoto \Tot(A\times B\times C)$.
 
One can generalize this to the case of tensor product of more
than two ``double'' complexes. 
If $A_1, \cdots, A_n$ are ``double'' complexes, there is an isomorphism
of complexes 
$$u_n:\Tot(A_1)\ts\cdots\ts\Tot(A_n)\to \Tot(A_1\times\cdots\times A_n)$$
which coincides with the above $u$ if $n=2$, and is in general a composition
of $u$'s in any order.
As in case $n=3$, one has commutative diagrams involving $u$'s;
we leave the details to the reader. 

Let $A$, $B$, $C$ be ``double'' complexes and 
 $\rho: A\dbullet\times B\dbullet\to 
C\dbullet$
be a map of ``double'' complexes, 
namely it is bilinear and 
for $\al\in A^{a,p}$ and $\be\in B^{b,q}$,
$$d' \rho (\al\ts\be)= \rho((-1)^b d'\al\ts\be+  \al\ts d'\be)$$
and 
$$d'' \rho (\al\ts\be)= \rho((-1)^q d''\al\ts\be+  \al\ts d''\be)\,\,.$$
Composing $\Tot(\rho): \Tot(A\times B)\to \Tot(C)$ with $u:
\Tot(A)\ts\Tot(B)\isoto \Tot(A\times B)$, one obtains the map 
$$\hat{\rho}: \Tot(A)\ts \Tot(B)\to \Tot(C)\,\,;$$
it is given given by $(-1)^{aq}\cdot 
\rho$ on the summand $A^{a,p}\ts B^{b,q }$. 

The same holds for a map of ``double'' complexes 
$\rho: A_1\times \cdots\times A_n\to C$.
\smallskip 

{\it Remark.}\quad One could discuss more general sign rules for 
the change of ordering of the set of gradings of multiple complexes. 
We have restricted our discussions to the case we will need in \S 2. 
\bigskip

\sss{\Finiteorderedset}
{\it Finite ordered sets, partitions and segmentations.}\quad
Let $I$ be a non-empty finite totally ordered 
set (we will simply  say a finite ordered set), so
$I=\{i_1, \cdots, i_n\}, i_1<\cdots <i_n$, where $n=|I|$. 
 The {\it initial} (resp. {\it terminal})
element of $I$ is $i_1$ (resp. $i_n$); let $\init(I)=i_1$, $\term(I)=i_n$.
If $n\ge 2$, let $\ctop{I}=I-\{\init(I), \term(I)\}$. 

If $I=\{i_1, \cdots, i_n\}$, a subset $I'$ of the form 
$[i_a, i_b]=\{i_a, \cdots, i_b\}$ is called a {\it sub-interval}. 

In the main body of the paper, for the sake of concreteness 
 we often assume $I=[1, n]=\{1, \cdots, n\}$, a subset of 
$\ZZ$. More generally a finite subset of $\ZZ$ is an example of a finite 
ordered set. 

A {\it partition} of $I$ is a disjoint decomposition into 
sub-intervals $I_1, \cdots, I_a$ such that there is a sequence of 
elememtns $\init(I)=i_0< i_1< \cdots < i_{a-1}< i_a=\term(I)$ so that 
$I_k=[i_{k-1}, i_k -1]$. 

 So far we have assumed $I$ and $I_i$ to be of cardinality $\ge 1$. 
In some contexts we allow only finite ordered sets 
with at least two elements. 
There instead of partition the following notion plays a role.
Given a subset of $\ctop{I}$, $\Sigma= \{i_1, \cdots, i_{a-1}\}$,
where $i_1<i_2<\cdots< i_{a-1}$, one has 
a decomposition of $I$ into the sub-intervals $I_1, \cdots, I_a$, where
$I_k=[i_{k-1}, i_k]$, with $i_0=\init(I)$, $i_a=\term(I)$. 
Thus the sub-intervals satisfy 
$I_k\cap I_{k+1}=\{i_k\}$ for $k=1, \cdots, a-1$. 
 The sequence of sub-intervals $I_1, \cdots, I_a$ is
  called the {\it segmentation } of $I$ corresponding to $\Sigma$. 
(The terminology is adopted  to distinguish it from 
 the  partition). 
 
 Finite ordered sets
of cardinality $\ge 1$ and partitions appear in connection with a 
sequence of fiberings. On the other hand, finite ordered sets
of cardinality $\ge 2$ and segmentations appear when we consider a sequence of varieties (or an associated
sequence of fiberings). See below.
\bigskip

\sss{0.4} {\it Sequence of fiberings} (\S 1).\quad  Let $n\ge 1$. A 
{\it sequence of fiberings} 
consists of smooth varieties $M_i$ ($1\le i\le n$) and $Y_i$ ($1\le i\le n-1$), together with smooth maps $M_i\to Y_i$ and $M_{i+1}\to Y_i$. 
For a sub-interval $I=[j, k]\subset [1, n]$ of cardinality $\ge 1$, 
one defines $M_I$ to be the fiber product $M_j\times_{Y_j}M_{j+1}\times\cdots
\times_{Y_{k-1}}M_k$. If $I_1, \cdots, I_c$ is a partition of $[1, n]$, 
then one has smooth varieties $M_{I_1}, \cdots, M_{I_c}$, which 
form a sequence of smooth varieties over appropriate $Y$'s. 
\bigskip

\sss{0.5}{\it Sequence of varieties} (\S 2).
\quad Let $n\ge 2$. A {\it sequence of smooth varieties} over $S$ 
is a set of smooth varieties 
$X_i$ indexed by $i\in [1, n]$, where each $X_i$ is equipped with 
a projective map to $S$. 
For a sub-interval $I=[j, k]$ of cardinality $\ge 2$, let 
$X_I$ be the direct product $\prod_{i\in I}X_i$. 
Given an segmentation $I_1, \cdots, I_c$ corresponding to $\Sigma=\{i_k\}$,
 the varieties $X_{I_t}$ and the projections to $X_{i_k}$ form a sequence of fiberings. 
\bigskip 



\section{ The \v Cech cycle complexes $\C(M, \cU)$}

Let $k$ be a fixed ground field. By a smooth variety over $k$ we mean 
a  smooth quasi-projective equi-dimensional variety over $k$. 
We refer to \S 0 for cycle complexes and finite ordered sets.
\bigskip 

\sss{\Icov} {\it $I$-coverings.}\quad Let $M$ be a smooth variety over $k$, $A\subset M$ a closed set, and $U:= M-A$.
Let $I$ be a finite ordered set. An open covering of $U$ indexed by $I$ 
(or just an {\em $I$-covering} of $U$) is a set of open sets $\cU =\{U_i\}_{i\in I}$, with $\cup_i U_i = U$.  It is also denoted by $(I, \cU)$. 

If $V\subset M$ is another open set, $J$ is another finite ordered set
 and $\cV=\{V_j\}_{j\in J}$ a $J$-covering of $V$,  a map of coverings $(I, \cU)\to (J, \cV)$, or just $\cU\to \cV$ for short, 
  is an order preserving
 map $\lambda: J\to I$ such that $U_{\la(j)}\supset V_j$ for $j\in J$. One then has $V\subset U$. We thus have the category of $I$-coverings of open sets
  of $M$, for varying $I$; it is denoted by $Cov(\cO(M))$. 
The subcategory of $I$-coverings of a given $U\subset M$ is denoted 
$Cov(U\subset M)$, or just $Cov(U)$. 

If $\cU$ is an $I$-covering of $U$ and $\lambda: J\to I$ is an order preserving map, define
$\lambda^*\cU$ to be the $J$-covering of $U'=\cup_j U_{\lambda (j)}$ given by 
$j\mapsto U_{\lambda(j)}$. There is a natural 
map of coverings $\lambda^*: 
(I, \cU )\to (J, \lambda^*\cU)$. For composition of maps 
$\lambda$, $\lambda^*\cU$ is contravariantly functorial. 
A map of coverings $\lambda: (I, \cU)\to (J, \cV)$ factors as 
$(I, \cU)\mapr{\lambda^*}(J, \lambda^*\cU)\to (J, \cV)$. 

If $\cU$ is an $I$-covering of $U$
and $\cU'$ an $I'$-covering of $U'$ then one has an
$I\amalg I'$-covering $\cU\amalg \cU'$ of $U\cup U'$. 
Here $I\amalg I'$ is ordered so that $i<i'$ for $i\in I$ and $i'\in I'$.

For the rest of this section, without so mentioning
 an indexing set $I$ is always finite ordered, and 
a map between them is always order preserving. 

The notion of coverings and maps can be defined for $I$ unordered or infinite.
For our purposes we restrict to finite ordered indexing sets. 
\bigskip 

\sss{\ZMU} {\it The complex $\Z(M, \cU)$.}\quad
For a quasi-projective variety $X$ and $s\in \ZZ$, let  $\Z_s(X, \cdot)$ denote the cubical cycle complex as recalled in \S 0; we will abbreviate it to $\Z_s(X)$ or $\Z(X)$,
as long as there is no confusion.  
We have the cycle complex $ \Z_s(U, \cdot)$ for an open set $U\subset M$. 

For an $I$-covering $\cU$ of $U$,  we define a complex denoted $\Z(M, \cU)$ 
to be the total complex associated to the double complex 
$Z^{\bullet \bullet}$ 
given as follows. Let 
$$Z^{0,q}=\Z(M, -q)\,\,,$$
and for $p\ge 0$,
$$Z^{p+1,q}=\bop_{i_0<\cdots <i_p} \Z(U_{i_0, \cdots, i_p}, -q)$$
where $U_{i_0, \cdots, i_p}=U_{i_0}\cap \cdots \cap U_{i_p}$. 
An element of $Z^{p+1,q}$ is $\al=(\al_{i_0, \cdots, i_p})$ with $\al_{i_0, \cdots, i_p}\in \Z(U_{i_0, \cdots, i_p}, -q)$.
It is convenient to set $U_{\emptyset}=M$, and when $p=-1$, we interpret
$(i_0, \cdots, i_p)=\emptyset$, so $\alpha_\emptyset\in \Z(M)$. 
With this convention, 
the differential $\delta$ of degree $(1, 0)$ is
given by sending 
$\al\in\bop \Z(U_{i_0, \cdots, i_p})$
to $$\delta(\al)_{i_0, \cdots , i_{p+1}}=\sum (-1)^{r}\,
\al_{i_0, \cdots, \widehat{
i_r}, \cdots, i_{p+1}}|U_{i_0, \cdots , i_{p+1}}\,\,.$$
The differential $\partial$ of degree $(0,1)$ is the boundary map of  
each cycle complex. The differential of the total complex is 
$\delta + (-1)^p\partial$ on $Z^{p,q}$. 
We shall denote it by $\partial$ by abuse of notation, as long 
as there is no confusion. 
Thus the total complex is of the form
$$\chZ_s(M, \cU)=\bigl[ \Z(M)\mapr{\delta} \bop\Z(U_{i_0}) \mapr{\delta}
 \bop_{i_0<i_1}\Z(U_{i_0\,i_1})\to \cdots \to\bop_{i_0<\cdots <i_p}
 \Z(U_{i_0, \cdots, i_p})\to \cdots ]\,\,.$$
Note that the natural map 
$$\iota: \Z_s(A)\to \chZ_s(M, \cU)$$
 is a quasi-isomorphism. This follows from the localization theorem
  for the cycle complex \RefBltwo. 
 
If $(J, \cV)$ covers $V$, and $\lambda:(I, \cU)\to (J, \cV)$ a map 
of coverings, 
there is the induced map of complexes (called the {\it restriction map})
$$\chZ(M, \lambda):\chZ(M, \cU)\to \chZ(M, \cV)\,\,;$$
thus we have a functor $\chZ(M, -)$ from the category $Cov(\cO(M))$ to 
$C(Ab)$. 
The following square commutes:
$$\begin{array}{ccc}
\chZ(M, \cU)&\mapr{\chZ(M, \lambda)} &\chZ(M, \cV) \\
\mapu{}& &\mapu{}  \\
\Z(A) &\mapr{} &\Z(B)\,\,.
\end{array}
$$
Here $B=M-V$, and the bottom is the map induced by inclusion. 

As special cases of $\chZ(M, \lambda)$, we have two types, as follows.
If $I=J$ and $\cU$, $\cV$ are coverings  such that 
$V_i\subset U_i$, one has the map 
$$\Cyc(M, \cU)\to \Cyc(M, \cV)\,\,.$$
If $\lambda: J\to I$ is a map and $\cV=\lambda^*\cU$, then 
$\lambda^*: (I, \cU)\to (J, \lambda^*\cU)$ induces the map 
$$\lambda^*: \Cyc(M, \cU)\to \Cyc(M, \lambda^*\cU)\,\,.$$
In general, the map $\chZ(M, \lambda)$ is a composition of 
these two types of maps. 
\bigskip 

\sss{\ZMUpush} {\it Push-forward and pull-back.}\quad 
If $p: M\to N$ is a projective map, $B\subset  N$ a closed set with  complement $V$ such that $p^{-1}V=U$ and $\cV\in Cov(V\subset N)$, then $p^{-1}\cV\in 
Cov(U\subset M)$ 
is defined in the obvious manner, 
and push-forward on  cycle  complexes  induces a map (also called the push-forward)
$$p_*: \chZ_s(M,p^{-1}\cV )\to \chZ_s(N, \cV)\,\,. $$
It is compatible with $p_*: \Z_s(A)\to \Z_s(B)$ via the maps $\iota$.   

The push-forward is covariantly functorial in $p$. 
It commutes with  the maps of functoriality for coverings:
For a map $\lambda: (I, \cU)\to (J, \cV)$ of coverings in $Cov(\cO(N))$, one has the induced map of coverings $\lambda: (I, p^{-1}\cU)\to (J, p^{-1}\cV)$
in $Cov(\cO(M))$, and the following square commutes:
$$\begin{array}{ccc}
\Z(M, p^{-1}\cU)&\mapr{\Z(M, \lambda)} &\Z(M, p^{-1}\cV) \\
\mapd{p_*}& &\mapdr{p_*} \\
\Z(N, \cU)&\mapr{\Z(N, \lambda)} &\Z(N, \cV)\,\,.
\end{array}
$$
\bigskip 

Let $p:M\to N$ be a smooth map of relative dimension $d$. For $\cV\in Cov(V\subset N)$, the pull-back map 
$$p^*: \Z_s(N, \cV)\to \Z_{s+d}(M, p^{-1}\cV)$$
is defined. It is compatible with the pull-back map 
$$p^*:\Z_s(B)\to \Z_{s+d}(p^{-1}B)$$
($B$ is the complement of $V$) via the maps $\iota$. 
The pull-back is contravariantly functorial in $p$. It commutes with 
the functoriality maps for coverings $\cV$. 
\bigskip 

\sss{\restrictedts} {\it Restricted tensor product and the product map\/.}  \quad 
Let $M$, $M'$ and $Y$ be smooth varieties  and $q: M\to Y$, $q':
M'\to Y$ be {\it smooth} maps of varieties. 
Let $M\diamond M':= M\times_Y M'$ and $p: M\diamond M'\to M$, 
$p': M\diamond M'\to M'$ be the 
projections. 

\vspace*{0.5cm}
\hspace*{4cm}
\unitlength 0.1in
\begin{picture}( 12.0800, 12.2400)(  4.6000,-20.8000)
\put(4.6000,-16.0800){\makebox(0,0)[lb]{$M$}}%
\put(16.6800,-16.0800){\makebox(0,0)[lb]{$M'$}}%
\put(10.7400,-22.1000){\makebox(0,0)[lb]{$Y$}}%
\put(8.7300,-9.8600){\makebox(0,0)[lb]{$M\diamond M'$}}%
%
{\color[named]{Black}{%
\special{pn 8}%
\special{pa 1116 1046}%
\special{pa 652 1404}%
\special{fp}%
\special{sh 1}%
\special{pa 652 1404}%
\special{pa 716 1380}%
\special{pa 694 1372}%
\special{pa 692 1348}%
\special{pa 652 1404}%
\special{fp}%
}}%
%
{\color[named]{Black}{%
\special{pn 8}%
\special{pa 1658 1630}%
\special{pa 1194 1990}%
\special{fp}%
\special{sh 1}%
\special{pa 1194 1990}%
\special{pa 1258 1964}%
\special{pa 1236 1956}%
\special{pa 1234 1932}%
\special{pa 1194 1990}%
\special{fp}%
}}%
%
{\color[named]{Black}{%
\special{pn 8}%
\special{pa 1216 1046}%
\special{pa 1638 1394}%
\special{fp}%
\special{sh 1}%
\special{pa 1638 1394}%
\special{pa 1600 1336}%
\special{pa 1598 1360}%
\special{pa 1574 1368}%
\special{pa 1638 1394}%
\special{fp}%
}}%
%
{\color[named]{Black}{%
\special{pn 8}%
\special{pa 660 1640}%
\special{pa 1084 1990}%
\special{fp}%
\special{sh 1}%
\special{pa 1084 1990}%
\special{pa 1044 1932}%
\special{pa 1042 1956}%
\special{pa 1020 1962}%
\special{pa 1084 1990}%
\special{fp}%
}}%
\end{picture}%
\vspace{0.8cm}

Let $a, b\in \ZZ$, and $c=a+b
-\dim Y$. We define the subcomplex 
$$\Z_a(M)\hts \Z_b(M')\subset \Z_a(M)\otimes \Z_b(M')$$
to be the submodule generated by 
elements $\alpha\ts\beta\in\Z_a(M)\ts \Z_b(M')$
($\al$ and $\beta$ are assumed to be irreducible, non-degenerate) such that 
$p^*\al$ and ${p'}^*\beta$ meet properly in $M\diamond M'$, and 
that the product $p^*\al\cdot {p'}^*\beta\in \Z_c(M\diamond M')$. 
We define 
$$\al\scirc_Y\beta=\al\scirc \beta:=p^*\al\cdot {p'}^*\beta\,\,.$$
(We say briefly that the condition is that 
$\al\scirc\beta\in \Z(M\diamond M')$ be
defined.)

Then the following conditions are satisfied; (i) is
 non-trivial and will be proved later in this section, 
  and the 
rest are immediate from the definitions.  
\smallskip 

(i) The inclusion of the subcomplex is a quasi-isomorphism. 

(ii) There is a map of complexes
$$\rho_Y=\rho: \Z_a(M)\hts \Z_b(M')\to \Z_c(M\diamond M')$$
which sends $\al\ts\beta$ to $\al\scirc_Y \beta$. 

(iii) If $\pi: N\to M$ is a smooth map of dimension $d$, the pull-back 
$\pi^*\ts id: \Z_a(M)\ts \Z_b(M')\to \Z_{a+d}(N)\ts \Z_b(M')$
takes the subcomplex $\Z_a(M)\hts \Z_b(M')$ into $\Z_{a+d}(N)\hts \Z_b(M')$. 
We thus have the induced map of complexes 
$$\pi^*\ts id: \Z_a(M)\hts \Z_b(M')\to \Z_{a+d}(N)\hts \Z_b(M')\,\,.$$
This applies in particular to open immersions $\pi: N\to M$. 
Similar property holds for pull-backs in $M'$. 

If $\pi: M\to N$ is a projective map, the push-forward
$\pi_*\ts id: \Z_a(M)\ts \Z_b(M')\to \Z_{a}(N)\ts \Z_b(M')$
takes the subcomplex $\Z_a(M)\hts \Z_b(M')$ into $\Z_{a}(N)\hts \Z_b(M')$:
$$\pi_*\ts id: \Z_a(M)\hts \Z_b(M')\to \Z_{a}(N)\hts \Z_b(M')\,\,.$$
Similar property holds for push-forward in $M'$. 
\bigskip

To state the next level of generalization, 
 let $M_1\to Y_1\gets M_2\to Y_2\gets M_3$
 be a sequence of smooth varieties and 
smooth maps. 
We have $M_i\diamond M_{i+1}=M_i\times_{Y_i} M_{i+1}$ as before,
and $M_1\diamond M_2\diamond M_3=M_1\times_{Y_1}M_2\times_{Y_2} M_3$. Note
$M_1\diamond M_2\diamond M_3=(M_1\diamond M_2)\diamond M_3= M_1\diamond (M_2\diamond M_3)$. Let $p_i: M_1\diamond M_2\diamond M_3\to M_i$ be the projection. 
We assume that the projection $M_1\diamond M_2\diamond M_3\to 
M_1\times M_3$ is smooth. 
(This is called a sequence of fiberings on $[1, 3]$.)

\vspace*{0.5cm}
\hspace*{3cm}
\unitlength 0.1in
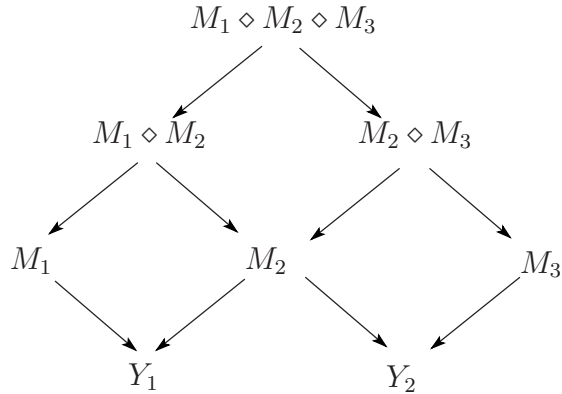
\begin{picture}( 26.7000, 18.8000)(  9.0000,-25.1000)
\put(9.0000,-20.2000){\makebox(0,0)[lb]{$M_1$}}%
\put(21.3000,-20.2000){\makebox(0,0)[lb]{$M_2$}}%
\put(15.2000,-26.3000){\makebox(0,0)[lb]{$Y_1$}}%
\put(13.3000,-13.6000){\makebox(0,0)[lb]{$M_1\diamond M_2$}}%
%
{\color[named]{Black}{%
\special{pn 8}%
\special{pa 1570 1440}%
\special{pa 1110 1810}%
\special{fp}%
\special{sh 1}%
\special{pa 1110 1810}%
\special{pa 1174 1784}%
\special{pa 1152 1778}%
\special{pa 1150 1754}%
\special{pa 1110 1810}%
\special{fp}%
}}%
%
{\color[named]{Black}{%
\special{pn 8}%
\special{pa 2130 2050}%
\special{pa 1670 2420}%
\special{fp}%
\special{sh 1}%
\special{pa 1670 2420}%
\special{pa 1734 2394}%
\special{pa 1712 2388}%
\special{pa 1710 2364}%
\special{pa 1670 2420}%
\special{fp}%
}}%
%
{\color[named]{Black}{%
\special{pn 8}%
\special{pa 1670 1440}%
\special{pa 2090 1800}%
\special{fp}%
\special{sh 1}%
\special{pa 2090 1800}%
\special{pa 2052 1742}%
\special{pa 2050 1766}%
\special{pa 2026 1772}%
\special{pa 2090 1800}%
\special{fp}%
}}%
%
{\color[named]{Black}{%
\special{pn 8}%
\special{pa 1140 2060}%
\special{pa 1560 2420}%
\special{fp}%
\special{sh 1}%
\special{pa 1560 2420}%
\special{pa 1522 2362}%
\special{pa 1520 2386}%
\special{pa 1496 2392}%
\special{pa 1560 2420}%
\special{fp}%
}}%
%
{\color[named]{Black}{%
\special{pn 8}%
\special{pa 2940 1470}%
\special{pa 2480 1840}%
\special{fp}%
\special{sh 1}%
\special{pa 2480 1840}%
\special{pa 2544 1814}%
\special{pa 2522 1808}%
\special{pa 2520 1784}%
\special{pa 2480 1840}%
\special{fp}%
}}%
%
{\color[named]{Black}{%
\special{pn 8}%
\special{pa 2450 2040}%
\special{pa 2870 2400}%
\special{fp}%
\special{sh 1}%
\special{pa 2870 2400}%
\special{pa 2832 2342}%
\special{pa 2830 2366}%
\special{pa 2806 2372}%
\special{pa 2870 2400}%
\special{fp}%
}}%
\put(28.7000,-26.4000){\makebox(0,0)[lb]{$Y_2$}}%
\put(27.2000,-13.6000){\makebox(0,0)[lb]{$M_2\diamond M_3$}}%
%
{\color[named]{Black}{%
\special{pn 8}%
\special{pa 3100 1470}%
\special{pa 3520 1830}%
\special{fp}%
\special{sh 1}%
\special{pa 3520 1830}%
\special{pa 3482 1772}%
\special{pa 3480 1796}%
\special{pa 3456 1802}%
\special{pa 3520 1830}%
\special{fp}%
}}%
%
{\color[named]{Black}{%
\special{pn 8}%
\special{pa 3570 2040}%
\special{pa 3110 2410}%
\special{fp}%
\special{sh 1}%
\special{pa 3110 2410}%
\special{pa 3174 2384}%
\special{pa 3152 2378}%
\special{pa 3150 2354}%
\special{pa 3110 2410}%
\special{fp}%
}}%
\put(35.7000,-20.4000){\makebox(0,0)[lb]{$M_3$}}%
%
{\color[named]{Black}{%
\special{pn 8}%
\special{pa 2220 830}%
\special{pa 1760 1200}%
\special{fp}%
\special{sh 1}%
\special{pa 1760 1200}%
\special{pa 1824 1174}%
\special{pa 1802 1168}%
\special{pa 1800 1144}%
\special{pa 1760 1200}%
\special{fp}%
}}%
%
{\color[named]{Black}{%
\special{pn 8}%
\special{pa 2420 840}%
\special{pa 2840 1200}%
\special{fp}%
\special{sh 1}%
\special{pa 2840 1200}%
\special{pa 2802 1142}%
\special{pa 2800 1166}%
\special{pa 2776 1172}%
\special{pa 2840 1200}%
\special{fp}%
}}%
\put(18.4000,-7.6000){\makebox(0,0)[lb]{$M_1\diamond M_2\diamond M_3$}}%
\end{picture}%
\vspace{0.8cm}

Define the subcomplex 
$$\Z(M_1)\hts\Z(M_2)\hts\Z(M_3)\subset \Z(M_1)\ts\Z(M_2)\ts\Z(M_3)$$
as follows. (In what follows we will not specify the dimensions $a_i$ 
for $\Z_{a_i}(M_i)$.)
It is generated by  $\al_1\ts\al_2\ts\al_3\in \Z(M_1)\ts\Z(M_2)\ts\Z(M_3)$
($\al_i$ are assumed irreducible, non-degenerate) 
satisfying the two conditions:\smallskip 

(i) The set of cycles 
$$\{p_1^*\al_1, p_2^*\al_2, p_3^*\al_3\}$$
is properly intersecting in $M_1\diamond M_2\diamond 
M_3$.  
(A set of cycles $\{z_1, \cdots, z_r\}$ on a smooth variety 
$M$ is defined to be properly intersecting if 
any for any sequence $1\le i_1 <\cdots < i_s\le r$,
the intersection $z_{i_1}\cap\, \cdots\, \cap z_{i_s}$ is empty or
has codimension equal to the sum of the codimensions of $z_{i_k}$. )

(ii) The intersections
$$p_1^*\al_1\,\cdot\, p_2^*\al_2, \,
\quad p_2^*\al_2\,\cdot\, p_3^*\al_3,
\quad \text{and}\quad  \,p_1^*\al_1\cdot p_2^*\al_2\cdot p_3^*\al_3
$$
are all in $\Z(M_1\diamond M_2\diamond M_3)$. 
\smallskip

One sets 
$$\al_1\underset{Y_1}{\scirc}
\al_2\underset{Y_2}{\scirc}\al_3:=
p_1^*\al_1\cdot p_2^*\al_2\cdot p_3^*\al_3\in \Z(M_1\diamond M_2\diamond M_3).
$$
We have $\al_1\scirc\al_2\in \Z(M_1\diamond M_2)$, 
and $(\al_1\scirc\al_2)\ts \al_3\in \Z(M_1\diamond M_2)\hts\Z(M_3)$. 
Similarly 
$\al_2\scirc\al_3\in \Z(M_2\diamond M_3)$ and 
$\al_1\ts(\al_2\scirc\al_3)\in \Z(M_1)\hts \Z(M_2\diamond M_3)$. 
Further, 
$$\al_1\scirc\al_2\scirc \al_3=(\al_1\scirc \al_2)\scirc \al_3=
\al_1\scirc (\al_2\scirc \al_3)\,\,.$$
The following statements hold, the only non-trivial one being the 
first assertion in (i).
\smallskip 

(i) The inclusion of the subcomplex is a quasi-isomorphism. 
There are also inclusions 
$$\Z(M_1)\hts\Z(M_2)\hts\Z(M_3)\subset(\Z(M_1)\hts\Z(M_2))\ts\Z(M_3)$$
 and 
$$\Z(M_1)\hts\Z(M_2)\hts\Z(M_3)\subset\Z(M_1)\ts(\Z(M_2))\hts\Z(M_3))\,\,.$$ 

(ii) There is a map of complexes 
$$\rho_{Y_1}:
 \Z(M_1)\hts\Z(M_2)\hts\Z(M_3)\to \Z(M_1\diamond M_2)\hts \Z(M_3)\,\,$$
which sends $\al_1\ts \al_2\ts\al_2$ to $(\al_1\scirc \al_2)\ts\al_3$.
Similarly there is a map  
$$\rho_{Y_2}:\Z(M_1)\hts\Z(M_2)\hts\Z(M_3)\to \Z(M_1)\hts \Z(M_2\diamond M_3)\,\,. $$
Also we have a map 
$$\rho_{Y_1Y_2}:
 \Z(M_1)\hts\Z(M_2)\hts\Z(M_3)\to \Z(M_1\diamond M_2\diamond M_3)\,\,$$
 which sends $\al_1\ts \al_2\ts\al_3$ to $\al_1\scirc\al_2\scirc\al_3$. 
One has $ \rho_{Y_1Y_2}=\rho_{Y_2}\rho_{Y_1}$, where 
$\rho_{Y_2}$ is the product map $\Z(M_1\diamond M_2)\hts \Z(M_3)\to \Z(M_1\diamond M_2\diamond M_3)$. Similarly $ \rho_{Y_1Y_2}=\rho_{Y_1}\rho_{Y_2}$. 

We may shorten the notation and write $\rho_i$ for $\rho_{Y_i}$ and 
$\rho_{12}$ for $\rho_{Y_1Y_2}$. We may alternatively write
$\rho([1, 2], [3])$ for $\rho_{Y_1}$, $
\rho([1], [2,3])$ for $ \rho_{Y_2}$, 
and $\rho([1, 3])$ for $\rho_{Y_1Y_2}$. 
Here $[j, k]$ denotes the set of integers between $j$ and $k$. 

 The inclusions  and the product maps are compatible 
 in the sense that the following square commutes:
$$\begin{array}{ccc}
\Z(M_1)\hts\Z(M_2)\hts\Z(M_3)&\mapr{\rho_{Y_1}} &\Z(M_1\diamond M_2)\hts \Z(M_3) \\
\mapd{\rm{incl}}& &\mapd{\rm{incl}} \\
\Z(M_1)\hts\Z(M_2)\ts\Z(M_3)&\mapr{\rho_{Y_1}\ts 1} &\Z(M_1\diamond M_2)\ts \Z(M_3)\,\,.
\end{array}
$$
Similarly for $\rho_{Y_2}$. 

(iii) If $\pi: N_1\to M_1$ is a smooth map of dimension $d$, the pull-back
$$\pi^*\ts 1\ts 1 : 
\Z(M_1)\ts \Z(M_2)\ts \Z(M_3)\to \Z(N_1)\ts\Z(M_2)\ts  \Z(M_3)$$
takes the subcomplex $\Z(M_1)\hts \Z(M_2) \hts\Z(M_3)$
into $\Z(N_1)\hts\Z(M_2)\hts\Z(M_3)$. So the map 
$\pi^*\ts 1\ts 1 : 
\Z(M_1)\hts \Z(M_2)\hts \Z(M_3)\to \Z(N_1)\hts\Z(M_2)\hts  \Z(M_3)$
is defined.
Similar property holds in each $M_i$.

If $\pi: M_1\to N_1$ is a projective map, the push-forward 
$$\pi_*\ts 1\ts 1 : 
\Z(M_1)\ts \Z(M_2)\ts \Z(M_3)\to \Z(N_1)\ts\Z(M_2)\ts  \Z(M_3)$$
takes the $\Z(M_1)\hts \Z(M_2) \hts\Z(M_3)$
into $\Z(N_1)\hts\Z(M_2)\hts\Z(M_3)$. 
Similar property holds in each $M_i$.

In the sequel of this section, we will give generalizations of 
these.  We first give the notion of a sequence of fiberings.
\bigskip 

\sss{\Seqfiberings} {\it Sequence of fiberings.}\quad 
Let $n\ge 1$.  A {\it sequence of fiberings} $(M, Y)$ indexed by $[1, n]$
consists of smooth varieties 
 $M_i$ ($1\le i\le n$) and  $Y_i$ ($1\le i\le n-1$)  
together with smooth maps $M_i\to Y_i$ and $M_{i+1}\to Y_i$.

\vspace*{0.5cm}
\hspace*{3cm}
\unitlength 0.1in
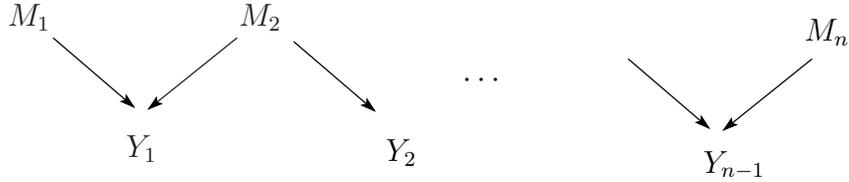
\begin{picture}(42.10,7.80)(3.60,-12.10)
\put(3.6000,-6.0000){\makebox(0,0)[lb]{$M_1$}}%
\put(15.7000,-6.0000){\makebox(0,0)[lb]{$M_2$}}%
\put(9.8000,-12.9000){\makebox(0,0)[lb]{$Y_1$}}%
%
\special{pn 8}%
\special{pa 1560 640}%
\special{pa 1100 1010}%
\special{fp}%
\special{sh 1}%
\special{pa 1100 1010}%
\special{pa 1164 984}%
\special{pa 1142 977}%
\special{pa 1139 953}%
\special{pa 1100 1010}%
\special{fp}%
%
\special{pn 8}%
\special{pa 600 640}%
\special{pa 1020 1000}%
\special{fp}%
\special{sh 1}%
\special{pa 1020 1000}%
\special{pa 982 941}%
\special{pa 980 965}%
\special{pa 956 972}%
\special{pa 1020 1000}%
\special{fp}%
%
\put(22.5000,-6.1000){\makebox(0,0)[lb]{}}%
\put(23.4000,-13.1000){\makebox(0,0)[lb]{$Y_2$}}%
%
\special{pn 8}%
\special{pa 1860 660}%
\special{pa 2280 1020}%
\special{fp}%
\special{sh 1}%
\special{pa 2280 1020}%
\special{pa 2242 961}%
\special{pa 2240 985}%
\special{pa 2216 992}%
\special{pa 2280 1020}%
\special{fp}%
\put(27.4000,-9.0000){\makebox(0,0)[lb]{$\cdots$}}%
\put(45.3000,-6.8000){\makebox(0,0)[lb]{$M_n$}}%
\put(40.0000,-13.8000){\makebox(0,0)[lb]{$Y_{n-1}$}}%
%
\special{pn 8}%
\special{pa 4570 750}%
\special{pa 4110 1120}%
\special{fp}%
\special{sh 1}%
\special{pa 4110 1120}%
\special{pa 4174 1094}%
\special{pa 4152 1087}%
\special{pa 4149 1063}%
\special{pa 4110 1120}%
\special{fp}%
%
\special{pn 8}%
\special{pa 3610 750}%
\special{pa 4030 1110}%
\special{fp}%
\special{sh 1}%
\special{pa 4030 1110}%
\special{pa 3992 1051}%
\special{pa 3990 1075}%
\special{pa 3966 1082}%
\special{pa 4030 1110}%
\special{fp}%
\end{picture}%
\vspace{0.5cm}

Given such, for any sub-interval $I=[j, k]\subset [1, n]$, let 
$$M_I=M_j\diamond M_{j+1}\diamond \cdots \diamond M_k\,\,.$$
There are projection maps $M_I\to M_{I'}$ for $I'\subset I$, and 
in particular, the map $p_i: M_I\to M_i $ for each $i\in I$.   

We will always assume the following condition:
\smallskip 

(*)\quad For each triple $1\le j <k<\ell\le n$, the product of projections 
$$M_{[j, \ell]} \to M_{[j, k-1]}\times M_{[k+1, \ell]}$$
\quad is smooth. 
\smallskip 

\noindent This is naturally 
satisfied for the sequence of fiberings arising from 
a sequence of varieties (see \S 2).
It simplifies some of our considerations 
(see (1.7.3) for example\,).  


Of course one can also define a sequence of fiberings indexed by any finite ordered set instead of $[1, n]$. 
\bigskip 

\sss{\distsubcpx} {\it Distinguished subcomplexes of cycle complexes.}\quad 
In \RefBlthree\, Bloch showed, for a smooth  variety $X$, the subcomplex of 
$\Z(X, \cdot)$ consisting of the cycles meeting a given set of subvarieties
of $X$ properly is a quasi-isomorphic  subcomplex.
We discuss a generalization of this. 

Let $X$ be a smooth quasi-projective variety. A finite set
$\{\al_1, \cdots , \al_n\}$  of irreducible 
subvarieties of $X$ is {\it properly intersecting\/} if 
for any subset $\{i_1, \cdots,  i_r\}$ of $\{1, \cdots , n\}$, 
the intersection $\al_{i_1}\cap\, \cdots\, \cap \al_{i_r}$ is empty or
has codimension equal to the sum of the codimensions of $\al_{i_k}$. 
A set of cycles $\{\al_1, \cdots , \al_n\}$ is properly intersecting if, for any choice of irreducible components 
$\be_i$ of $\al_i$, the set $\{\be_1, \cdots, \be_n\}$ is 
properly intersecting. 
In this case, the intersection cycle 
$\al_{i_1}\cdot\al_{i_2}\cdot\, \cdots\,
 \cdot \al_{i_r}$ is well-defined, independent of the 
order of taking intersections. 

Let $X_1, \cdots , X_r$ be smooth quasi-projective varieties.
 For a sequence of integers $s_1, \cdots, s_r$, 
 we have the cycle complexes
$\Z_{s_i}(X_i)$ and their tensor product $\Z_{s_1}(X_1)\ts\cdots\ts\Z_{s_r}(X_r)$.
We will consider a class of subcomplexes of this tensor 
product complex. 
From now we will usually drop the dimensions from the notation. 

Consider now triples of the form $(T; W; P)$ where
\smallskip 

(a) $T$ is another smooth variety,
and $W$ is an admissible cycle on $X_1\times \cdots \times X_r\times T\times \sq^{\ell} $ for some $\ell$
({\it admissible} means meeting faces properly).

(b) $P$ is a subset of $[1, r]$. 
\smallskip 

Let $W=\{(T_\lambda; W_\lambda; P_\lambda)\}_{\lambda=1, \cdots, c}$ be a finite set of triples;
thus $W_\lambda$ is an admissible cycle on $X_1\times \cdots \times X_r\times T_\lambda\times \sq^{\ell_\lambda} $. 
We define an $r$-tuple subcomplex, called 
the {\it distinguished subcomplex}
with respect to $W$, 
$$[\Z(X_1)\ts\cdots \ts\Z(X_r)]_{W}
\subset \Z(X_1)\ts\cdots \ts\Z(X_r)\,, $$
to be the subgroup generated by elements 
$$\al_1\ts\cdots\ts\al_r\in \Z(X_1, n_1)\ts\cdots \ts\Z(X_r, n_r)$$
where $\al_i$ are irreducible non-degenerate subvarieties 
of $X_i \times \sq^{n_i}$ satisfying the following condition:

 (PI)\quad For each $\lambda$,
the set of cycles 
$$\{\pi_i^*\al_i \,\,(\mbox{for $i\in P_\lambda$}),\,\, W_\lambda,\,
 \text{faces} \,\}\,\,$$
 is properly intersecting in $X_1\times\cdots\times X_r\times
 T_\lambda \times \sq^{n_1+\cdots +n_r}\times \sq^{\ell_\lambda}$. 
Here we employ the following obvious abuse of notation: 

\quad$\bullet$\quad $\pi_i$  denotes also 
 the projection $X_1\times \cdots \times X_r\times T_\lambda\times 
\sq^{n_1+\cdots +n_r}\times \sq^{\ell_\lambda}\to X_i\times\sq^{n_i}$;

\quad$\bullet$\quad
$W_\lambda$ denotes its pull-back by the projection
$X_1\times \cdots \times X_r\times T_\lambda\times 
\sq^{n_1+\cdots +n_r}\times \sq^{\ell_\lambda}\to 
X_1\times \cdots \times X_r\times T_\lambda\times\sq^{\ell_\lambda}$;

\quad$\bullet$\quad
a face means the pull-back of a face, namely a closed set 
of the form $X_1\times \cdots \times X_r\times T_\lambda\times F$
where $F$ is a face of $\sq^{n_1+\cdots +n_r}\times \sq^{\ell_\lambda}$. 
\smallskip 

{\it Remarks.}\quad (1) Let $T=T_1\times\cdots\times T_c$ and 
$p_\lambda : T\to T_\lambda$ be the projection.
Then the set of triples 
$$\{(T; p_\lambda^{-1}W_\lambda; P_\lambda)\}$$
gives the same distinguished subcomplex. 
Thus we will usually assume $T=T_\lambda$ in our discussion of 
distinguished subcomplexes. 

(2) The class of distinguished subcomplexes is closed under
finite intersections and tensor product. 
For the latter, if $W=\{(T; W_\lambda; P_\lambda)\}$
specifies a distinguished subcomplex of $\Z(X_1)\ts\cdots\ts\Z(X_n)
$, and $W'=\{(T'; W'_\mu ; P'_\mu)\}$
specifies a distinguished subcomplex of 
$\Z(X'_1)\ts\cdots\ts\Z(X'_m)$, then the ``union" of them,  
$$\{(T; W_\lambda; P_\lambda)\}\cup \{(T'; W'_\mu ; P'_\mu)\}$$
(where $W_\lambda$, say, stands for its pull-back by the projection
$X_1\times \cdots X_n\times X'_1\times\cdots\times X'_m
\to X_1\times \cdots X_n$, and $P_\lambda$ is viewed as a subset of 
$[1, n]\amalg [1', m']$)
specifies the tensor product of the two subcomplexes. 
\bigskip 

(\distsubcpx.1)\Thm{\it  The inclusion $[\Z(X_1)\ts\cdots \ts\Z(X_r)]_W\subset \Z(X_1)\ts\cdots \ts\Z(X_r)$ is a quasi-isomorphism.
}\smallskip 

This is proved in case $X_i$ are smooth projective in \RefHaonetwo, \S 1, generalizing \RefBlthree.  
(To be precise, the proof in \RefHaonetwo
\, is for the case $r\le 2$, 
but the proof is the same for any $r$). 

The case $X_i$ are smooth quasi-projective is similar. 
We sketch here an argument due to A. Krishna (as communicated to us by M. Levine).
Assume $r=1$, so we must show 
$\Z_W(X)\injto \Z(X)$ is a quasi-isomorphism. 
Take a projective closure $\bar{X}$ of $X$, let $Z=\bar{X}-X$, and 
consider the following commutative diagram:
$$\begin{array}{ccccc}
\Z(Z) &\mapr{}  &\Z(\bar{X}) &\mapr{} &\Z(X) \\
\Vert &       &\mapu{}       &      &\mapu{} \\
\Z(Z) &\mapr{}  &\Z_W(\bar{X}) &\mapr{} &\Z_W(X) \,\,.
\end{array}$$
By the localization theorem \RefBltwo\, we know the map 
$\Z(\bar{X})/ \Z(Z)\to \Z(X)$ is a quasi-isomorphism. 
The same proof shows 
$\Z_W(\bar{X})/ \Z(Z)\to \Z_W(X)$ is a quasi-isomorphism. 
The argument in \RefBlthree\, shows $\Z_W(\bar{X})\to \Z(\bar{X})$
is a quasi-isomorphism; although $\bar{X}$ is singular, the same
proof works since $W$ is contained in the smooth locus of $\bar{X}$. 
Hence one obtains the conclusion.
We leave the case $r\ge 2$ to the reader. 
\bigskip  

{\it Remarks.}\quad (1) For simplicity the defining condition (PI) may be phrased as follows, 
dropping $\pi^*$ and $\sq^{n_1+\cdots +n_r}\times \sq^{\ell_\lambda}$:
For each $\lambda$, the set 
$$\{ \al_i \,(\mbox{for $i\in P$}),\,\, W_\lambda,\,\text{faces} \,\}\,\,$$
 is properly intersecting in $X_1\times\cdots\times X_r \times T$. 
 
(2) The condition (PI) is equivalent to: For each $\lambda$ and each face $F$
of $\sq^{n_1+\cdots +n_r}\times \sq^{\ell_\lambda}$, the set of cycles 
$\{\al_i \,\,(i\in P)\,,\,\,W_\lambda \cap F\,\}$ is properly intersecting
in $X_1\times\cdots\times X_r \times T$. This follows from the following lemma.
\bigskip 

(\distsubcpx.2)\Lem{\it Let $X$ be a smooth variety, $\al_1, \cdots,\al_n$ be 
cycles on $X$, and $z_1, \cdots, z_m$ be properly intersecting cycles on 
$X$. 
Then the following are equivalent:

(i) $\{\al_1, \cdots,\al_n, z_1, \cdots, z_m\}$ is properly intersecting in
$X$. 

(ii) For each intersection $z_{j_1}\cap\cdots\cap z_{j_p}$, where 
$1\le j_1<j_2<\cdots<j_p\le m$, the set 
$$\{\al_1, \cdots,\al_n, z_{j_1}\cap\cdots\cap z_{j_p}\}$$
is properly intersecting in $X$.}
\bigskip

(\distsubcpx.3){\bf Example.}\quad 
Distinguished subcomplexes may be given as follows. 
Let $\{V_1, \cdots, V_k\}$ be a finite set of admissible cycles 
 $V_j$ on $X_1\times \cdots \times X_r\times T\times \sq^{\ell_j} $. We assume
the set
$$\{V_1, \cdots, V_k, \text{ faces}\}$$
 is {\it properly intersecting}. 
For a subset $P\subset [1, r]$, consider 
the subcomplex of $[\Z(X_1)\ts\cdots\ts\Z(X_r)]$
generated by elements
$$\al_1\ts\cdots\ts\al_r\in \Z(X_1, n_1)\ts\cdots \ts\Z(X_r, n_r)$$
where $\al_i$ are irreducible non-degenerate subvarieties 
satisfying the following condition:  the set
$$\{\pi_i^*\al_i\,\,(i \in P)\,\,, V_1, \cdots, V_k,
\text{faces} \,\}\,\,$$
 is properly intersecting in $X_1\times\cdots\times X_r\times T\times \sq^*$.  
 Then it is a distinguished subcomplex. 

Indeed, by the above lemma, the subcomplex coincides with the 
distinguished subcomplex with respect to 
$$W=\{(T; V_{j_1,\cdots, j_p}; P)\}$$
for partial intersections $V_{j_1,\cdots, j_p}=V_{j_1}\cap
\cdots\cap V_{j_p}$. 
 Note that the 
subcomplex differs from the distinguished subcomplex with respect to 
$\{(T; V_j; P)\}$, in which the proper intersection property is 
required with respect to each of $V_j$ separately.
\bigskip

\sss{\ZMtensor} {\it The subcomplex $\Z(M_1)\hts\cdots \hts\Z(M_n)$.}\quad 
Let $(M, Y)$ be a sequence of fiberings on $[1, n]$.   
We have  
$M_{[1, n]}=M_1\diamond \cdots \diamond M_n\subset M_1\times \cdots \times M_n$.
 There are  projections $p_{[1, n], i}=p_i: M_{[1, n]}\to M_i$, 
$\pi_i: M_1\times \cdots \times M_n\to M_i$,
and we have a commutative diagram:

\vspace*{0.5cm}
\hspace*{3cm}
\unitlength 0.1in
\begin{picture}(25.60,16.71)(2.40,-21.01)
\put(8.7000,-6.0000){\makebox(0,0)[lb]{$M_{[1, n]}$}}%
\put(22.9000,-6.0000){\makebox(0,0)[lb]{$M_1\times \cdots \times M_r$}}%
\put(12.1000,-12.1000){\makebox(0,0)[rt]{$M_I$}}%
\put(12.2000,-20.0000){\makebox(0,0)[rt]{$M_i$}}%
\put(16.7000,-6.0000){\makebox(0,0)[lb]{$\injto$}}%
%
\special{pn 8}%
\special{pa 1090 680}%
\special{pa 1090 1120}%
\special{fp}%
\special{sh 1}%
\special{pa 1090 1120}%
\special{pa 1110 1053}%
\special{pa 1090 1067}%
\special{pa 1070 1053}%
\special{pa 1090 1120}%
\special{fp}%
%
\special{pn 8}%
\special{pa 1080 1490}%
\special{pa 1080 1970}%
\special{fp}%
\special{sh 1}%
\special{pa 1080 1970}%
\special{pa 1100 1903}%
\special{pa 1080 1917}%
\special{pa 1060 1903}%
\special{pa 1080 1970}%
\special{fp}%
\put(7.9000,-18.2000){\makebox(0,0)[lb]{$p_{I, i}$}}%
%
\special{pn 8}%
\special{ar 840 1391 390 710  1.7659105 4.4792404}%
%
\special{pn 8}%
\special{pa 720 2061}%
\special{pa 720 2061}%
\special{fp}%
%
\special{pn 8}%
\special{pa 730 2071}%
\special{pa 790 2101}%
\special{fp}%
\special{sh 1}%
\special{pa 790 2101}%
\special{pa 739 2053}%
\special{pa 742 2077}%
\special{pa 721 2089}%
\special{pa 790 2101}%
\special{fp}%
\put(2.4000,-13.0000){\makebox(0,0)[lb]{$p_i$}}%
%
\special{pn 8}%
\special{pa 2800 700}%
\special{pa 1450 2060}%
\special{fp}%
\special{sh 1}%
\special{pa 1450 2060}%
\special{pa 1511 2027}%
\special{pa 1488 2022}%
\special{pa 1483 1999}%
\special{pa 1450 2060}%
\special{fp}%
\put(23.2000,-14.6000){\makebox(0,0)[lb]{$\pi_i$}}%
\put(6.7000,-9.8000){\makebox(0,0)[lb]{$p_{[1, n], I}$}}%
\end{picture}%
\vspace{0.5cm}

Let $\{I_1, \cdots , I_r\}$ be a partition of an interval $I=[j, k]$ 
into sub-intervals,
 namely  there is an increasing sequence $j=i_1< \cdots < i_{r+1}=k+1$
such that $I_a=[i_a, i_{a+1}-1]$. Then there are projections 
$$M_{I_a}\to Y_{i_{a+1}-1}\gets M_{I_{a+1}}\,\,.$$
So after renumbering
$$M'_a=M_{I_a},\quad Y'_a=Y_{i_{a+1}-1}$$
we have another sequence of fiberings indexed by $[1, r]$. 
Thus 
$M_{I_1}\diamond \cdots \diamond M_{I_r}$ makes sense and coincides with 
$M_I$. 
\bigskip 

In what follows we fix a sequence of integers $a_i\in \ZZ$, and 
take the complexes $\Z_{a_i}(M_i)$. To an interval $I=[j, k]$ we 
assign the integer
$$a_I=\sum_{i=j}^{k} a_i -\sum_{i=j}^{k-1} \dim Y_i$$
and take it as the dimension of the cycle complex $\Z(M_I)$. 
With this agreement we will drop the dimensions from the notation.
\bigskip 

(\ZMtensor.1)\Prop{ For a set of elements  $\al_i\in \Z(M_i, m_i), i\in [1, n]$, 
the following conditions are equivalent:
\smallskip 

(i) The set of cycles $\{p_i^*\al_i (i=1, \cdots, n),\, \text{faces} \,\}$
is properly intersecting in $M_{[1, n]}\times \sq^{m_1+\cdots m_n}$. 

(ii) The set of cycles $\{\pi_i^*\al_i (i=1, \cdots, n), M_{[1, n]},
\, \text{faces} \,\}$
is properly intersecting in $M_1\times \cdots\times M_n
\times \sq^{m_1+\cdots m_n}$.  
}\bigskip 

When this condition is satisfied we  will just say that the set 
$\{\al_i (i=1, \cdots, n),\, \text{faces} \,\}$
is properly intersecting in $M_{[1, n]}$. 
The equivalence  follows from the obvious
\bigskip 

(\ZMtensor.2)\Lem{\it Let $X$ be a smooth variety and $Y\subset X$ a smooth subvariety
(both assumed to be ). 
For a set of cycles $\al_1, \cdots, \al_n$ on $X$, the following are equivalent.\smallskip 

(i) The set $\{\al_1, \cdots, \al_n, Y\}$ is properly intersecting in $X$. 

(ii) The set $\{\al_1, \cdots, \al_n \}$ is properly intersecting in $X$, 
the intersection $\al_i\cap Y$ is proper for each $i$, and 
the set $\{\al_i\cdot Y\}$ is properly intersecting in $Y$. 
}\bigskip 

We define the (multiple) subcomplex 
$$\Z(M_1)\hts\cdots\hts\Z(M_n)\subset \Z(M_1)\ts\cdots\ts\Z(M_n)
\eqno{(\ZMtensor.a)}$$
to be the one generated by elements $\al_1\ts\cdots\ts\al_n$, with
each $\al_i$ irreducible non-degenerate, and $\{\al_1, \cdots, \al_n,\text{
faces}\}$ properly intersecting
in $M_{[1, n]}$. It is also denoted by $\widehat{\bts}_{i\in[1, n]}\Z(M_{i})$.

The proposition shows that it coincides with the distinguished subcomplex
with respect to the triple (see (\distsubcpx))
 $W=\{(T=pt; M_{[1, n]}; P=[1, n])\}$
 
Slightly more generally, given a subset $P$ of $[1, n]$, the subcomplex generated by elements $\al_1\ts\cdots\ts\al_n$, 
with 
$$\{p_i^*\al_i\,(i\in P),\, \text{
faces}\}$$
 properly intersecting, is a distinguished subcomplex 
given by the triple $W=\{(pt; M_{[1, n]}; P)\}$. 
We show that this coincides with a tensor product of complexes of the form (\ZMtensor.a). 
Let $\{I_1, \cdots, I_r\}$ be the partition of $[1, n]$, where 
each $I_j$ is either a maximal sub-interval (of $[1, n]$) contained in $P$, or a single point not in $P$. 
[Example: If $n=8$ and $P=\{1, 2, 3, 6, 7 \}$, then 
the partition consists of $[1, 3], \{4\}, \{5\}, [6, 7], \{8\}$.\,
]
With this we have:
\bigskip 

(\ZMtensor.3)\Prop{\it The above distinguished subcomplex 
coincides with the tensor product
$$(\bhts_{I_1} \Z(M_i)\,) \ts 
\cdots \ts (\bhts_{I_r} \Z(M_i)\,)\,\,.$$
}\smallskip 

{\it Proof.}\quad  The projection 
$$M_{[1, n]}\to M_{I_1}\times M_{I_2}\times\cdots\times M_{I_r}$$
is smooth by assumption (*) of (\Seqfiberings). 
So the collection $\{\al_i\,(i\in P),\, \text{
faces}\}$ is properly intersecting in $M_{[1, n]}$ if and only if 
for each $j$, $\{\al_i\,(i\in I_j),\, \text{faces}\}$ is properly intersecting in $M_{I_j}$. 
The assertion follows from this. 
\bigskip

Let  $\al_i\in \Z(M_i, m_i), i\in [1, n]$ be elements such that 
$\{\al_1, \cdots, \al_n, \text{faces}\}$ is properly intersecting in 
$M_{[1, n]}$. 
 Then for each interval $I=[j, k]\subset [1, n]$, the set $\{\al_j, \cdots, \al_k, 
\text{faces}\}$
is also properly intersecting, thus 
$$\al_j\ts\cdots\ts \al_k\in \Z(M_j)\hts\cdots\hts\Z(M_k)\,\,.$$
[To see this note the projection $M_{[1, n]}\to M_I$ is smooth, and the 
pull-back by a smooth map preserves the proper intersection property of cycles.
]
Thus the intersection in $M_I\times\sq^{m_j+\cdots+ m_k}$ 
$$(p_{I, j}^*\al_j)\cdot\,\,\,  \cdots\,\, \, \cdot( p_{I, k}^*\al_k)$$
is defined and $\in \Z(M_I)$. 
This is denoted by 
$$\al_j\underset{Y_j}{\scirc}\cdots\underset{Y_{k-1}}{\scirc}\al_k=\al_j\scirc
\cdots\scirc \al_k\,\,$$
or just by $\al_I$.

If $I_1, \cdots , I_r$ is a partition of $I$ then 
$$\al_{I_1}\ts \cdots \ts\al_{I_r}
\in \Z(M_{I_1})\hts\cdots\hts\Z(M_{I_r})$$
and 
$$\al_{I_1}\scirc \cdots\scirc \al_{I_r}=\al_I$$
in $\Z(M_I)$. 
So one has the product map
$$\rho(I_1, \cdots , I_r):
\Z(M_1)\hts\cdots\hts\Z(M_n)\to \Z(M_{I_1})\hts\cdots\hts\Z(M_{I_r})\,\,,$$
that maps $\al_1\ts\cdots \ts\al_n$ to $\al_{I_1}\ts \cdots\ts \al_{I_r}$.
It is a map of $r$-tuple complexes where the source 
$\Z(M_1)\hts\cdots\hts\Z(M_n)$ is viewed as an $r$-tuple complex by 
appropriate totalization (0.1). 
\bigskip 

\sss{\prophtsZM} {\it Functorial properties of $\bhts\Z(M_i)$.}\quad
In what follows by a complex we mean a multiple complex, and 
by a map of complexes we mean a map of multiple complexes
(with appropriate totalization as needed). 
We list the maps between the complexes $\bhts\Z(M_i)$.
\smallskip 

(1) One clearly has, for any partition 
$I_1, \cdots, I_r$ of $[1, n]$, the inclusion
$$\bhts_{[1, n]} \Z(M_i)\subset (\bhts_{I_1} \Z(M_i)\,) \ts 
\cdots \ts (\bhts_{I_r} \Z(M_i)\,)\,\,.$$
\smallskip

(2) To a partition $I_1, \cdots , I_r$ of $[1, n]$ there corresponds a 
map $$\rho(I_1, \cdots , I_r):
\bhts\Z(M_i)\to \bhts\Z(M_{I_a})\,,$$
defined in the previous subsection.
\smallskip 

(3) If $\pi_i: N_i\to M_i$ are smooth maps, there is the corresponding map 
$$\pi^*=\ts \pi_i^*: \bhts \Z(M_i)\to \bhts\Z(N_i)\,.$$
If $\pi_i: N_i\to M_i$ are projective maps, there is the corresponding map
$$\pi_*=\ts \pi_{i\,*}:  \bhts\Z(N_i)\to \bhts \Z(M_i)
\,.$$
\smallskip 

These maps are subject to the following compatibilities, all easily verified.
\smallskip 

$\bullet$\quad The inclusion (1) and the map $\rho$ of (2) commute
(the vertical maps are the inclusions):
$$\begin{array}{ccc}
\bhts_{[1, n]}\Z(M_i)&\mapr{\rho(I_1, \cdots, I_r)} &\bhts_{[1, r]} \Z(M_{I_a})\\
\mapd{}& &\mapd{}\\
\bhts_{I_1}\Z(M_i)\ts\cdots \ts \bhts
_{I_r}\Z(M_i)
&\mapr{\rho(I_1)\ts\cdots\ts\rho(I_r)}&\bts_{[1, r]}\Z(M_{I_a})\,\,.
\end{array}
$$

$\bullet$\quad The map $\rho$ in (2) safisfies associativity.
To formulate it, it will be convenient to think of the partition $I_1, \cdots, I_r$ of $[1, n]$
as follows. 
Let $f: [1, n]\to [1, r]$ be the map that takes $I_j$ to $j$.
This gives a one-to-one correspondence between the set of partitions of 
$[1,n]$ and order-preserving surjective maps from $[1, n]$ to another 
finite ordered set. We also write $\rho(f)$ for
 $\rho(I_1, \cdots, I_r)$. 
 
Let $f: [1, n]\to [1, r]$ and $g: [1, r]\to [1, \ell]$ be order-preserving surjections; 
let $I_1, \cdots, I_r$ be the partition of $[1, n]$ given by $f$, 
$K_1, \cdots, K_\ell$ the partition of $[1, r]$ given by $g$, 
and $J_1, \cdots, J_\ell$ the partition of $[1, n]$ given by $gf$. 
Note that $M'_j=M_{I_j}$ for $j\in [1, r]$ gives a sequence of fibrations on 
$[1, r]$; this sequence and the partition $K_1, \cdots, K_\ell$ gives 
another sequence of fibrations $M'_{K_k}$ on $[1, \ell]$, which coincides with 
the sequence of fibrations  $M_{J_k}$ on $k\in [1, \ell]$. 
Then the following diagram commutes:

\vspace*{0.3cm}
\hspace*{1cm}
\unitlength 0.1in
\begin{picture}( 29.5500, 14.3000)(  6.2000,-18.4000)
%
\put(6.2000,-9.5900){\makebox(0,0)[lb]{}}%
\put(10.0700,-7.1400){\makebox(0,0)[lb]{$\bhts_{[1, n]}\Z(M_i)$}}%
\put(35.7500,-7.0500){\makebox(0,0)[lb]{$\bhts_{[1, r]}\Z(M_{I_j})$}}%
\put(24.3500,-19.7000){\makebox(0,0)[lb]{$\bhts_{[1, \ell]}\Z(M_{J_k})$}}%
\put(23.2000,-5.4000){\makebox(0,0)[lb]{$\rho(f)$}}%
\put(33.5900,-13.7300){\makebox(0,0)[lb]{$\rho(g)$}}%
\put(14.4300,-13.7300){\makebox(0,0)[lb]{$\rho(gf)$}}%
%
{\color[named]{Black}{%
\special{pn 8}%
\special{pa 1960 600}%
\special{pa 3536 600}%
\special{fp}%
\special{sh 1}%
\special{pa 3536 600}%
\special{pa 3468 580}%
\special{pa 3482 600}%
\special{pa 3468 620}%
\special{pa 3536 600}%
\special{fp}%
}}%
%
{\color[named]{Black}{%
\special{pn 8}%
\special{pa 1714 814}%
\special{pa 2466 1724}%
\special{fp}%
\special{sh 1}%
\special{pa 2466 1724}%
\special{pa 2440 1660}%
\special{pa 2432 1684}%
\special{pa 2408 1686}%
\special{pa 2466 1724}%
\special{fp}%
}}%
%
{\color[named]{Black}{%
\special{pn 8}%
\special{pa 3572 796}%
\special{pa 2810 1716}%
\special{fp}%
\special{sh 1}%
\special{pa 2810 1716}%
\special{pa 2868 1676}%
\special{pa 2844 1674}%
\special{pa 2838 1652}%
\special{pa 2810 1716}%
\special{fp}%
}}%
\end{picture}%
\vspace{0.5cm} 

$\bullet$\quad The maps $\pi^*$, $\pi_*$ in (3) commute with the inclusion map (1). The pull-back map $\pi^*$ commutes with $\rho$ of (2).
\bigskip

\sss{\htsZMU} {\it Restricted tensor product of \v Cech cycle complexes.}\quad
Let $M$, $M'$ and $Y$ be as in (\restrictedts).
For open subsets $U\subset M$ and $U'\subset M'$, let $U\diamond U'= 
U\times_Y U'=p^{-1}(U)\cap {p'}^{-1}(U')\subset M\diamond M'$. If 
$A, A'$ are the complements of $U$, $U'$, $A\diamond A':=A\times_Y A'$ is the 
complement of $p^{-1}(U)\cup {p'}^{-1}(U')$. 

Given coverings $\cU\in Cov(U\subset M)$ and $\cV\in Cov(U'\subset M')$, 
we define the open covering
$$p^{-1}\cU\amalg {p'}^{-1}\cU' \in Cov(p^{-1}(U)\cup {p'}^{-1}(U')\subset M\diamond
M')\eqno{(\htsZMU.a)}$$
as follows. 
Let $I$ (resp. $I'$) be the index set for $\cU$ (resp. $\cU'$).
The disjoint union $I\amalg I'$ is a finite ordered set by $i<j$ for 
$i\in I$, $j\in I'$. The covering $(\htsZMU.a)$ is given by 
$$I\ni i\mapsto p^{-1}(U_i), \qquad I'\ni j\mapsto {p'}^{-1} (V_j)\,.$$

Let 
$$\Z(M, \cU)\hts \Z(M', \cU')\subset \Z(M, \cU)\ts \Z(M', \cU')$$
be the quasi-isomorphic subcomplex 
defined 
as the direct sum  
$$\bop \Z(U_{i_0,\cdots , i_p})\hts
\Z(V_{j_0,\cdots , j_q})
\subset \bop \Z(U_{i_0,\cdots , i_p})\ts
\Z(V_{j_0,\cdots , j_q})\,\,.$$
For the open covering $p^{-1}\cU\,\amalg\, {p'}^{-1}\cU'$,  the open set corresponding
to a sequence $(i_0\cdots i_pj_0\cdots j_q)$ is 
$$U_{i_0,\cdots , i_p}\diamond V_{j_0,\cdots , j_q}=p^{-1}(U_{i_0,\cdots , i_p})
\cap {p'}^{-1}(V_{j_0,\cdots , j_q})\,.$$
Then we have the map 
$$\rho_Y=\rho:\Z(M, \cU)\hts \Z(M', \cU')\to 
\Z(M\diamond M',  p^{-1}\cU\amalg {p'}^{-1}\cU')\,\,$$
which sends $\al\ts \al'\in \Z(M, \cU)\hts \Z(M', \cU')$ to 
$\al\scirc_Y \al'$ given by 
$$(\al\scirc_Y \al')_{i_0\cdots i_pj_0\cdots j_q}
=\al_{i_0\cdots i_p}\scirc_Y \al'_{j_0\cdots i_q}\,\,.$$
Here $\al$ consists of 
components $\al_{i_0\cdots i_p}\in \Z(U_{i_0,\cdots , i_p})$, 
where we interpret $\al_\emptyset\in \Z(M)$ if $p=-1$. 
Recall $\Z(M, \cU)$ and $\Z(M', \cU')$ are ``double" complexes, so 
their tensor product may also be viewed as a ``double" complex. 
One verifies that
$\rho_Y$ is a map of ``double" complexes. 

If $\lambda: \cU\to \cV$ and $\lambda': \cU'\to \cV'$ are maps 
of coverings of open sets of $M$, $M'$, respectively, there is an 
induced map of coverings $p^{-1}\cU\amalg {p'}^{-1}\cU'
\to p^{-1}\cV\amalg {p'}^{-1}\cV'$, and one easily verifies that 
the following diagram of ``double" complexes commutes:
$$\begin{array}{ccc}
\Z(M, \cU)\hts \Z(M', \cU')&\mapr{\rho} &\Z(M\diamond M',  p^{-1}\cU\amalg {p'}^{-1}\cU') \\
\mapd{\lambda\ts \lambda'}& &\mapdr{} \\
\Z(M, \cV)\hts \Z(M', \cV')&\mapr{\rho} &\phantom{\,.}\Z(M \diamond M',  p^{-1}\cV\amalg {p'}^{-1}\cV') \,.
\end{array}
$$

Since $A\diamond A'$ is the complement of $ p^{-1}U\cup
{p'}^{-1}U'$, the above $\rho$
 gives rise to a map in the derived category 
$$\rho_Y: \Z_a(A)\otimes \Z_b(A') \to \Z_{a+b-\dim Y}(A\diamond A')\eqno{}$$
which makes the following diagram commute:
$$\begin{array}{ccccc}
\Z(M, \cU)\otimes \Z(M', \cU') &\hookleftarrow &\Z(M, \cU)\hts \Z(M', \cU')
&\mapr{\rho}&\Z(M\diamond M',  p^{-1}\cU\amalg{p'}^{-1}\cU') \\
\mapu{ {\iota\ts\iota}}&&&&\mapur{{\iota}} \\
\Z(A)\otimes \Z(A')&&\mapr{\rho} &&\Z(A\diamond A')\,\,.
\end{array}$$

All this can be generalized as follows. 
Let $(M, Y)$ be a sequence of fiberings on $[1, n]$; let 
$U_i\subset M_i$  be open sets, and $A_i=M_i-U_i$. 
For an interval $I=[j, k]$, one has $M_I$ and the projections $p_i:
M_I\to M_i$. 
The complement of the union of $p_i^{-1}U_i$ for $i\in I$ is 
$A_I=A_j\diamond \cdots \diamond A_k$.

Assume given open coverings $\cU_i\in Cov (U_i\subset M_i)$.
One has a quasi-isomophic subcomplex 
$$\Z(M_1, \cU_1)\hts \cdots \hts\Z(M_n, \cU_n)\subset 
\Z(M_1, \cU_1)\ts \cdots \ts\Z(M_n, \cU_n)\,\,.$$
There is the open covering
$$p_1^{-1}\cU_1\amalg p_2^{-1}\cU_2\amalg \cdots\amalg p_n^{-1}\cU_n
\in Cov(\bigcup_i p_i^{-1}U_i\subset M_{[1, n]})\,,$$
which will be abbreviated to $\cU_{[1, n]}$. 

Let  $\al_i\in \Z(M_i, \cU_i), i\in [1, n]$ be elements such that 
$\{\al_1, \cdots, \al_n, \text{faces}\}$ is properly intersecting in 
$M_{[1, n]}$. 
 Then for each interval $I=[j, k]\subset [1, n]$, the set $\{\al_j, \cdots, \al_k, 
\text{faces}\}$
is also properly intersecting in $M_I$, so that
$$\al_I:=\al_j\scirc\cdots\scirc \al_k\in \Z(M_I,\cU_I)\,\,$$
is defined.

For a partition $I_1, \cdots I_r$ of $[1, n]$, we have open coverings 
$$\cU_{I_a}=\coprod_{i\in I_a} p_i^{-1}\cU_i\in
 Cov(\bigcup_{i\in I_a} p_i^{-1}U_i\subset M_{I_a})
\,\,,$$
and a map of ``double" complexes 
$$\rho(I_1, \cdots , I_r):
\Z(M_1, \cU_1)\hts \cdots \hts\Z(M_n, \cU_n)\to \Z(M_{I_1}, \cU_{I_1})\hts\cdots\hts
\Z(M_{I_r}, \cU_{I_r})$$
which sends $\al_1\ts\cdots \ts\al_n$ to $\al_{I_1}\ts \cdots\ts \al_{I_r}$.
(Recall from \S 0 that if $A_1\dbullet, \cdots, A_n\dbullet$ are ``double"
complexes, the tensor product $A_1\dbullet\ts\cdots\ts A_n\dbullet$ may 
be viewed as a ``double" complex.)
\bigskip 

\sss{\prophtsZMU} {\it Functorial properties of $\bhts \Z(M_i,\cU_i)$.}\quad
Continuing the previous subsection, we list the maps between the complexes 
$\Z(M_1, \cU_1)\hts \cdots \hts\Z(M_n, \cU_n)$.
\smallskip  

(1) For any partition 
$I_1, \cdots, I_r$ of $[1, n]$, there is inclusion
$$\bhts_{[1, n]} \Z(M_i, \cU_i)\subset (\bhts_{I_1} \Z(M_i, \cU_i)\,) \ts 
\cdots \ts (\bhts_{I_r} \Z(M_i, \cU_i)\,)\,\,.$$
\smallskip 

(2) The subcomplex is functorial in $\cU_i$. If $\cU_i\to \cV_i$ are maps in 
$Cov(M_i)$, there is an induced map 
$$\Z(M_1, \cU_1)\hts \cdots \hts\Z(M_n, \cU_n)
\to \Z(M_1, \cV_1)\hts \cdots \hts\Z(M_n, \cV_n)\,\,.$$
\smallskip 

(3) To a partition $I_1, \cdots , I_r$ of $[1,n]$ there corresponds a map 
$$\rho(I_1, \cdots , I_r): \Z(M_1, \cU_1)\hts \cdots \hts\Z(M_n, \cU_n)
\to \Z(M_{I_1}, \cU_{I_1})\hts\cdots\hts\Z(M_{I_r}, \cU_{I_r})\,.$$
\smallskip

(4) If $\pi_i: N_i\to M_i$ are smooth maps, there is the corresponding map 
$$\pi^*=\ts \pi_i^*: \bhts \Z(M_i, \cU_i)\to \bhts\Z(N_i, p_i^{-1}\cU_i)\,.$$
If $\pi_i: N_i\to M_i$ are projective maps, there is the corresponding map
$$\pi_*=\ts \pi_{i\,*}:  \bhts\Z(N_i, p_i^{-1}\cU_i)\to \bhts \Z(M_i, \cU_i)
\,.$$
\smallskip 

These maps satisfy the following compatibilities.
\smallskip 

$\bullet$\quad The maps (1) and (2) commute.

$\bullet$\quad The map $\rho(I_1,\cdots, I_r)$ in (3) safisfies associativity, 
formulated as in (\prophtsZM). 
It also commutes with the maps (1) and (2);
the commutativity with (1) means that the following diagram commutes
(where the vertical maps are inclusions):
$$\begin{array}{ccc}
\bhts_{[1, n]}\Z(M_i, \cU_i)&\mapr{\rho(I_1, \cdots, I_r)} &
\bhts_{[1, r]} \Z(M_{I_a}, \cU_{I_a})\\
\mapd{}& &\mapd{}\\
\bhts_{I_1}\Z(M_i, \cU_i)\ts\cdots \ts \bhts
_{I_r}\Z(M_i,\cU_i)
&\mapr{\rho(I_1)\ts\cdots\ts\rho(I_r)}&\bts_{[1, n]}\Z(M_{I_a, \cU_{I_a}})\,\,;
\end{array}
$$
the commutativity with the map (2) means the commutativity of the following diagram:
$$\begin{array}{ccc}
\Z(M_1, \cU_1)\hts \cdots \hts\Z(M_n, \cU_n)&\mapr{\rho(I_1, \cdots , I_r)}
&\Z(M_{I_1}, \cU_{I_1})\hts\cdots\hts\Z(M_{I_r}, \cU_{I_r}) \\
\mapd{}& &\mapd{} \\
\Z(M_1, \cV_1)\hts \cdots \hts\Z(M_n, \cV_n)&\mapr{\rho(I_1, \cdots , I_r)}
&\Z(M_{I_1}, \cV_{I_1})\hts\cdots\hts\Z(M_{I_r}, \cV_{I_r})\,\,.
\end{array}
$$

$\bullet$\quad The maps $\pi^*$, $\pi_*$ in (4) commute with the maps (1) and (2). The pull-back map $\pi^*$ commute with $\rho$ in (3). 
\bigskip

(5) The quasi-isomorphisms $\iota: \Z(A_i)\to \Z(M_i, \cU_i)$ induce a 
quasi-isomorphism
$$\ts \iota: \bts_{[1, n]}\Z(A_i)\to \bts_{[1, n]}\Z(M_i, \cU_i)\,\,.$$
Composing with the inverse of the inclusion $\bhts\Z(M_i, \cU_i)\injto
\bts\Z(M_i, \cU_i)$, 
 one obtains an isomorphism in the {\it derived category\/}
$$\iota:\bts_{[1, n]}\Z(A_i)\to \bhts_{[1, n]}\Z(M_i, \cU_i)\,\,.$$
(It is a slight abuse of notation to use the same $\iota$ for a map in 
the derived category.)  
 For a partition $I_1, \cdots, I_r$ of $[1, n]$, there is a unique map 
in the derived category
$$\rho(I_1, \cdots, I_r): \Z(A_1)\ts\cdots\ts\Z(A_n)\to \Z(A_{I_1})\ts
\cdots\ts\Z(A_{I_r})$$
which makes the following diagram commute:
$$\begin{array}{ccc}
\bhts_{[1,n]}\Z(M_i, \cU_i)&\mapr{\rho(I_1, \cdots, I_r)} &
\bhts_{[1,r]}\Z(M_{I_a}, \cU_{I_a}) \\
\mapu{\iota}&&\mapur{\iota} \\
{\bts_{[1,n]}}\Z(A_i)&\mapr{\rho(I_1, \cdots, I_r)} &
{\bts_{[1,r]}}\Z(A_{I_a})\,\,.
\end{array}
$$
This map $\rho(I_1, \cdots, I_r)$ also satisfies associativity (in the 
derived category), formulated as that for the map $\rho$ in (3).  
\bigskip


\section{ Function complexes $F(X_1, \ddd , X_n)$}

For notions regarding finite ordered sets, see \S 0.  
In particular, for an integer $n\ge 2$, we have the 
finite ordered set $[1, n]=\{1, \ldots, n\}$. 
We often use that to a finite ordered set $I$ of (cardinality $\ge 2$) 
and a subset $\Sigma\subset \ctop{I}$, 
there corresponds a segmentation of $I$. 
\bigskip 

\sss{\ZXU}{\it The complex $\Z(X_{[1, n]}, \cU(\cJ))$.}\quad 
Let $S$ be a quasi-projective variety and $X_1, \cdots , X_n$, $n\ge 2$, be smooth 
quasi-projective varieties, each equipped with a 
projective map to $S$ (we call such $X_i$ a {\it sequence of varieties} over $S$).
For a subset $I\subset [1, n]$, let
$X_{I}=\prod_{i\in I} X_i$ (product over $k$). 
  So $X_{[1, n]}=X_1\times \ddd\times X_n$.

For a non-empty subset $I\subset [1, n]$, let $X_{[1, n]}\isoto 
\underset{i\in I}{\prod}X_i \times \underset{i\not\in I}{\prod}X_i$ be the 
natural isomorphism (switching factors); define the closed subset 
$A_I\subset X_{[1, n]}$  by the Cartesian square (up to taking reduced structure)
$$\begin{array}{ccc}
X_{[1, n]} &\mapr{\sim} 
&\underset{i\in I}{\prod}X_i \times \underset{i\not\in I}{\prod}X_i \\
\mapu{}  &  &\mapu{} \\
{A_I} &\mapr{\sim}
&\underset{i\in I\,\,\,}{\prod_S}
X_i \times \underset{i\not\in I}{\prod}X_i 
\end{array}$$
where $\prod_S$ denotes fiber product over $S$. 
\smallskip 

$\bullet$\quad For example, if $I$ consists of a single element, $A_I= X_{[1, n]}$; if $I=\{1, 2\}$, $A_I= (X_1\times_S X_2)\times X_3\times \cdots \times X_n$; if $I=[1, n]$, $A_{[1, n]}=  X_1\times_S X_2\times_S\cdots \times_S X_n$. 

$\bullet$\quad If $I\subset I'$, then $A_I\supset A_{I'}$. 
For two subsets $I$ and $I'$ with non-empty intersection, 
$A_{I\cup I'}= A_I\cap A_{I'}$. 
\smallskip 

Let $U_I= X_{[1, n]}- A_I$. $U_{[1, n]}$ is the complement of $X_1\times_S X_2\times_S\cdots \times_S X_n$. 
If $I\subset I'$, then 
$U_I\subset U_{I'}$. If $I\cap I'$ is non-empty, 
$U_{I\cup I'}= U_I\cup U_{I'}$.

Let $\J$ be a subset of $ (1, n)=[2, n-1]$. 
If  $\J=\{j_1, \cdots , j_r\}$, the 
corresponding segmentation consists of sub-intervals
$J^k =[j_k, j_{k+1}]$ for $k=0, \ddd , r$ with $j_0=1$ and $j_{r+1}=n$. 
 To each $J^k$
there corresponds the closed set $A_{J^k} \subset X_{[1, n]}=X_1\times\cdots
\times X_n$ and its complement $U_{J^k}$. 
The intersection of $A_{J^k}$'s is $A_{[1, n]}$, 
and the union of $U_{J^k}$'s is $U_{[1, n]}$. We thus have a covering 
of $U_{[1, n]}$ indexed by $[0, r]$:
$$\cU(\cJ)=\{U_{J^0}, U_{J^1}, \cdots, U_{J^r}\}\,\,.$$

Taking $M=X_{[1, n]}$ and $\cU=\cU(\J)$ in the construction of the previous 
section, 
one obtains the complex
$$\C_s(X_{[1, n]}, \cU(\J))\,\,.$$
The differential of this complex is denoted by $d$, and when necessary we write
$\C_s(X_{[1, n]}, \cU(\J))^\bullet$ where the upper indexing 
is  the cohomological degree.
There is a natural quasi-isomorphism 
$$\iota:\cZ_s(A_{[1, n]})=\cZ_s(X_1\times_S\times \cdots \times_S X_n)\to 
\C_s(X_{[1, n]}, \cU(\J))\,\,.$$
We will fix (and usually drop) $s\in \ZZ$ for the dimension of the cycle complex. 
\bigskip 

Note in the discussion so far, one can replace $[1, n]$ by any 
 subset  $\BI$ (with cardinality $\ge 2$)
of $[1, n]$,  and $\cJ\subset (1,n)$ by a subset $\J\subset \ctop{\BI}$.
More specifically:
\smallskip 

$\bullet$\quad  One has the product 
$$X_{\BI}=\prod_{i\in \BI} X_i\,\,.$$
Associated to a subset $I\subset \BI$ is
a closed set $A_I\subset X_{\BI}$ and its complement $U_I$ (to be specific, we
write $A_{I\subset\BI}$ and $U_{I\subset\BI}$). In particular, $A_{\BI}$ is the fiber product of all $X_i$ over $S$, and $U_{\BI}$ its complement. 

$\bullet$\quad  
For a set $\J\subset \ctop{\BI}$ of cardinality $r$, 
there corresponds a set of intervals 
$$J^i=J^i(\cJ\subset\BI)\,\,,\quad 0\le i\le r\,\,,$$
of $\BI$. Thus we have an $[0,r]$-covering of $U_{\BI}$
$$\cU(\cJ)=\cU(\cJ\subset \BI)=\{U_{J^0}, \cdots, U_{J^r}\}\,\,.$$
This  gives us the complex $\C(X_{\BI}, \cU(\J\subset\BI))$ equipped with
 a quasi-isomorphism from $\Z(A_{\BI})$. 

$\bullet$\quad Note $U_{I\subset\BI}$ is an open set of $X_{\BI}$; it should
 be distinguished from the open set $U_{I\subset [1, n]}\subset X_{[1, n]}$. 
\bigskip 

{\it Assume now that $S$ is projective.} Then $X_i$ are also all projective. 
We have natural maps between the complexes above, the 
restriction and the projection.
\smallskip 

(1) For $\J\subset \J'$, one has
the {\it restriction map\/}, which is a quasi-isomorphism:
$$\C(X_{\BI}, \cU(\J)) \to \C(X_{\BI}, \cU(\J'))\,\,.\eqno{(\ZXU.a)}$$

To define it, assume $\BI=[1, n]$ for simplicity. 
Let $\cJ=\{j_1, \cdots, j_r\}$, $J^k=[j_k, j_{k+1}]$ for $0\le k\le r$ as above.
Similarly to  $\cJ'=\{j'_1, \cdots, j'_{r'}\}$
there corresponds a set of intervals 
$\{ {J'}^0, \cdots, {J'}^{r'}\}$, 
and the open covering 
$$\cU(\cJ')=\{U_{{J'}^0}, \cdots, U_{{J'}^{r'}} \}\,.$$
Define an order-preserving map $\lambda: [0, r']\to
[0,r]$ by the condition  $J'_t\subset J_{\lambda (t)}$ for $t\in [0, r']$. 
Since $U_{J'_t}\subset U_{J_{\lambda (t)}}$, one has a map of coverings 
$\lambda: \cU(\J)\to \cU(\cJ')$. 
It induces the map between the \v Cech cycle complexes as stated.

(2) For $\ell\in \ctop{\BI}-\J$,
we have the {\it projection\/} along $X_\ell$
$$\proj_{\ell}: \C(X_{\BI}, \cU(\J\subset \BI))\to
 \C( X_{\BI-\{\ell\}}, \cU(\J\subset \BI-\{\ell\}))\,\,.\eqno{(\ZXU.b)}$$

The definition in case $\BI=[1, n]$ is as follows. 
If $\ell\in (j_k, j_{k+1})$, we set 
$$\bar{J}^i=\left\{
\begin{array}{cl} 
 J^i  &\mbox{if $0\le i\le r$, $i\neq k$,}     \\
 J^k-\{\ell\} &\mbox{if $i=k$}.      \\
\end{array}
\right.$$
The corresponding open covering is 
$\{U_{\bar{J}^i}\}$ in $X_{[1, n]-\{\ell\}}$. 
Since $U_{\bar{J}^i}\subset U_{{J}^i}$ for each $i$, one has the 
restriction map 
$$\Z(X_{[1, n]}, \{U_{{J}^0}, \cdots, U_{{J}^r}\})
\to \Z( X_{[1, n]}, \{U_{\bar J^0}, \cdots, U_{\bar J^r}\})\,\,.
\eqno{(\ZXU.c)}$$

Let  $p: X_{[1, n]}\to  X_{[1, n]-\{\ell\}}$ be the projection;
$p$ is a projective map since $X_\ell$ is projective. 
The segmentation of $[1, n]-\{\ell\}$ corresponding 
its subset $\cJ$ is $\{\bar{J}^0, \cdots, \bar{J}^r\}$, so the corresponding 
open cover in $X_{[1, n]-\{\ell\}}$ is 
$$\cU(\cJ\subset [1, n]-\{\ell\}) =\{U_{\bar{J}^0}, \cdots, U_{\bar{J}^r}\}.$$
Since one has $p^{-1}(U_I)=U_I$ 
(more precisely, $p^{-1}(U_{I\subset[1, n]-\{\ell\} })=U_{I\subset [1, n]}$)
for a subset $I\subset [1, n]-\{\ell\}$, 
there is the push-forward map
$$p_*: \Z(X_{[1, n]}, \{U_{\bar J^0}, \cdots, U_{\bar J^r}\})
\to \Z(X_{[1, n]-\{\ell\} }, \{U_{\bar J^0}, \cdots, U_{ \bar J^r}\})\,.
\eqno{(\ZXU.d)}
$$
We define the map (\ZXU.b) as the 
composition of (\ZXU.c) and (\ZXU.d).

More generally for a subset $K\subset \ctop{\BI}-\J$ one has 
the corresponding projection 
$$\pi_K: 
\C(X_{\BI}, \cU(\J\subset \BI))\to
 \C( X_{\BI-K}, \cU(\J\subset \BI-K))\,\,.\eqno{(\ZXU.e)}$$
If $K=K'\amalg K''$, then $\pi_K=\pi_{K''}\pi_{K'}$, namely the following 
diagram commutes.

\vspace*{0.5cm}
\hspace*{2cm}
\unitlength 0.1in
\begin{picture}( 33.0000,  7.8000)(  2.9000,-15.5000)
\put(2.9000,-10.0000){\makebox(0,0)[lb]{$\Z(X_\BI, \cU(\J))$}}%
\put(35.9000,-9.9000){\makebox(0,0)[lb]{$\Z(X_{\BI-K}, \cU(\J))$}}%
\put(20.0000,-16.8000){\makebox(0,0)[lb]{$\Z(X_{\BI-K'}, \cU(\J))$}}%
%
{\color[named]{Black}{%
\special{pn 8}%
\special{pa 1500 950}%
\special{pa 3240 950}%
\special{fp}%
\special{sh 1}%
\special{pa 3240 950}%
\special{pa 3174 930}%
\special{pa 3188 950}%
\special{pa 3174 970}%
\special{pa 3240 950}%
\special{fp}%
}}%
%
{\color[named]{Black}{%
\special{pn 8}%
\special{pa 1300 1040}%
\special{pa 2010 1500}%
\special{fp}%
\special{sh 1}%
\special{pa 2010 1500}%
\special{pa 1966 1448}%
\special{pa 1966 1472}%
\special{pa 1944 1482}%
\special{pa 2010 1500}%
\special{fp}%
}}%
%
{\color[named]{Black}{%
\special{pn 8}%
\special{pa 2550 1520}%
\special{pa 3490 1050}%
\special{fp}%
\special{sh 1}%
\special{pa 3490 1050}%
\special{pa 3422 1062}%
\special{pa 3442 1074}%
\special{pa 3440 1098}%
\special{pa 3490 1050}%
\special{fp}%
}}%
\put(19.3000,-9.0000){\makebox(0,0)[lb]{$\pi_K$}}%
\put(12.5000,-14.3000){\makebox(0,0)[lb]{$\pi_{K'}$}}%
\put(30.7000,-14.1000){\makebox(0,0)[lb]{$\pi_{K''}$}}%
\end{picture}%
\vspace{0.8cm}

(3) The quasi-isomorphism
$\iota: \Z(A_\BI)\to \Z(X_\BI, \cU(\J))$ is compatible with restriction maps 
and projections. It means, for projection, the commutativity of the 
following diagram:
$$\begin{array}{ccc}
\Z(X_\BI, \cU(\J))&\mapr{\pi_K} &\Z(X_{\BI-K}, \cU(\J)) \\
\mapu{\iota}& &\mapur{\iota}\\
\Z(A_\BI)&\mapr{\pi_K} &\Z(A_{\BI-K})\,\,.
\end{array}
$$
Here $\pi_K$ at the bottom is the map induced by the projection 
$A_\BI\to A_{\BI-K}$. 
\bigskip

\sss{\ZXUgeneral}{\it The complex $\Z(X^{\cJ}_{[1, n]}, \cU(\cJ))$.}\quad
We now generalize all in the previous subsection under the general 
assumption that $S$ is quasi-projective and $X_i\to S$ are projective. 
With slight modifications,  we will have a quasi-isomorphic complex to
$\cZ_s(X_1\times_S\times \cdots \times_S X_n)$, one for each $\cJ\subset
(1, n)$, and restriction and projection maps as above.

 Let $S\injto \bar S$ be a compactification, namely  an open immersion 
to a projective variety.   For each $X_i$ take a projective variety 
$\bar X_i$ with a projective  map $\bar X_i\to \bar S$ extending $p_i$. 
(We say $\bar X_i / \bar S$ is a compactification of $X_i/S$.)

Then one has   
$$X_{[1, n]}=\prod X_i\qquad \text{and}\qquad 
\bar{X}_{[1, n]}:=\prod \bar{X}_i\,\,.$$ 
To non-empty $I\subset [1, n]$, there corresponds a closed set $A_I\subset X_{[1, n]}$ 
and its complement $U_I$ as in (\ZXU); similarly we have the
closed set  $\bar{A}_I\subset \bar{X}_{[1, n]}$ given by 
$$\bar{A}_I\mapr{\sim}
\underset{i\in I\,\,\,}{{\prod}_{\bar S}}
\bar{X}_i \times \underset{i\not\in I}{\prod}\bar{X}_i 
$$
 and its 
complement $\bar{U}_I$. 

Given $\J\subset (1, n)$, we define a partial compactification by 
$$X_{[1, n]}^\J:= \prod_{i\in [1, n]}X'_i\quad\mbox{with}\quad 
X'_i= 
\begin{cases}
\bar X_i & \mbox{if $i\in (1, n)-\J$} \\
     X_i & \mbox{if $i\in\{1, n\}\cup \J$\,. }
\end{cases} 
$$
It is the compactification of $X_{[1, n]}$ in the factor
$X_i$ with $i\in (1, n)-\J$. 

For $I\subset [1, n]$, 
define the closed subset $A_I^\J\subset X_{[1, n]}^\J$
by the following diagram:
$$\begin{array}{ccc}
X_{[1, n]}^\J &\mapr{\sim} 
&\underset{i\in I}{\prod}X'_i \times \underset{i\not\in I}{\prod}X'_i \\
\mapu{}  &  &\mapu{} \\
A_I^\J &\mapr{\sim}
&\underset{i\in I\,\,\,}{\prod_{\bar{S}}}
X'_i \times \underset{i\not\in I}{\prod}X'_i \,\,.
\end{array}$$
Note that we have Cartesian squares
$$\begin{array}{ccccc}
 X_{[1, n]}   &\injto   &X^\cJ_{[1, n]} &\injto &\bar{X}_{[1, n]}  \\
\mapu{} &     &\mapu{} &&\mapur{}   \\
 A_{[1, n]}   &\injto   &A^\cJ_{[1, n]} &\injto &\,\,\bar{A}_{[1, n]} \,. 
 \end{array}
$$  
If $\cJ\subset \cJ'$, then 
one has an open immersion $X^{\cJ'}_{[1, n]} \injto X^\cJ_{[1, n]} $
and the following square is Cartesian:
$$\begin{array}{ccc}
 X^{\cJ'}_{[1, n]}   &\injto   &X^\cJ_{[1, n]}  \\
\mapu{} &     &\mapu{}    \\
 A^{\cJ'}_{[1, n]}   &\injto   &\,\, A^\cJ_{[1, n]} \,.
 \end{array}
$$  
If $\cJ=(1, n)$, then $X^\cJ_{[1, n]}=X_{[1, n]}$. 

If $I\subset I'$ then $A_I^\cJ\supset A_{I'}^\cJ$. 
For two subsets $I$ and $I'$ with non-empty intersection,  
$A_{I\cup I'}^\J= A_I^\J\cap A_{I'}^\J$. 
Further, if $I\supset (1, n)-\cJ$ and 
$I\cap \,(\{1, n\}\cup \J)\, \neq\emptyset$, then  $A_I^\J= A_I$.
Indeed the first condition implies that $\underset{i\not\in I}{\prod}X'_i
=\underset{i\not\in I}{\prod}X_i$, and the second condition implies 
$\underset{i\in I\,\,\,}{\prod_{\bar{S}}}
X'_i =
\underset{i\in I\,\,\,}{\prod_{{S}}}
X_i$.  In particular, $A_{[1, n]}^\J= A_{[1, n]}=X_1\times_S\ddd \times_S X_n$. 

Let $U_I^\J= X_{[1, n]}^\J- A_I^\J$. 
If $I\subset I'$ then $U_I^\J\subset U_{I'}^\J$; if $I\cap I' \neq \emptyset$, then
$U_{I\cup I'}^\J= U_I^\J\cup U_{I'}^\J$. Note $U_{[1, n]}^\J= U_{[1, n]}$. 
One has Cartesian squares 
$$\begin{array}{ccc}
X_{[1, n]}&\subset  &X_{[1, n]}^\J \\
\cup& &\cup  \\
U_I&\subset &U_I^\J 
\end{array}
$$
and 
$$\begin{array}{ccc}
X^{\cJ'}_{[1, n]}&\subset  &X_{[1, n]}^\J \\
\cup& &\cup  \\
U^{\cJ'}_I&\subset &U_I^\J 
\end{array}
$$
for $\cJ\subset \cJ'$.

The $\J$ specifies a covering $\cU(\J)$ of 
 $U_{[1, n]}^\J=U_{[1, n]}$ defined now as
$$\cU(\J)=\{U^{\cJ}_{J^0}, \cdots, U^{\cJ}_{J^r}\}\,.$$
So we have the complex 
$\C(X_{[1, n]}^\J, \cU(\J))$ and a quasi-isomorphism
$\cZ(X_1\times_S\times \cdots \times_S X_n)\to 
\C(X_{[1, n]}^\J, \cU(\J))$. 

As before
the same construction can be applied to a  subset $\BI$ of $[1, n]$, 
and a subset $\J\subset \ctop{\BI}$. 
  One has the product $X_{\BI}$ and its partial 
 compactification $X^\cJ_{\BI}$. To each subset $I\subset \BI$ there corresponds a closed set $A^\cJ_{I}$ and its complement 
 $U^\cJ_{I}$. 
 A subset $\cJ\subset \ctop{\BI}$ gives a covering 
 $\cU(\cJ)=\cU(\cJ\subset \BI)$ of $U_{\BI}$, and thus 
 the complex $\C(X_{\BI}^\J, \cU(\J))$
 quasi-isomorphic to $\Z(A_{\BI})$.
 
Parallel to (1)-(3) of the previous subsection, we have:
\smallskip 

(1) For $\J\subset \J'$, we have a map, called the  
restriction map
$$\C(X_{\BI}^\J, \cU(\J)) \to \C(X_{\BI}^{\J'}, \cU(\J'))\,\,.\eqno{(\ZXUgeneral.a)}$$
It is a quasi-isomorphism. 

To obtain this map, note that 
the covering $\cU(\J\subset \BI)=\{U^{\cJ}_{J^0}, \cdots, U^{\cJ}_{J^r}\}$ 
restricts to the covering 
$$\cU(\J\subset \BI)\cap X_{\BI}^{\cJ'}=\{U^{\cJ'}_{J^0}, \cdots, U^{\cJ'}_{J^r}\}$$
of $U_{\BI}$ in $X_{\BI}^{\cJ'}$. 
On the other hand,  recall that 
$\cU(\J'\subset \BI)=\{U^{\cJ'}_{{J'}^0}, \cdots, U^{\cJ'}_{{J'}^{r'}}\}$
where $\{{J'}^0, \cdots, {J'}^{r'}\}$ are the associated intervals to $\cJ'$. 
The map $\lambda: [0, r']\to [0,r]$ such that $J'_t\subset J_{\lambda (t)}$ for $t\in [0, r']$ gives a refinement 
$$\lambda: \cU(\J\subset \BI)\cap X_{\BI}^{\cJ'} \to \cU(\J'\subset \BI)\,.$$
The map (\ZXUgeneral.a) is defined to be the composition of the restriction map 
and the refinement map
$$\Z(X_{\BI}^\J, \cU(\J\subset \BI))\to \Z(X_{\BI}^{\J'}, \cU(\J) \cap X_{\BI}^{\cJ'}) 
\to  \Z(X_{\BI}^{\J'},\cU(\J'\subset \BI))\,.$$

(2) For $\ell\in \ctop{\BI}-\J$ 
$$\proj_{\ell}: \C(X_{\BI}^\J, \cU(\J\subset \BI))\to \C( X_{\BI-\{\ell\}}^\J, 
\cU(\J\subset \BI-\{\ell\}))\,\,.\eqno{(\ZXUgeneral.b)}$$
\noindent This is defined in the same way as before, since
 the projection $p: X_{\BI}^\J=
X_{\BI-\{\ell\}}^\J\times \bar{X}_\ell \to X_{\BI-\{\ell\}}^\J$
is projective. 

More generally for $K\subset \ctop{\BI}-\J$ one has the projection
$\pi_K: \C(X_{\BI}^\J, \cU(\J\subset \BI))\to
 \C( X_{\BI-K}^\J, \cU(\J\subset \BI-K))$. 
 If $K=K'\amalg K''$, then $\pi_K=\pi_{K''}\pi_{K'}$. 
 
(3) The quasi-isomorphism $\iota: \Z(A_\BI)\to \Z(X_\BI^\J, \cU(\J))$ is compatible with restrictions
and projections.  
\bigskip 

\sss{\cFIcJ} {\it The complex $\cF(\BIc, \cJ)$.}\quad 
For simplicity let
$$\cF([1, n], \cJ)=\Z(X_{[1, n]}^\J, \cU(\J))\,\,;\eqno{(\cFIcJ.a)}$$
this is a complex of free $\ZZ$-modules
with differential $\partial $, and
write $\cF([1, n], \cJ)^\bullet$ to specify degrees.
The same construction applies to any finite subset
(of cardinality $\ge 2$) $\BIc\subset [1, n]$ and 
a subset $\cJ\subset\ctop{I}$,
so that one has the complex $\cF(\BIc, \cJ)$.  

For $\cJ\subset \cJ'$, there is the
corresponding restriction map 
$$r_{\cJ, \cJ'}: \cF(I, \cJ)\to \cF(I, \cJ')\,.\eqno{(\cFIcJ.b)}$$
The $r_{\cJ, \cJ'}$ is transitive for inclusions $\cJ\subset \cJ'\subset 
\cJ''$. 
There is a quasi-isomorphism $\iota:
\Z(A_I)\to \cF(I, \cJ)$, compatible with the restrictions and projections. 

To a subset $K$ of $\ctop{I}$, there corresponds the projection
$$\pi_K: \cF(I, \cJ)\to \cF(I-K, \cJ)\,.\eqno{(\cFIcJ.b)}$$
If $K=K'\amalg K''$, then $\pi_K=\pi_{K''}\pi_{K'}$. 
\bigskip 

\sss{\cFIcJSigma}{\it The complex $\cF(I, \cJ|\Sigma)$.}\quad
Let $I$ be the subset of $[1, n]$, and $I_1, \cdots, I_c$ be a segmentation of 
$I$ corresponding to $\Sigma\subset \ctop{I}$. 
Assume given subsets $\cJ_j\subset \ctop{I_j}$ for $j=1, \cdots, c$.
There correspond a sequence of fiberings on $[1, c]$ consisting of 
$M_j=X_{I_j}^{\cJ_j}$, $Y_j=X_{t_j}$ with $t_j=\term{I_j}$, and the natural 
projections $M_j\to Y_j\gets M_{j+1}$.
The condition (*) in Definition (\Seqfiberings) is obviously satisfied.  
One has open coverings $\cU_j=\cU(\cJ_j)$ of open sets $U^{\cJ_j}_{I_j}$ of $X^{\cJ_j}_{I_j}$.  

Let $\cJ=\cJ_1\cup\{t_1\}\cup \cJ_2\cup
\cdots \cup\{t_{c-1}\}\cup \cJ_c$. 
Then we have
$$X^{\cJ_1}_{I_1}\diamond\cdots\diamond X^{\cJ_c}_{I_c}=X^\cJ_I\,.$$
Let $p_j: X^\cJ_I\to X^{\cJ_j}_{I_j}$ be the projection. 
Then one has (easily verified) 
$$p_1^{-1}\cU(\cJ_1)\amalg\cdots \amalg p_c^{-1}\cU(\cJ_c)=\cU(\cJ)$$
which is an open covering of $U_I$ in $X^\cJ_I$. 

According to (\htsZMU) one has a $c$-tuple subcomplex
$$\bhts \Z(M_j, \cU_j) \subset \bigotimes_{j\in [1, c]} \Z(M_j, \cU_j)$$
which is also written 
$$\cF(\BIc_1, \cJ_1)\hts\cF(\BIc_2, \cJ_2
)\hts\cdots\hts\cF(\BIc_c, \cJ_c)
\subset 
\cF(\BIc_1, \cJ_1)\ts\cF(\BIc_2, \cJ_2
)\ts\cdots\ts\cF(\BIc_c, \cJ_c)\,.
$$
Recall that the differential of the tensor product complex, also denoted $d$, is given by
$$d(\al_1\ts\cdots\ts\al_c) =\sum (-1)^{\sum_{j>i} \deg \al_j}\,\al_1\ts\cdots \ts\al_{i-1}\ts 
d(\al_i)\ts\cdots \ts\al_c\,\,. $$
For convenience we set
$$\cF(I, \cJ\tbar \Sigma)=\cF(\BIc_1, \cJ_1)\ts\cF(\BIc_2, \cJ_2
)\ts\cdots\ts\cF(\BIc_c, \cJ_c)\,\,,\eqno{(\cFIcJSigma.a)}$$
if $\cJ\subset\ctop{I}-\Sigma$, $I_1, \cdots, I_c$ is the segmentation of 
$I$ by $\Sigma$, and $\cJ_i=\cJ\cap \ctop{I_i}$. 
Similarly we have 
$$\cF(I, \cJ\dbar \Sigma):=\cF(\BIc_1, \cJ_1)\hts\cF(\BIc_2, \cJ_2
)\hts\cdots\hts\cF(\BIc_c, \cJ_c)\,\,.\eqno{(\cFIcJSigma.b)} $$
The following properties hold by (\prophtsZMU). 
\smallskip

(1) There is a quasi-isomorphic 
inclusion $\iota_{\Sigma}:\cF(I, \cJ\dbar \Sigma)\to \cF(I, \cJ\tbar \Sigma)$. 
More generally, if $\Sigma\supset T$, $T$ gives the segmentation 
$I_1, \cdots, I_c$ of $I$, and $\cJ_j=\cJ\cap \ctop{I_j}$, 
$\Sigma_j=\Sigma\cap I_j$, 
then 
one has a quasi-isomorphic inclusion of $c$-tuple complexes
$$\iota_{\Sigma/T}: \cF(I, \cJ|\Sigma)
\injto \cF(I_1, \cJ_1|\Sigma_1)\ts\cdots\ts\cF(I_c, \cJ_c|\Sigma_c)
\eqno{(\cFIcJSigma.c)}$$
If $T=\emptyset$, $\iota_{\Sigma/\emptyset}$ is the identity, and if $T=\Sigma$,
$\iota_{\Sigma/\Sigma}$ coincides with $\iota_\Sigma$. 

(2) For $\cJ\subset\cJ'$ there is the corresponding quasi-isomorphism of complexes
$$r_{\cJ, \cJ'}: \cF(I, \cJ\dbar \Sigma)\to \cF(I, \cJ'\dbar \Sigma)\,; 
\eqno{(\cFIcJSigma.d)}$$
it is transitive in $\cJ$.  The map $r_{\cJ, \cJ'}$ is also compatible with 
the inclusion $\iota_{\Sigma/T}$ above, namely the 
following square commutes:
$$\begin{array}{ccc}
\cF(I, \cJ|\Sigma)&\injto  &
\cF(I_1, \cJ_1|\Sigma_1)\ts\cdots\ts\cF(I_c, \cJ_c|\Sigma_c) \\
\mapd{r_{\cJ, \cJ'}}& &\mapdr{} \\
\cF(I, \cJ'|\Sigma)&\injto  &\phantom{\,\,.}
\cF(I_1, \cJ'_1|\Sigma_1)\ts\cdots\ts\cF(I_c, \cJ'_c|\Sigma_c)\,\,
\end{array}
$$
where the right vertical map is $\ts r_{\cJ_i, \cJ'_i} $.
This compatibility will be referred to as the compatibility of 
the map $r$ and tensor product.

(3) 
For each $t=1, \cdots, c-1$ (which corresponds to $k=\term(I_t)\in \Sigma$),
one has the product map 
$$\begin{array}{c}
\rho_k: \cF(I_1, \cJ_1)\hts\cdots\hts\cF(I_t, \cJ_t)\hts\cF(I_{t+1}, \cJ_{t+1})\hts\cdots\hts
\cF(I_c, \cJ_c)  \\
\qquad \qquad \to\cF(I_1, \cJ_1)\hts\cdots\hts\cF(I_t\cup I_{t+1}, \cJ_t\cup\{k\}\cup
\cJ_{t+1})
\hts\cdots\hts\cF(I_c, \cJ_c)\,\,,
\end{array} $$
or 
$$\rho_k: \cF(I, \cJ |\Sigma)\to \cF(I, \cJ\cup\{k\} |\Sigma-\{k\})\,.\eqno{(\cFIcJSigma.e)} $$
This sends $u_1\ts\cdots\ts u_c$ to $u_1\ts\cdots\ts (u_t\scirc u_{t+1})\ts \cdots\ts u_c$. 
In the following 
 subsection we will discuss how to  specify dimensions of cycle complexes in a manner compatible with the product map. 

For distinct elements $k, k'\in \Sigma$, 
the following diagram commutes:
$$\begin{array}{ccc}
\cF(\BIc,\cJ \dbar \Sigma)&\mapr{\rho_k} &\cF(\BIc, \cJ \cup\{k\}\dbar \Sigma-\{k\}) \\
\mapd{\rho_{k'}}& &\mapdr{\rho_{k'}} \\
\cF(\BIc, \cJ\cup\{k'\}\dbar \Sigma-\{k'\})&\mapr{\rho_{k}} &\cF(\BIc, 
\cJ\cup\{k, k'\} \dbar \Sigma-\{k, 
k'\})\,\,.
\end{array}
$$
For $K \subset \Sigma$, let 
$$\rho_K: \cF(I, \cJ|\Sigma)\to \cF(I, \cJ\cup K|\Sigma-K)\,\,.
\eqno{(\cFIcJSigma.f)}$$
be the composition of $\rho_k$ for $k\in K$ in any order. 
If $K=K'\amalg K''$ then $\rho_K=\rho_{K''}\rho_{K'}$.
The map $\rho_K$ is compatible with the inclusion
 $\iota_{\Sigma/T}$, 
namely the following diagram commutes, where 
$K_i=K\cap \Sigma_i$:
$$\begin{array}{ccc}
\cF(I, \cJ|\Sigma)&\injto &\otimes \cF(I_i, \cJ_i|\Sigma_i)  \\
\mapd{\rho_K}& &\mapdr{\ts\rho_{K_i}} \\
\cF(I, \cJ\cup K|\Sigma-K )&\injto &\otimes \cF(I_i, \cJ_i\cup K_i|\Sigma_i-K_i)\,.
\end{array}
$$
 
(4) To $K\subset \ctop{I}-\Sigma$ disjoint from $\cJ$, 
there corresponds a map of  complexes
$$\pi_K: \cF(I, \cJ|\Sigma)\to \cF(I-K, \cJ|\Sigma)\,\,.
\eqno{(\cFIcJSigma.g)}$$
If $K=K'\amalg K''$ then $\pi_K=\pi_{K''}\pi_{K'}$. 
The map $\pi_K$ is compatible with the 
inclusion $\iota_{\Sigma/T}$, meaning the commutativity of a diagram as the 
one for $\rho_K$. 

(5) The maps $r$, $\rho$, and $\pi$ commute with each other. 
The commutativity of $r$ and $\rho$ means the commutativity of the 
following square:
$$\begin{array}{ccc}
\cF(\BIc,\cJ\dbar\Sigma)&\mapr{r_{\cJ, {\cJ}'}
} &\cF(\BIc, \cJ'\dbar\Sigma ) \\
\mapd{\rho_K}& &\mapdr{\rho_K} \\
\cF(\BIc, \cJ\cup K\dbar\Sigma-K )&\mapr{r_{\cJ\cup K, {\cJ}'\cup K}
} &\cF(\BIc-K, \cJ'\cup K\dbar\Sigma-K)\,.
\end{array}
$$
The reader can write down the commutative diagrams 
expressing the commutativity of $r$ and $\pi$, and of $\rho$ and $\pi$. 

(6) The maps $r$, $\rho$ and $\pi$ provide  another map
in the derived category. 
Let $K\subset \Sigma$. We have the maps 
$$\begin{array}{cccc}
 & &  &\cF(I, \emptyset|\Sigma)  \\
 &  & &\mapdr{\rho_K}  \\
&\cF(I, \emptyset|\Sigma-K) &\mapr{r_{\emptyset, K}} &\cF(I, K|\Sigma-K) \\
&\mapd{\pi_K} & & \\
&\cF(I-K, \emptyset|\Sigma-K) &  &
\end{array}
$$
Since the map $r$ is a quasi-isomorphism, inverting it gives a map 
in the derived category of abelian groups
$$\psi_K:\cF(I, \emptyset \tbar \Sigma)\to \cF(I-K, \emptyset\tbar
\Sigma-K)\,\,.\eqno{(\cFIcJSigma.h)}$$
We call this the {\it composition map}. 

The map $\psi_K$ satisfies transitivity in $K$,
which says $\psi_K=\psi_{K''}\psi_{K'}$ if $K=K'\amalg K''$.
This follows from the compatibility of $r$, $\rho$, and $\pi$. 

Also $\psi$ is compatible with tensor product: 
Assume $m\in \Sigma$, $m\not\in K$, and 
 $K$ is a subset of $\Sigma$.
Let $I'$, $I''$ be the segmentation of $I$ by $m$, $\Sigma'=\Sigma\cap I'$, $\Sigma''=\Sigma\cap I''$, and 
$K'=K\cap I'$, $K''=K\cap I''$. Then the following diagram commutes:
$$\begin{array}{ccc}
\cF(I, \emptyset\tbar \Sigma)&= &\cF(I', \emptyset\tbar \Sigma')\ts \cF(I'', \emptyset\tbar \Sigma'') \\
\mapd{\psi_K}& &\mapdr{\psi_{K'}\ts \psi_{K''} }\\
\cF(I-K, \emptyset\tbar \Sigma-K)&= &\cF(I'-K', \emptyset\tbar \Sigma'-K')
\ts \cF(I''-K'', \emptyset\tbar \Sigma''-K'') \,\,.
\end{array}
$$
This follows from the compatibility of $\iota_{\Sigma/T}$ with 
$r$, $\rho$, and $\pi$.
\bigskip 

\sss{\dimcyclcpx}{\it Dimensions of the cycle complexes.}\quad
The dimensions of the cycle complexes can be specified as follows. 
Assume to each interval $[i, i+1]\subset [1, n]$, an integer $a_i\in \ZZ$ 
is assigned. For a subset $I\subset [1, n]$, if $j=\init(I)$, $k=\term(I)$, 
let 
$$a_I=\sum_{j\le i\le k-1} a_i -\sum_{j\le i\le k-1} \dim X_i\,\,.$$
We then have the following property: 
For $I=[i, i+1]$, one has $a_I=a_i$. 
The number $a_I$ depends only on $\init(I)$ and $\term(I)$. 
If 
$\term(I)=\init(I')=c$, then 
$a_{I\cup I'}=a_I+a_{I'}-\dim X_c$. 

Correspondingly let $a_I$ be the dimension of the cycle complex for 
$\cF(I, \cJ):=\Z_{a_I}(X_{I}^\cJ, \cU(\J))$. 
The dimensions for $\cF(I, \cJ|\Sigma)$ are also specified by additivity. 

Then the maps
$\rho_k: \cF(I, \cJ |\Sigma)\to \cF(I, \cJ\cup\{k\} |\Sigma-\{k\})$
and  $\pi_K: \cF(I, \cJ|\Sigma)\to \cF(I-K, \cJ|\Sigma)$ 
are compatible with the dimensions thus specified. 
\bigskip 

\sss{\cFI}{\it The complex $\cF(I)$.}\quad
Let $I$ be a subset of $[1, n]$. 
We will define a complex denoted $\cF(I)$. 

There is the restriction map, for $\J\subset \J'$ with $|\cJ'|=|\cJ|+1$, 
$$r_{\J, \J'}: \cF(I, \cJ)\to \cF(I, \cJ')
\,\,.$$
This is a quasi-isomorphism.
Let 
$$A^{a, p}=\bop_{a=|\cJ|+1} \cF(I,\cJ)^p\,, $$
the sum over $\cJ$ with $a=|\cJ|+1$. 
One has differential $\partial: A^{a, p}\to A^{a, p+1}$. 
 If 
 $\J'=\J\cup\{k\}$,  let 
$\cJ_{>k}=\{i\in \cJ| i>k\}$, and $| \cJ_{>k}|$ be its cardinality. 
We define a map $r: A^{a, p}\to A^{a+1, p}$ as the sum of the maps
$$(-1)^{| \cJ_{>k}|}r_{\cJ, \cJ'}:\cF(I, \cJ)\to \cF(I, \cJ')$$
for pairs $\cJ, \cJ'$ with $\cJ\subset \cJ'$, 
$a=|\cJ|+1$, $a+1=|\cJ'|+1$. 
Then one has $rr=0$, so that 
 $A^{a, p}$ forms a ``double" complex with differentials $r, \partial $. 
The total complex $\Tot (A)$ is the complex with 
 differential $d_\cF$ (also denoted $d$ as long as there is no confusion), 
 which is  $r+(-1)^a \partial$ on $A^{a, p}$.
Now define
$$\cF(I)=(\Tot(A); d)\,\,. \eqno{(\cFI.a)}$$
It has complexes 
$\cF(I, \cJ)[-(|\cJ|+1)]$ as subquotients.
Here recall for a complex $(K^\bullet, d)$, the shift $K^\bullet[1]$ is 
defined by $(K^\bullet[1])^p=K^{p+1}$ and $d_{K[1]}=-d$.  
\bigskip 

(\propcFI){\bf Proposition.}\quad 
(1) 
If $|\BIc|=2$, $\cF(\BIc)=\cF(\BIc,\emptyset)[-1]$, 
so there is a quasi-isomorphism
$\Z(A_{\BIc})[-1]\to \cF(\BIc)$. 

(2) If $|\BIc|\ge 3$, then $\cF(\BIc)$ is acyclic. 
\smallskip 

{\it Proof.}\quad 
(1) Obvious. 

(2) Follows from the following lemma (apply to $T=\ctop{I}$ and $\cJ\mapsto 
\cF(I, \cJ)$). 
\bigskip 

\sss{\lemmaFI}\quad Let $T$ be a non-empty finite  
set.  Suppose to each subset $S$ of $T$ (including the empty set) there
corresponds a complex $C_S\in C(Ab)$, 
 and to each inclusion $S\subset S'$
there corresponds a map of complexes $f_{S\, S'}: C_S\to C_{S'}$, satisfying 
$f_{S\, S}=id$, and 
 $f_{S'\, S''}f_{S\, S'}=f_{S\, S''}$ for $S\subset S'\subset S''$.
 
We then have a ``double" complex  
$$0\to C_\emptyset \to \bop_{|S|=1}C_S  \to \bop_{|S|=2}C_S
\to \cdots \to C_T\to 0\,\,.$$
Here the maps are  signed sums of the maps $f_{S\, S'}$ with $|S'|=|S|+1$, 
the signs being specified as in (\cFI).  
We can form its  total complex
$\Tot(C_S)$.  One proves 
the following lemma by induction on $n$. 
\bigskip 

(\lemmaFI.1)\Lem {\it  Assume for each $S\subset S'$ the map $f_{S\, S'}$
 is a quasi-isomorphism. Then $\Tot(C_S)$ is acyclic. 
}\bigskip 

\sss{\cFISigma}{\it The complex $\cF(I|\Sigma)$.}\quad 
Let $I$ be the subset of $[1, n]$, and $I_1, \cdots, I_c$ be a segmentation of 
$I$ corresponding to $\Sigma\subset \ctop{I}$. 
We have the $c$-tuple complex
$$\cF(\BIc \tbar \Sigma):=\cF(\BIc_1)\ts\cF(\BIc_2)\ts\cdots\ts
\cF(\BIc_c)\,\,,\eqno{(\cFISigma.a)} $$
which is the tensor product of the simple complexes $\cF(I_i)$. 
One has  as a module
$$\cF(\BIc \tbar \Sigma)=\oplus  \cF(\BIc, \cJ \tbar \Sigma)\,,\eqno{(\cFISigma.b)} $$
the sum for subsets $\cJ$ of $\ctop{I}-\Sigma$;
the differential is a signed sum of the maps 
$$1\ts\cdots\ts 1\ts \partial \ts 1\ts\cdots \ts 1:\cF(\BIc, \cJ \tbar \Sigma)\to \cF(\BIc, 
\cJ \tbar \Sigma)
$$
where $\partial $ is the differential of $\cF(I_i, \cJ_i)$, and the maps
$$1\ts\cdots\ts 1\ts r_{\cJ_i, \cJ'_i}\ts 1\ts\cdots \ts 1:\cF(\BIc, \cJ \tbar \Sigma)\to \cF(\BIc, \cJ' \tbar \Sigma)\,,
$$
for $\cJ\subset \cJ'$ with $|\cJ'|=|\cJ|+1$ (then there is $i$ with 
$|\cJ_i'|=|\cJ_i|+1$). 
Note there is a filtration of $\cF(\BIc \tbar \Sigma)$
 by subcomplexes whose successive quotients are direct sums of 
the $c$-tuple complexes $\cF(\BIc, \cJ \tbar \Sigma)$ . 
Similarly we have another $c$-tuple complex $\cF(\BIc | \Sigma)$
which is defined, as the sum of modules, by  
$$\cF(\BIc | \Sigma)=\oplus  \cF(\BIc, \cJ | \Sigma)\,.\eqno{(\cFISigma.c)} $$
It is a subcomplex of $\cF(I\tbar \Sigma)$, since the map 
 $1\ts\cdots\ts \partial \ts\cdots \ts 1$ takes  $\cF(\BIc, \cJ | \Sigma)$ to itself and 
the $1\ts\cdots\ts r_{\cJ_i, \cJ'_i}\ts\cdots \ts 1$ takes 
$\cF(\BIc, \cJ | \Sigma)$ to $\cF(\BIc, \cJ' | \Sigma)$. 

There is a filtration by subcomplexes whose successive quotients are direct sums of 
the $c$-tuple complexes $\cF(\BIc, \cJ | \Sigma)$. 
We also use the notation 
$\cF(\BIc_1)\hts\cF(\BIc_2)\hts\cdots\hts\cF(\BIc_c)$ for
$\cF(\BIc | \Sigma)$. \bigskip 

(1) The inclusion $\iota: \cF(\BIc | \Sigma)\to \cF(\BIc \tbar \Sigma)$
is a quasi-isomorphism. 
More generally, if  $\Sigma\supset T$ and $T$ gives the segmentation 
$I_1, \cdots, I_c$ of $I$, then 
one has a quasi-isomorphic inclusion of $c$-tuple complexes
$$\iota_{\Sigma/T}: \cF(I|\Sigma)
\injto \cF(I_1|\Sigma_1)\ts\cdots\ts\cF(I_c|\Sigma_c)
\eqno{(\cFISigma.d)}$$

(2) Recall that $\cF(I)$ was defined to be the total complex of a ``double" complex.
Likewise there is a ``double"
complex  whose total complex is canonically
isomorphic to  $\cF(I\dbar \Sigma)$. 
We explain this in the case $|\Sigma|=1$. 

Let $A^{\dbullet}$, $B^{\dbullet}$
 be the ``double" complexes as above defining 
$\cF([1,m])$, $\cF([m, n])$, respectively.  
Then the tensor product of them, $A^{\dbullet}\ts B^{\dbullet}$, 
is a ``quadruple" complex. 
There is a quasi-isomorphic 
 multiple subcomplex in it, given by 
$$A^{a, p}\hts B^{b, q}=\bop_{a=|\cJ|+1, b=|\cJ'|+1} \cF([1, m],  \cJ)^p\hts 
 \cF([m, n],  \cJ')^q\,\,.$$

Recall from (\Tensorproductcpx) that  
$A^{\dbullet}\times 
B^{\dbullet}$ is the ``double" complex $E^{\dbullet}$
obtained by applying the totalization $\Tot_{13}\Tot_{24}$ to the 
``quadruple" complex $A^{\dbullet}\ts B^{\dbullet}$; specifically
$$E^{c, r}=\bop_{a+b=c\,,p+q=r}A^{a,p}\ts B^{b,q}\,\, $$
with appropriate differentials $r$, $\partial$. 
Likewise one defines a quasi-isomorphic ``double" complex
$$A^{\dbullet}\htimes 
B^{\dbullet}:= \Tot_{13}\Tot_{24}(A^{\dbullet}\hts B^{\dbullet})
\,.$$ 
By (\Tensorproductcpx) 
there is  an isomorphism of complexes 
$$u: \cF([1,m])\hts \cF([m,n])=\Tot(A)\hts\Tot(B)\isoto \Tot(A^{\dbullet}\htimes 
B^{\dbullet})\,\,.$$ 
Although the groups on the two sides are identical, the differentials
are different, and $u$ is not the identity.

More generally, let $A_1^{\dbullet}, \cdots, A_c^{\dbullet}$
 be the ``double" complexes for $\cF(I_1), \cdots, \cF(I_c)$, respectively.
Then we have the ``$2c$-tuple" tensor product complex 
$A_1^{\dbullet}\otimes \cdots \otimes A_c^{\dbullet}$, 
a ``double'' complex 
$$A_1^{\dbullet}\times \cdots \times A_c^{\dbullet}
:= \Tot_{1, 3, \cdots, 2c-1}\Tot_{2, 4, \cdots, 2c}(A_1^{\dbullet}\ts
 \cdots \ts A_c^{\dbullet})\,$$
and  its quasi-isomorphic subcomplex
$$A_1^{\dbullet}\htimes \cdots \htimes A_c^{\dbullet}
:= \Tot_{1, 3, \cdots, 2c-1}\Tot_{2, 4, \cdots, 2c}(A_1^{\dbullet}\hts
 \cdots \hts A_c^{\dbullet})\,.$$
There is an isomorphism of complexes 
$$u: \cF(I|\Sigma)=\Tot(A_1^{\dbullet})\hts\cdots\hts\Tot(A_c^{\dbullet})
\isoto \Tot(A_1^{\dbullet}\htimes \cdots \htimes A_c^{\dbullet})\,\,.$$

(3) We next define a map of complexes
$\rho:\cF([1, m])\hts \cF([m, n])\to \cF([1, n])$, which we still call the product map. 
Let $A^{\dbullet}$, $B^{\dbullet}$, 
$C^{\dbullet}$ be the ``double" complexes as above defining 
$\cF([1,m])$, $\cF([m, n])$ and $\cF([1,n])$, respectively.  
The  product maps in the previous subsection
$$\rho:\cF([1,m], \J)\hts \cF([m,n], \J')\to \cF([1, n], \J\cup\{m\}\cup\J')$$
give a map of  ``double" complexes 
$$\rho:
A^{\dbullet}\htimes B^{\dbullet}\to C^{\dbullet}\,\,.$$
Namely we can verify that, for $u\in A^{a, p}$ and $v\in B^{b, q}$, 
$$\partial\rho(u\ts v)=\rho(\, (-1)^q \partial u\ts v + u\ts \partial v)\,,$$
with respect to the second differential $\partial$. 
Also we have
$$r\rho(u\ts v)=\rho(\, (-1)^b r(u)\ts v + u\ts r(v)\,)\,,$$
which follows immediately from the definition of $r$.
Taking $\Tot$ and using the isomorphism $u$ above we get a map of complexes
${\rho}: \cF([1,m])\hts \cF([m,n])\to \cF([1, n])$;
it sends $u\ts v$, $u\in A^{a, p}$, $v\in B^{b, q}$, to $(-1)^{aq} u\scirc v$. 

More generally if $I_1, \cdots, I_c$ is a segmentation of $I\subset [1, n]$, 
and $I_t\cap I_{t+1}=k$, there is the corresponding product map 
$$\rho_k: \cF(I_1)\hts\cdots\hts\cF(I_t)\hts\cF(I_{t+1})\hts\cdots\hts
\cF(I_c)\to \cF(I_1)\hts\cdots\hts\cF(I_t\cup I_{t+1})
\hts\cdots\hts
\cF(I_c)\,\,, $$
or $\rho_k: \cF(\BIc\dbar \Sigma)\to \cF(\BIc\dbar \Sigma-\{k\})$. 
This is defined just as above, changing  the order
of totalization and tensor product in
the factors $\cF(I_t)$, 
$\cF(I_{t+1})$, only.  Thus
$$\rho_k (\al_1\ts\cdots\ts\al_c)= (-1)^{aq}\al_1\ts\cdots\ts (\al_t\scirc \al_{t+1})\ts\cdots\ts\al_c\eqno{(\cFISigma.e)}$$
if $\al_t\in \cF(I_t, \cJ_t)^p$, $a=|\cJ_t|+1$, and $\al_{t+1}\in \cF(I_{t+1}, \cJ_{t+1})^q$, $b=|\cJ_{t+1}|+1$. 
The same notation $\rho_k$ is used for the map (\cFIcJSigma.e) and the above map, so we will avoid possible confusions 
by specifying which map we mean. 

The following diagram commutes (for distinct $k, k'\in\Sigma$):
$$\begin{array}{ccc}
\cF(\BIc\dbar \Sigma)&\mapr{\rho_k} &\cF(\BIc\dbar \Sigma-\{k\}) \\
\mapd{\rho_{k'}}& &\mapdr{\rho_{k'}} \\
\cF(\BIc\dbar \Sigma-\{k'\})&\mapr{\rho_{k}} &\cF(\BIc\dbar \Sigma-\{k, 
k'\})\,\,.
\end{array}
$$
For $K\subset \Sigma$ let 
$$\rho_K: \cF(I|\Sigma)\to \cF(I|\Sigma-K)\eqno{(\cFISigma.f)}$$
be the composition of $\rho_k$ for $k\in K$ in any order. 
The map $\rho_K$ is compatible with the inclusion $\iota_{\Sigma/T}$, 
namely the following diagram commutes, where 
$K_i=K\cap \Sigma_i$:
$$\begin{array}{ccc}
\cF(I|\Sigma)&\injto &\cF(I_1|\Sigma_1)\ts\cdots\ts \cF(I_c|\Sigma_c) \\
\mapd{\rho_K}& &\mapdr{\ts\rho_{K_i}} \\
\cF(I|\Sigma-K )&\injto &\cF(I_1-K_1|\Sigma_1)\ts\cdots\ts \cF(I_c-K_c|\Sigma_c)\,.
\end{array}
$$

(4) Recall from (\secfuncpx.4) that 
for a subset $K\subset \ctop{I}$ disjoint from $\cJ\cup\Sigma$, one has 
the map $\pi_K: \cF(I, \cJ|\Sigma)\to \cF(I-K, \cJ|\Sigma)$. 
Using this we will produce, for $K\subset \ctop{I}-\Sigma$, a 
map of complexes 
$$\pi_K: \cF(I|\Sigma)\to \cF(I-K|\Sigma)\,.\eqno{(\cFISigma.g)}$$
Let $\pi_K: \bop \cF(I, \cJ|\Sigma)\to \bop \cF(I-K, \cJ|\Sigma)$
be the sum of the maps $\pi_K: \cF(I, \cJ|\Sigma)\to \cF(I-K, \cJ|\Sigma)$ 
for $K$ disjoint from $\cJ$, and the zero maps on $\cF(I, \cJ|\Sigma)$ if 
$K\cap \cJ\neq \emptyset$. 
The repeated use of $\pi_K$ will not lead to a confusion. 
The map $\pi_K$ is compatible with the inclusion $\iota_{\Sigma/T}$.

(5) The maps $\rho_{K'}$ and $\pi_{K}$ (for $K$ and $K'$ disjoint)
commute with each other, namely the following square
commutes:
$$\begin{array}{ccc}
\cF(\BIc\dbar\Sigma)&\mapr{\rho_{K'} } &\cF(\BIc\dbar\Sigma-K') \\
\mapd{\pi_K}& &\mapdr{\pi_K} \\
\cF(\BIc-K\dbar\Sigma)&\mapr{\rho_{K'} } &\cF(\BIc-K-K'\dbar\Sigma-K')\,\,.
\end{array}
$$

(6) It follows from (\cFI.1) that $\cF(I|\Sigma)$ is acyclic unless $\Sigma=\ctop{I}$. 
If $I=[1,n]$ and $I_i=[i, i+1]$, 
\begin{eqnarray*}
\cF(I|\ctop{I})&=&\cF(I_1)\hts\cdots\hts\cF(I_{n-1}) \\
     &=&  \cF(I_1,\emptyset)[-1]\hts\cdots\hts\cF(I_{n-1},\emptyset)[-1]\,\,.
\end{eqnarray*}
So one has quasi-isomorphisms
$$\cF(I|\ctop{I})\injto
\cF(I_1)\hts\cdots\hts\cF(I_{n-1}) 
\hookleftarrow
\Z(A_{I_1})[-1]\ts\cdots\ts\Z(A_{I_{n-1}})[-1]\,\,.$$
\bigskip 
\sss{\barcpx}{\it Variant of the bar complex.}\quad  
We give a variant of the so-called bar complex. 
Let $n\ge 2$ and assume:
\smallskip

(a) To each sub-interval  $I\subset [1,n]$ of cardinality $\ge 2$, 
a complex of abelian groups
$ (A(I)^\bullet, d_A)$ is assigned.

(b) For a sub-interval $I$ and a 
segmentation of $I$ into $I'$ and $I''$, there corresponds a map 
of complexes $\rho: A(I')\ts A(I'')\to A(I)$. If $I$ is segmented
into three intervals $I', I'', I'''$, then the following commutes:
$$\begin{array}{ccc}
A(I')\ts A(I'')\ts A(I''')&\mapr{\rho\ts 1} &
 A(I'\cup I'')\ts A(I''')\\
\mapd{\rho}& &\mapdr{\rho} \\
A(I')\ts A(I''\cup I''')&\mapr{1\ts \rho} &A(I)\,\,.
\end{array}
$$
In the following we write $\al\cdot\be$ for $\rho(\al\ts\be)$. 
\smallskip 

For a segmentation $(I_1, \cdots ,I_c)$ of $[1, n]$, one forms 
 the tensor product of modules 
$A(I_1)\ts A(I_2)\ts\cdots\ts A(I_c)$. Let 
$$B([1, n])=\bop A(I_1)\ts A(I_2)\ts\cdots\ts A(I_c)\eqno{(\barcpx.a)}$$
be the sum over all segmentations. Give a grading on $B([1, n])$ by 
$$\overline{\deg} (\al_1\ts\cdots\ts\al_c)=\sum (\deg \al_i-1)$$
and give  differentials  by (put $\eps_j=\deg(\al_j)-1$ )
$$\bar{d}(\al_1\ts\cdots\ts\al_c) =-\sum (-1)^{\sum_{j>i} \eps_j}\,\al_1\ts\cdots \ts\al_{i-1}\ts 
d_A(\al_i)\ts\cdots \ts\al_c\,\,,$$
$$\bar{\rho}(\al_1\ts\cdots\ts\al_c)= \sum 
(-1)^{\sum_{j\ge i} \eps_j}\,\al_1\ts\cdots \ts\al_{i-2}\ts 
(\al_{i-1}\cdot\al_i)\ts\cdots \ts\al_c\,\,. $$
One verifies $\bar{d}\bar{d}=0$, $\bar{\rho}\bar{\rho}=0$ and
 $\bar{d}\bar{\rho}+\bar{\rho}\bar{d}=0$, so that 
$d_{B}=\bar{d}+\bar{\rho}$ is a differential. 
Thus with the grading $\overline{\deg}$ and the differential $d_B$, 
the module $B([1, n])$ is a complex; we call it the {\it bar complex} associated to 
$(A(I); d_A, \rho)$. 

To a subset $\Sigma\subset (1, n)$ there corresponds a segmentation
$I_1, \cdots, I_c$ of $[1, n]$. 
Let 
$$A([1, n]\tbar \Sigma)=A(I_1)\ts A(I_2)\ts\cdots\ts A(I_c)\,\,.$$
Then the above definition can be rewritten
$$B([1, n])=\bop_{\Sigma}A([1, n]\tbar \Sigma)$$
as a group. 
The differential $\bar{d}$ is the sum of 
the maps $\bar{d}:A([1, n]\tbar \Sigma)
\to A([1, n]\tbar \Sigma)$, so each $A([1, n]\tbar \Sigma)$ is a complex 
with respect to $\bar{d}$. 
In view of the definition of $\bar{d}$, one has an identity of complexes
$$A([1, n]\tbar \Sigma)=A(I_1)[1]\ts\cdots\ts A(I_c)[1]\,.$$
The differential $\bar{\rho}$ is the sum of the maps 
$$\bar{\rho}_k: A([1, n]\tbar \Sigma)\to A([1, n]\tbar \Sigma- \{k\})\,\,\quad 
\mbox{for $k\in \Sigma$}\,.$$

The complex $B([1, n])$ is of the form 
$$\cdots \mapr{\bar{\rho}}\bop_{|\Sigma|=2}A([1, n]\tbar \Sigma)
\mapr{\bar{\rho}}\bop_{|\Sigma|=1}A([1, n]\tbar \Sigma)
\mapr{\bar{\rho}} A([1, n])\to 0\,.$$
There is a decreasing filtration $F^k$ of $B([1, n])$ by subcomplexes given by
$$F^kB([1, n])=\bop_{|\Sigma|\ge k}A([1, n]\tbar \Sigma)\,, \qquad k\ge 0\,.$$
A graded quotient in the filtration $\Gr^k_F$ is the direct sum of the complexes 
$A([1, n]\tbar \Sigma)$ for $\Sigma$ with $|\Sigma|=k$. 
\bigskip

(1) There are quotient complexes of  $B([1, n])$ defined as follows. 
For a subset $S\subset (1, n)$, 
$\bop_{\Sigma\not\supset S} A([1,n]\tbar \Sigma)$ is 
a subcomplex of $B([1, n])$; the quotient complex is denoted  
$B([1, n]\tbar S)$:
$$B([1, n]\tbar S)= \bop_{\Sigma\supset S} 
A([1,n]\tbar \Sigma)\,\,\qquad\mbox{as a module.}\eqno{(\barcpx.b)}$$
For $S=\emptyset$, $B([1, n]\tbar \emptyset)=B([1, n])$. 
For $S=(1, n)$, one has $B([1, n]\tbar (1, n))=A([1,n]\tbar (1, n)\,)$. 
If $S\subset S'$ there is a natural surjection of complexes 
$$\tau_{S\, S'}: B([1, n]\tbar S)\to B([1, n]\tbar S')\,\,.\eqno{(\barcpx.c)}$$
If $S\subset S'\subset S''$ then $\tau_{S\,S''}=\tau_{S'\,S''}
\tau_{S\,S'}$.

Note the above construction can be applied to any sub-interval  $I\subset [1, n]$,
with $|I|\ge 2$, in place of $[1, n]$. 
So one has the complex $B(I)$, 
the quotient complexs  $B(I\tbar S)$ for $S\subset \ctop{I}$,
and the maps $\tau_{S\, S'}: B(I\tbar S)\to B(I\tbar S')$. 

(2) If $S$ corresponds to a segmentation $I_1, \cdots, I_c$ of $I=[1, n]$, there 
is an identity of complexes 
$$B(I \tbar S)=B(I_1)\ts\cdots\ts B(I_c)\,\,,\eqno{(\barcpx.d)}$$
where the right hand side is the usual tensor product of complexes. 
This is obtained by taking the sum, for $\Sigma\supset S$, of the identities of modules 
$$A(I\tbar\Sigma)=A(I_1\tbar\Sigma_1)\ts\cdots\ts A(I_c\tbar\Sigma_c)$$
 where $\Sigma_i=\Sigma\cap I_i$. 
One verifies immediately that (\barcpx.d) is an isomorphism of complexes.
\bigskip 

\sss{\cpxFI} {\it The complex $F(I)$.}\quad  
The complexes $I\mapsto \cF(I)$ in (\cFI) almost satisfy the assumptions 
(a) and (b) at the beginning of the previous subsection; instead of a product map  
$\rho$ we have only a partially defined product map. 
But one can define the bar complex in a parallel manner, as we explain now. 

Let $I=[1, n]$ (or more generally a subset of $[1, n]$), and let 
$$F(I)=F([1, n])=\bop_{\Sigma}\cF([1, n]| \Sigma)\eqno{(\cpxFI.a)}$$
 as a group. 
An element of $\cF(I|\Sigma)$ of multi-degree $(\delta_1, \cdots,  \delta_c)$ is 
a sum of elements $\al=\al_1 \ts\cdots\ts \al_c$ with $\al_i\in \cF(I_i)$, $\deg (\al_i)=
\delta_i$. 
For such $\al$, let $\eps_j=\delta_j-1$ and 
$\overline{\deg} (\al_1\ts\cdots\ts\al_c):=\sum \eps_j$. 
Define a map  $\bar{d}:\cF([1, n]|\Sigma)\to \cF([1, n]|\Sigma)$ by 
$$\bar{d}(\al_1 \ts\cdots\ts \al_c) =- \sum_i (-1)^{\sum_{j>i} \eps_j}\,
\al_1\ts\cdots \ts\al_{i-1}\ts d(\al_i)\ts\cdots \ts\al_c \,\,.\eqno{(\cpxFI.b)} $$
(Since $\cF([1, n]|\Sigma)$ is a multiple complex,  
the $i$-th differential of $\al$, which equals 
$ \al_1\ts\cdots \ts\al_{i-1}\ts d(\al_i)\ts\cdots \ts\al_c$, is in $\cF([1, n]|\Sigma)$. ) 
One has $\bar{d}\bar{d}=0$, so $\cF([1, n]|\Sigma)$ is a complex with respect
to $\overline{\deg}$ and $\bar{d}$. 
Note $(\cF([1, n]|\Sigma), \bar{d})$ differs from the complex $\cF([1, n]|\Sigma)$ 
in (\cFISigma) by shift, namely 
$$(\cF([1, n]|\Sigma), \bar{d})=\cF([1, n]|\Sigma)[c]$$
 where $c=|\Sigma|+1$. 
We will often write $\cF([1, n]|\Sigma)^{\shift}$ for the former to avoid
confusion. 
Let  $\bar{d}:F([1,n])\to F([1, n])$ be the sum of them. 
For $k\in \Sigma$, define a map $\bar{\rho}_k: 
\cF([1, n]|\Sigma)\to \cF([1, n]|\Sigma-\{k\})$ by
$$\bar{\rho}_k(\al_1 \ts\cdots\ts \al_c)=  
(-1)^{\sum_{j\ge i} \eps_j} \rho_{k}(\al_1 \ts\cdots\ts \al_c) \,\,;\eqno{(\cpxFI.c)}$$
where $k=\term(I_{i-1})$ and $\rho_k$ is the map in (\cFISigma.e).
Let $\bar{\rho}:F([1,n])\to F([1, n])$ be the sum of $\bar{\rho}_k$. 
One has $\bar{d}\bar{d}=0$, $\bar{\rho}\bar{\rho}=0$ and
 $\bar{d}\bar{\rho}+\bar{\rho}\bar{d}=0$, so that 
$$d_{F}:=\bar{d}+\bar{\rho}\eqno{(\cpxFI.d)}$$
 is a differential on 
$F([1, n])$. 

There is a decreasing filtration $F^k$ of 
$F([1, n])$ by subcomplexes given by
$$F^k F([1, n])=\bop_{|\Sigma|\ge k}\cF([1, n] | \Sigma)\,, \qquad k\ge 0\,.$$
A graded quotient in the filtration $\Gr^k_{F}$ is the direct sum of the complexes 
$\cF([1, n]| \Sigma)^{\shift}$ for $\Sigma$ with $|\Sigma|=k$. 
Since $\cF([1, n]| \Sigma)$ is acyclic unless $\Sigma=\ctop{I}$, the 
quotient map $F([1, n])\to \Gr^{n-2}_{F}=\cF([1, n]| (1, n))^{\shift}$
is a quasi-isomorphism. 
\smallskip

(1) For $S\subset (1, n)$ one has the corresponding quotient 
complex $F([1, n]|S)$ such that 
$$F([1, n]| S)= \bop_{\Sigma\supset S} 
\cF([1,n]\tbar \Sigma)\eqno{(\cpxFI.e)}$$
as a module. 
For $S=\emptyset$, $F([1, n]|\emptyset)=F([1, n])$. 
For $S=(1, n)$, one has $F([1, n]| (1, n))=\cF([1,n]| (1, n)\,)$. 
If $S\subset S'$ there is a natural surjection of complexes 
$$\sigma_{S\, S'}: F([1, n]| S)\to F([1, n]| S')\,\,.\eqno{(\cpxFI.f)}$$
If $S\subset S'\subset S''$ then $\sigma_{S\,S''}=\sigma_{S'\,S''}
\sigma_{S\,S'}$. 

The constructions so far apply  to any subset  $I\subset [1, n]$,
in place of $[1, n]$; one has complexes $F(I|S)$ and maps 
$\sigma_{S\,,  S'}$.

The maps $\sigma_{S\, S'}:
F(I|S)\to F(I|S')$ are quasi-isomorphisms. Indeed by 
considering a filtration $F^k$ as that of $F(I)$, one shows that
 the quotient map $F(I|S)\to \cF(I|\ctop{I})^{\shift}$ is 
a quasi-isomorphism. 

(2) If $S\supset T$,
$I_1, \cdots, I_c$ is the corresponding segmentation of $I$, 
and $S_i=S\cap I_i$, there is a quasi-isomorphic
inclusion of $c$-tuple complexes 
$$\iota_{S/T}: F(I|S)\injto F(I_1|S_1)\ts\cdots\ts F(I_c|S_c)\,.
\eqno{(\cpxFI.g)}$$
This is defined as the sum, for $\Sigma\supset S$,  of the inclusions 
$\iota_{\Sigma/T}: \cF(I|\Sigma)
\injto \cF(I_1|\Sigma_1)\ts\cdots\ts\cF(I_c|\Sigma_c)$.
If $T=\emptyset$, then $\iota_{S/\emptyset}=id$. 
In the particular case $T=S$, one has a quasi-isomorphic inclusion
$\iota_S:  F(I|S)\injto F(I\tbar S)$. 

The composition of the map 
$\iota_{S/T}$
followed by
$$\ts\iota_{S_i}: F(I_1|S_1)\ts\cdots\ts F(I_c|S_c)\to F(I_1\tbar S_1)\ts\cdots\ts F(I_c\tbar S_c)=F(I\tbar S)$$
coincides with $\iota_S$. 

The maps $\iota_{S/T}$ and $\sigma_{S\, S'}$ are compatible, 
namely 
if $S\subset S'$ and $S'_i=S' \cap I_i $, 
the following commutes:
$$\begin{array}{ccc}
F(I|S)&\mapr{\iota_{S/T}} & F(I_1|S_1)\ts\cdots\ts F(I_c|S_c) \\
\mapd{\sigma_{S\, S'}}& &\mapdr{\ts \sigma_{S_i\, S'_i} } \\
F(I|S')&\mapr{\iota_{S'/T}} & F(I_1|S'_1)\ts\cdots\ts F(I_c|S'_c)\,\,.
\end{array}
$$

(3) For $S\subset S'$, also define a map 
$$\tau_{S\, S'}: F(I\tbar S)\to F(I\tbar S')\eqno{(\cpxFI.h)}$$
 as follows. 
First define $\tau_S=\tau_{\emptyset\, S}: F(I)\to F(I\tbar S)$ 
as the composition 
$$F(I)\mapr{\sigma_{\emptyset, S} } F(I|S)\mapr{\iota_S} F(I\tbar S)\,.$$
For $S\subset S'$, let $I_1, \cdots, I_c$ be the 
 segmentation of $I$ by $S$, $S'_i=S'\cap I_i$, and 
$$\tau_{S\,S'}:=\ts \tau_{S'_i}: F(I_1)\ts\cdots\ts F(I_c)
\to F(I_1\tbar S'_1)\ts\cdots\ts F(I_c\tbar S'_c)=F(I\tbar S')\,.$$

One verifies the following two 
properties (using the properties stated in (2)).  
If $S\subset S'\subset S''$ then $\tau_{S\,S''}=\tau_{S'\,S''}
\tau_{S\,S'}$. 
The maps $\tau_{S\, S'}$ and $\sigma_{S\, S'}$ are compatible, namely
 the following diagram commutes:
$$\begin{array}{ccc}
F(I| S)&\mapr{\iota_S} &F(I \tbar S) \\
\mapd{\sigma_{S\, S'}}& &\mapdr{\tau_{S\, S'}} \\
F(I|S')&\mapr{\iota_{S'}} &F(I\tbar S')\,\,.
\end{array}
$$

(4)
For $K\subset \ctop{I}$ disjoint from $S$, one has a
map of complexes
$$\vphi_K: F(I\dbar S)\to F(I-K\dbar S)$$ 
given as follows. 
If $K$ is disjoint from $\Sigma$, the map 
$\pi_K: \cF(I|\Sigma)\to \cF(I-K|\Sigma)$ commutes with 
$\rho_k$ ($k\in \Sigma)$, thus also with $\bar{\rho}_k$;
likewise $\pi_K$ commutes with $\bar{d}$. 
It follows that if we define $\vphi_K$ to be the sum of 
$\pi_K:\cF(I|\Sigma)\to \cF(I-K|\Sigma)$ for $\Sigma$ such that $\Sigma\supset S$  and 
disjoint from $K$, and the zero map on $\cF(I|\Sigma)$ otherwise, 
then $\vphi_K$ commutes with $d_F=\bar{d}+\bar{\rho}$. 

The following properties are easily verified.
 If $K=K'\amalg K''$ then 
$\vphi_K=\vphi_{K''}\vphi_{K'}:
F(I\dbar S)\to  F(I-K\dbar S)$. 
If $K$ and $S'$ are disjoint, the
 maps $\sigma_{S\, S'}$ and $\vphi_K$ commute with each other, namely the 
following diagram commutes.
$$\begin{array}{ccc}
F(I\dbar S)&\mapr{\sigma_{S\, S'}} &F(I|S') \\
\mapd{\vphi_K}& &\mapdr{\vphi_K} \\
F(I-K\dbar S)&\mapr{\sigma_{S\, S'}} &F(I-K | S')
\end{array}
$$
\bigskip

\sss{\phiphicompatible}
{\bf Proposition.}\quad{\it  
The following square commutes in the derived category
($n=|I|$ and $c=|K|$):
$$\begin{array}{ccc}
 F(I) &\mapr{} & \cF(I | \ctop{I})[n-1] \\
 \mapd{\vphi_K} & & \mapdr{\psi_K} \\
 F(I-K) &\mapr{} & \cF(I-K |\ctop{I}-K)[n-c-1]
\end{array}
$$
where the horizontal maps are the canonical surjective maps, and 
the right vertical map is the one defined in (\cFIcJSigma), (6).
 }
\bigskip

\sss{\acyclicitysigma}
{\bf Proposition.}\quad{\it  
Let $R$, $J$ be disjoint subsets of $\ctop{I}$, with $J$ non-empty. 
Then the following  sequence of complexes is exact 
(the maps are alternating sums of 
the quotient maps $\sigma$) 
$$F(I|R)\mapr{\sigma}
\bop_{\stackrel{S\subset J}{|S|=1}\ } F(I| R\cup S)
\mapr{\sigma} \bop_{\stackrel{S\subset J}{|S|=2}} 
F(I|R\cup S)
\mapr{\sigma} \cdots \to  F(I|R\cup J)\to 0\,\,.$$}
\bigskip

\sss{\proofphiphicompatible}{\bf Proof of Proposition (\phiphicompatible).}
\quad 
We need two lemmas.
\bigskip

(\proofphiphicompatible.1) {\bf Lemma.}\quad{\it Let $\al: A\to C$ and $\be: B\to C$ be maps of complexes. 
Define another complex by 
$$F=\Cone(A\oplus B\mapr{(\al, -\be)} C)\,[-1]\,.$$ 
Let $p: F\to A$ and $q: F\to B$ be the projections. 
Then the square of complexes
$$\begin{array}{ccc}
  F  &\mapr{p}  &A  \\
\mapd{q} & &\mapdr{\al} \\
 B  &\mapr{\be}  &C
 \end{array}
$$  
is commutative up to homotopy. 
Further, if $\be$ is a quasi-isomorphism, then $p$ is also a quasi-isomorphism. }
\smallskip 

{\it Proof.}\quad By definition there is a distinguished triangle
in the homotopy category of complexes
$$F\to A\oplus B\to C\mapr{[1]}\,\,.$$
The claim hence follows.
\bigskip

For the next lemma, let  $p$ be a positive integer, and assume given a family of complexes $A_{i, j}$ for $i, j\in \ZZ$, that are zero unless 
$0\le i, j\le p$ and $0\le i+j\le p$. 
Also given maps of complexes $f: A_{i, j}\to A_{i-1, j}$ and $g: A_{i, j}\to A_{i, j-1}$
such that $fg=gf$ and $g$ are all quasi-isomorphisms. We thus have a commutative diagram
of complexes:
$$\begin{array}{cccccccc} 
  &&   &  &    &       &       &A_{p,0}  \\
  &&   &  &    &       &        &\mapdr{f}  \\
  &&   &  &   &A_{p-1,1}  &\mapr{g} &A_{p-1, 0} \\
  &&   &  &    &\mapd{f}&        &\mapdr{f}  \\ 
  && &A_{p-2,2} &\mapr{g} &A_{p-2,1} &\mapr{g} &A_{p-2,0} \\  
        &&   &\mapd{}  &     &\mapd{}  &    &\mapdr{}  \\ 
        &&\to&\vdots&&\vdots&&\vdots                                   \\
        &&\mapd{f}\phantom{\to\cdots\to}   &\mapd{f}  &     &\mapd{f}  &    &\mapdr{f}  \\ 
A_{0,p}&\mapr{g}&A_{0,p-1}\to\cdots\to&A_{02} &\mapr{g} &A_{01} &\mapr{g} &\,\,A_{00}\,\,.
  \end{array}$$
Let 
$$A_p=\bop_{i+j=p} A_{i, j}\,, \quad  A_{p-1}=\bop_{i+j=p-1} A_{i, j} $$
be the sum of the complexes; there are maps $f, g: A_p\to A_{p-1}$ that are 
sums of $f$ and $g$, respectively. 
We form a complex 
$$C_1:=\Cone(f-g: A_p\to A_{p-1})[-1]\,.$$
There is projection $C_1\to A_p$, hence also projections 
from $C_1$ to $A_{p, 0}$ and $A_{0, p}$. 

Consider also the diagram 
$$\begin{array}{ccc}
    &  &A_{p,0}  \\
 &  &\mapdr{f} \\
A_{0,p}   &\mapr{g}  &A_{00}
 \end{array}
$$  
where $f$ is the composition of the maps $f$, 
$A_{p, 0}\to A_{p-1, 0}\to \cdots \to A_{00}$, and 
similarly for $g$. 
Let 
$$C_2=\Cone(A_{p, 0}\oplus A_{0, p}\mapr{f-g} A_{00})[-1]\,.$$
There are projections from $C_2$ to $A_{p, 0}$ and $A_{0, p}$.
\bigskip 

(\proofphiphicompatible.2){\bf Lemma.}\quad{\it There is a quasi-isomorphism of complexes 
$C_1\to C_2$ compatible with the projections to $A_{p, 0}$ and $A_{0, p}$.}
\smallskip 

{\it Proof.}\quad 
For simplicity we give the argument in case $p=2$, the general case being 
similar. 
We have 
$$C_1=\Cone\left(A_{20}\oplus A_{02}\to \Cone(A_{11}\to A_{10}\oplus A_{01})\,\right)[-1]\,.$$
Since $g: A_{11}\to A_{10}$ is a quasi-isomorphism, the projection 
$\Cone(A_{11}\to A_{10}\oplus A_{01})\to A_{01}$ is a quasi-isomorphism
by the previous lemma, so the induced map 
$\Cone(A_{11}\to A_{10}\oplus A_{01})\to A_{00}$
is also a quasi-isomorphism. 
This last map induces a quasi-isomorphism $C_1\to C_2$. 
\bigskip

We now prove the commutativity of the diagram in (\phiphicompatible). 
We consider the case $K=\ctop{I}$ and 
$I=[1, n]$, and discuss later the necessary changes for a general $K$. 
Recall that $F(I)$ is the sum $\bop  \cF(I, \cJ|\Sigma)$ for all $(\cJ, \Sigma)$ with 
$ \cJ, \Sigma\subset \ctop{I}$ and $ \cJ\cap \Sigma=\emptyset$. 
So there are projections from $F(I)$ to $ \cF(I,\emptyset|\ctop{I})
= \cF(I|\ctop{I})[n-1]$ and 
$\cF(I, \emptyset|\emptyset)[n-1]=\cF(I, \emptyset)[n-1]$.
Let $S$ be the sum of $ \cF(I,  \cJ|\Sigma)$ with $( \cJ, \Sigma)$ satisfying either of
following conditions:
\smallskip 

$\bullet$\quad If $\Sigma=[i, n-1]$ for $3\le i\le n$ (so $\Sigma=\emptyset$ if $i=n$), 
and $ \cJ\cap [2, i-2]\neq \emptyset$ (equivalently, $ \cJ\neq \emptyset$ and $ \cJ\neq \{i-1\}$),

$\bullet$\quad If $\Sigma\neq [i, n-1]$ for $3\le i\le n$, and $ \cJ$ is any subset 
disjoint from $\Sigma$. 
\smallskip 

\noindent Using that $\sigma$ are quasi-isomorphisms, one easily verifies that $S$ is an acyclic subcomplex. 
Let $F_1(I)$ be the quotient complex of $F(I)$ by $S$; the projection 
$F(I)\to F_1(I)$ is a quasi-isomorphism compatible with projections to 
$ \cF(I|\ctop{I})[n-1]$ and $\cF(I, \emptyset)[n-1]$. 
We may display $F_1(I)$ as follows. 
$$\begin{array}{cccccccc}
   &&  &  &    &       &        &\scriptstyle{ \cF(I, \emptyset|[2, n-1])}  \\
   &&  &  &    &       &        &\mapdr{\rho}  \\
   &&  &  &   &\scriptstyle{ \cF(I, \emptyset|[3, n-1])} &\mapr{r} &{\scriptstyle
 \cF(I, \{2\}|[3, n-1])}\\
   &&  &  &    &\mapd{\rho}&        &  \\ 
   &&  &\scriptstyle{ \cF(I, \emptyset|[4, n-1])}   &\mapr{r}   &\scriptstyle{ \cF(I, \{3\}|[4, n-1])}  & &    \\
   &&  &\mapd{\rho}                 &    &                        &  &   \\ 
  &&\cdots  &   &  &              & &    \\
  &&\mapd{\rho}  &        &   &   & &     \\
\scriptstyle{ \cF(I, \emptyset)} &\mapr{r}&\scriptstyle{ \cF(I, \{n-1\})} &        &   &   & &   
  \end{array}
  $$
It consists of the two diagonals of a bigger diagram   
$$\begin{array}{ccccccccc}
 &  &&  &  &    &       &        &\scriptstyle{ \cF(I, \emptyset|[2, n-1])}  \\
 &  &&  &  &    &       &        &\mapdr{\rho}  \\
 &  &&  &  &   &\scriptstyle{ \cF(I, \emptyset|[3, n-1])} &\mapr{r} &{\scriptstyle
 \cF(I, \{2\}|[3, n-1])}\\
 &  &&  &  &    &\mapd{\rho}&        &\mapdr{\rho}  \\ 
&  
&
&  
&\scriptstyle{ \cF(I, \emptyset|[4, n-1])}   
&\mapr{r}  
&\scriptstyle{ \cF(I, \{3\}|[4, n-1])}  
&\mapr{r} 
&\scriptstyle{ \cF(I, [2, 3]|[4, n-1])}   \\
  & &&  &\mapd{\rho} &    &\mapd{\rho}                &  &\mapdr{\rho}   \\ 
  & &&  &\vdots &    &\vdots                &  &\vdots  \\ 
  & &\mapd{\rho}&  &\mapd{\rho} &    &\mapd{\rho}                &  &\mapdr{\rho}   \\   
\scriptstyle{ \cF(I, \emptyset)}
&\mapr{r}
&\scriptstyle{ \cF(I, \{n-1\})}
&\to\cdots\to  
&\scriptstyle{ \cF(I,[4, n-1] )}   
&\mapr{r}  
&\scriptstyle{ \cF(I, [3, n-1])}  
&\mapr{r} 
&\phantom{\,.}\scriptstyle{ \cF(I, [2, n-1])} \,.  
\end{array}
  $$
In other words, it is a triangular diagram of the kind mentioned above where we set
$$A_{i,j}= \cF(I, [2+j, n-i-1]\,|\,[n-i, n-1])$$
for $0\le i, j\le n-2$ and $i+j\le n-2$; the  $F_1(I)$ is 
the complex $C_1$ associated to this triangular diagram
(up to degree shift by 1).

Applying the second lemma, it is quasi-isomorphic to the complex $C_2$, 
which is the total complex 
 of the diagram 
$$\begin{array}{ccc}
   &  & \cF(I, \emptyset|[2, n-1])  \\
 &    &\mapdr{\rho}  \\
  \cF(I, \emptyset) &\mapr{r}  &\phantom{\, .} \cF(I, [2, n-1])\ .
 \end{array}
$$  
where $\rho=\rho_{[2, n-1]}$ and $r=r_{\emptyset, [2, n-1]}$. 
We denote the latter by $F_2(I)$. 
Thus we have  quasi-isomorphisms 
$$F(I)\to F_1(I)\to F_2(I)$$
compatible with the projections.

By the first lemma, the following square commutes in the derived category 
$$\begin{array}{ccc}
 F_2(I)   &\mapr{p}  & \cF(I|\ctop{I})[n-1]  \\
\mapd{q} &  &\mapdr{\rho}   \\
  \cF(I, \emptyset)[n-1]  &\mapr{r}  &\,\, \cF(I, \ctop{I})[n-1] \,\,.
 \end{array}
$$  
One may replace $F_2(I)$ on the upper right corner by $F(I)$, and still gets
a commutative square, which is part of the following diagram.
$$\begin{array}{ccc}
 F(I)   &\mapr{p}  & \cF(I|\ctop{I})[n-1]  \\
\mapd{q} & &\mapdr{\rho} \\
  \cF(I, \emptyset)[n-1]  &\mapr{r}  & \cF(I, \ctop{I})[n-1] \\
\mapd{\pi} & & \\ 
 \cF(\{1, n\}, \emptyset)[1]
 \end{array}
$$  
Here $\pi=\pi_{\ctop{I}}:  \cF(I, \emptyset)[n-1]\to  \cF(\{1, n\}, \emptyset)[1]$. 
Note one has $F(\{1, n\})= \cF(\{1, n\}, \emptyset)[1]$ and the composition 
of $q$ and $\pi$ is nothing but the map $\vphi: F(I)\to F(\{1, n\})$. 
On the other hand, recall that $\psi$ is defined as $\pi r^{-1}\rho$. 
The assertion hence follows.

In the case $K$ is a subset of $\ctop{I}$ in general, we argue similarly, and reduce to consider
the following commutative diagram:
$$\begin{array}{cccc}
& F(I)   &\mapr{p}  & \cF(I|\ctop{I})[n-1]  \\
& \mapd{q} & &\mapdr{\rho_K} \\
&  \cF(I, \emptyset|\ctop{I}-K)[n-1]   &\mapr{r}  & \cF(I, K|\ctop{I}-K)[n-1] \\
&\mapd{\pi_K} & & \\ 
F(I-K|\ctop{I}-K)=& \cF(I-K, \emptyset|\ctop{I}-K)[n-c-1]
 \end{array}
$$  
This concludes the proof. 
\bigskip

\sss{\proofacyclicitysigma}{\bf Proof of (\acyclicitysigma).}\quad
We recall that the complex $F(I|S)$ has a decreasing  
filtration by subcomplexes $F^k$, such that $\Gr^k F(I|S)$ is the direct sum of the complexes
$\cF(I|\Sigma)$ with $\Sigma \supset S$ and $|S|=k$.

The map $\sigma_{S\, S'}: F(I|S)\to F(I|S')$ respects the filtration, 
so we can take $\Gr^k$ of the sequence. 
For the term $F(I|R\cup S)$, the group $\cF(I|\Sigma)$ appears in its graded quotient
$\Gr^k$ if and only if 
$\Sigma\supset R\cup S$, equivalently, 
if $S\subset (\Sigma -R)\cap J$. 
So the graded quotient of the sequence is the sum, for $\Sigma$, of the following sequences
(where $T:=(\Sigma -R)\cap J$\,):
$$ \cF(I|\Sigma) \mapr{\sigma}
\bop_{\stackrel{S\subset T}{|S|=1} }  \cF(I|\Sigma)
\mapr{\sigma} \bop_{\stackrel{S\subset T}{|S|=2}} 
 \cF(I|\Sigma)
\mapr{\sigma} \cdots \to   \cF(I|\Sigma)\to 0\,\,.$$
Here the maps $\sigma$ are the alternating sums of the identity maps. 
The acyclicity of this sequence is consequence of the following obvious fact.

Let $A$ be an abelian group, $T$ be a (possibly empty) finite set. 
For each subset $S$ of $T$ let $A_S=A$ be a copy of $A$, and for 
an inclusion $S\subset S'$, let $\al_{S\, S'}: A_S\to A_{S'}$ be the identity map. 
Then the sequence 
$$A_\emptyset \to 
\bop_{|S|=1}  A_S\to 
\bop_{|S|=2}  A_S\to \cdots \to A_T\to 0 \,, $$
with the maps alternating sums of $\al_{S\, S'}$, is exact. 

Indeed, if $T$ is empty, the sequence is of the form $A_\emptyset \to 0$, 
which is trivially exact.  If $T$ is non-empty, the sequence is exact even with 
$0$ at left. 
\bigskip 

\sss{\symbols}{\it The class of symbols over $S$.}\quad      
Given a quasi-projective variety $S$, let $({\rm Smooth}/k\,, {\rm Proj}/S)$
be the category of smooth varieties $X$ equipped with projective maps to $S$. 
A {\it symbol} over $S$ is an object the form 
$$\bop_{\al\in A}(X_\al/S, r_\al)$$
where $X_\al$ is a collection of objects in $({\rm Smooth}/k\,, {\rm Proj}/S)$
indexed by a finite set $A$, and $r_\al\in \ZZ$. The class of objects over 
$S$ is denoted $\Symb(S)$. 
\bigskip

For a finite sequence of symbols of the form $(X_i/S, r_i)$, 
let
$$F\bigl((X_1/S, r_1), \cdots, (X_n/S, r_n)\bigr)$$
 be the complex $F(X_1, \cdots , X_n)$ with 
 respect to the sequence of dimensions
$$[i, i+1]\mapsto \dim X_{i+1}-r_{i+1}+r_i\,\,, i=1, \cdots, n-1.$$
For any sequence of symbols $K_i$, define the complex $F(K_1, \cdots, K_n)$
by linearity.

We thus have the complexes
$F(K_1, \cdots , K_n|S)$, and maps $\sigma_{S\, S'}$, $\vphi_K$ 
satisfying the properties as before. 
The class of objects $Symb(S)$, together with these complexes and maps, still denoted ${Symb}(S)$.
We shall define in \RefHathree\, the notion of a {\it quasi DG category}.   The content of this paper amounts to the verification 
that $Symb(S)$ forms a quasi DG category. 
The axioms for a quasi DG category correspond to the 
statements (\cpxFI), (\acyclicitysigma) proven in
 this section and the additional structure - 
generating set for the 
complex $F(K_1, \cdots , K_n)$, notion of properly intersecting elements, and distinguished 
subcomplexes with respect to constraints (see (\distcpxFI)\,)
as well as  diagonal cycles and 
diagonal extension (see \S 4).
\bigskip 

{\bf Remark.}\quad 
There is the structure of a category on $\Symb(S)$ as in \RefCH.
This is not used in this paper, and we refer the reader to \RefCH for details. 
Let us only say that the homomorphism group is 
$$\Hom ( (X/S, r), (Y/S, s)\,)=\CH_{\dim Y-s+r}(X\times_S Y)$$
and the composition is to be defined appropriately. 
\bigskip


\section{Distinguished subcomplexes}

In this section, complexes mean multiple complexes, and 
maps of complexes mean maps of multiple complexes
(appropriate totalization may be needed). 
\bigskip 

\sss{\distcpxZM}{\it Distinguished subcomplexes of $\bts_{i\in [1, r]}\Z(M_i)$.}
\quad  

(\distcpxZM.1)\quad Let $M$ be a smooth variety and $\Z(M, \cdot)$ its cycle complex. 
Recall that  $\Z(M, m)$ is $\ZZ$-free on the set  $ \set \Z(M, m)$ 
of irreducible non-degenerate admissible cycles. 
Thus a  non-zero element $u\in \Z(M, m)$ can be written as $\sum c_\nu \al_\nu$, with $c_\nu$ non-zero integers and $\al_\nu\in \set\Z(M, m)$;  we will say that $\al_\nu$ are the generators that appear in the basis expansion of $u$. 
We let 
$$\set\Z(M)=\amalg_m \set\Z(M, m)\,\eqno{(\distcpxZM.a)}$$
and refer to it as the set of generators for $\Z(M)$. 
(In this section we will use the symbol $\set C$  for the generating 
set of a free complex $C$.)
\bigskip

(\distcpxZM.2){\it Proper intersecting property.}\quad 
Let $(M, Y)=(M_i, Y_i)$ be a sequence of fiberings on  $[1, n]$:
$M_i \to Y_i \gets M_{i+1}$ for $i=1, \cdots, n-1$. 
Recall that we have the subvariety $M_{[1, n]}$ of $M_1\times\cdots
\times M_n$, and the 
projection $p_i: M_{[1, n]}\to M_i$ for each $i=1, \cdots, n$. 
 (We may take a sequence of fiberings on 
a finite ordered set $I$, but for convenience assume 
$I=[1, n]$.)

(1) Let $A$ be a subset of $[1, n]$. 
We define an element 
$(\al_i)_{i\in A}\in \prod_{A} \set\Z(M_i)$ 
to be {\it properly intersecting}
 with respect to $(M, Y)$ on $[1, n]$
if 
$\{p_i^* \al_i\,,\text{faces}\}_{i\in A}$ is a properly intersecting set of cycles in $M_{[1, n]}$. 
(See (\ZMtensor).)

(2) Generalizing this, an element 
$(u_i)_{i\in A}\in \prod_{A} \Z(M_i)$
 is defined to be  properly intersecting
if the following condition holds:
Let $A'=\{i\in A\mid u_i\neq 0\}$. 
For any choice of elements $\al_i\in \set\Z(M_i)$ for $i\in A'$, 
that appear in the basis expansion of $u_i$, the element 
$(\al_i)_{i\in A'}\in \prod_{A'} \set\Z(M_i)$
is properly intersecting.
\smallskip

The notion of proper intersection is compatible with 
boundary $\partial$. 
 Assume $(u_i )_{i\in [1, n]}$ is properly intersecting. 
For any $i\in [1, n]$ and any element $\be_\nu\in \Z(M_i)$ 
 appearing in the 
basis expansion of $\partial u_i$ (namely $\partial u_i=\sum c_{i\nu}\be_\nu$ with 
$\be_\nu\in \set\Z(M_i)$),  the element obtained by replacing 
$u_i$ by $\be_\nu$,
$$(u_1, \cdots, u_{i-1}, \be_\nu, u_{i+1}, \cdots, u_n)\in \prod_{[1, n]} \Z(M_i)$$
is properly intersecting.
The same holds for a properly intersecting set 
$(u_i )_{i\in A}$
indexed by $A\subset [1, n]$. 
\bigskip 

We will need more than  sequences of fiberings. 
 Let $\Tee_1, \cdots, \Tee_c$ be a finite sequence of finite ordered sets (each with cardinality $\ge 1$). 
Then $T:=\Tee_1\amalg \Tee_2\amalg \cdots\amalg \Tee_c$ is a finite ordered set equipped with a partition.
(Conversely a finite ordered set with a partition gives a finite sequence of finite ordered sets.)
We will assume $T=[1, n]$ for simplicity. 

Assume given, for each $\lambda$,  a sequence of fiberings 
$(M, Y)$ on $\Tee_\lambda$. 
[Example:  If $T=[1, 5]$ is partitioned into $[1, 3]$ and $[4, 5]$, we have two sequences of fiberings 
$$M_1 - M_2 - M_3 \qquad M_4-M_5$$
where each string represents a sequence of fiberings.]
We have the variety $M_{T_\lambda}\subset \prod_{i\in T_\lambda} M_i$ associated to the sequence $(M,Y)$ on $T_\lambda$, so let 
$$M_{[1, n]}=M_{T_1}\times\cdots \times M_{T_c}\subset \prod_{i\in [1, n]} M_i\,$$
be the product of them, and let $p_i: M_{[1, n]}\to M_i$ be the projections. 

For an element $(\al_i)_{i\in A}\in \prod_{A} \set\Z(M_i)$
or for  $(u_i)_{i\in A}\in \prod_{A} \Z(M_i)$ (where 
$A$ is a subset of $[1, n]$), 
we have the condition of proper intersection with respect to $(M, Y)$ on $\{T_\lambda\}$,
 defined exactly in the same manner as that with respect to 
 a sequence of fibering's. 
\smallskip

Given an ordered set of smooth varieties $M_1, \cdots, M_r$, we will consider a class of subcomplexes of 
$\Z(M_1)\ts\cdots\ts\Z(M_r)$ specified as follows. 
\bigskip 

(\distcpxZM.3){\bf Definition.}\quad  Let $M=(M_i)_{i\in [1, r]}$ be an ordered set of smooth varieties. 

(1) A {\it single condition of constraint} on the $M$ consists of (a)-(d):
\smallskip 

(a) An embedding of finite ordered sets $[1, r]\injto \BT$, and a partition $\{\BT_1, \cdots, \BT_c\}$ of $\BT$;

(b) A sequence of fiberings $(M, Y) $ on each $\BT_\lambda$ (with the varieties $M_i$ as given);

(c) A subset $P$ of $[1, r]$;

(d) A set of elements 
$f_k\in \Z(M_k)$ for $k\in \BT':=\BT-[1, r]$, such that $\{f_k\}_{k\in \BT'}$ is properly intersecting
with respect to $(M, Y)$ on $\{\BT_\lambda\}$.  
\smallskip 

\noindent 
We denote such data by 
$$\cC=([1, r]\injto \BT=\cup \BT_\lambda; (M,Y)\,\,\mbox{on $\{\BT_\lambda\}$}\, \, ; P\subset [1, r]; \{f_k\in \Z(M_k) \}_{k\in \BT'})\,.
\eqno{(\distcpxZM.b)}$$

(2) A {\it (general) condition of constraint} $\cC$ on $(M_i)_{i\in [1, r]}$ 
is a finite set of single conditions of constraint:  $\cC=\{\cC_j\}$.
\bigskip

%

To explain data (a) and (b) by an example, let  $r=5$, 
$\BT=\{1, 2, 3, a, b, 4, 5, c\}$
(elements in this order), with partition 
$\BT_1=\{1, 2\}$, $\BT_2=\{3, a, b, 4\}$, and $\BT_3=\{5, c\}$. 
For (b) one must give three sequences of fiberings
 $$M_1 - M_2 \quad M_3-M_a-M_b-M_4\quad M_5-M_c\,. $$
\bigskip

Given a single constraint $\cC$, consider the ($r$-tuple) subcomplex of 
$\Z(M_1)\ts\cdots\ts\Z(M_r)$
generated by 
$\al_1\ts\cdots\ts \al_r$ with $\al_i\in \set\Z(M_i)$ (see
(\distcpxZM.1)\,)
 satisfying the condition that 
$$\mbox{the element}\quad
(\,\al_i\,\, (i\in P),\quad f_k\,\, (k\in T) \,)\in\prod_{P\cup \BT'} \Z(M_i)\qquad \mbox{
be properly intersecting}$$
with respect to $(M,Y)$ on $\{\BT_\lambda\cap [1, r]\}$. 
We denote it by $[\Z(M_1)\ts\cdots\ts\Z(M_r)]_\cC$. 
For a general constraint $\cC=\{\cC_j\}$, we
define 
$$[\Z(M_1)\ts\cdots\ts\Z(M_r)]_\cC:=\bigcap_j [\Z(M_1)\ts\cdots\ts\Z(M_r)]_{\cC_j}\,\,$$
the intersection of the subcomplexes associated to each $\cC_j$, and call it
the  {\it distinguished subcomplex} of 
$\otimes_{i\in [1, r]} \Z(M_i)$
with respect to $\cC$. This terminology is justified by the proposition below.

We note that the class of subcomplexes of the form $[\otimes_{i\in [1, r]} \Z(M_i)]_\cC$
is closed under taking tensor product and finite intersection. 
\smallskip

As a special case of a distinguished subcomplex, assume we are given a sequence 
of fiberings $(M, Y)$ on $[1, r]$, and 
a subset  $P\subset [1, r]$. 
As we already considered in (\ZMtensor), let $C_P$ be the subcomplex of $\Z(M_1)\ts\cdots\ts\Z(M_r)$
generated by $\al_1\ts\cdots\ts\al_r$ with
 $\al_i\in \set\Z(M_i)$, where 
$\{\al_i\,(i\in P)\}$ is properly intersecting. 
We verified it to be a distinguished subcomplex in (\ZMtensor).
One can give an alternative proof using the proposition below, since $C_P$ is the subcomplex associated to a single condition of constraint $\cC$
consisting of $[1, r]=\BT$ with the trivial partition, the  given sequence of fibrerings $(M, Y)$, 
the given set $P$, and the empty set for $\BT'$. 
\bigskip


(\distcpxZM.3)\Prop\quad{\it The subcomplex $[\Z(M_1)\ts\cdots\ts\Z(M_r)]_\cC$ associated to $\cC$ is a 
distinguished subcomplex in the sense of (\distsubcpx). In particular, 
the inclusion 
$[\Z(M_1)\ts\cdots\ts\Z(M_r)]_\cC\injto \Z(M_1)\ts\cdots\ts\Z(M_r)$ is a quasi-isomorphism.
}\smallskip 

{\it Proof.}\quad For a single constraint $\cC$ as in (\distcpxZM.b), 
the subcomplex $[\Z(M_1)\ts\cdots\ts\Z(M_r)]_\cC$ is 
generated by 
$\al_1\ts\cdots\ts \al_r$ with $\al_i\in \set\Z(M_i)$ 
such that 
$$(\,p_i^*\al_i\,\, (i\in P),\quad M_{\BT}, \quad f_k\,\, (k\in \BT'), \quad\mbox{faces}\, \,)\in\prod_{P\cup \BT'} \Z(M_i)$$
(recall $M_{\BT}:= M_{\BT_1}\times\cdots\times M_{\BT_c}$)
is properly intersecting in $\prod_{i\in \BT} M_i$. 
By taking as $\{V_j\}$ the properly intersecting set of cycles 
$$\{ M_{\BT}, f_k\, (k\in \BT')\}$$ and applying (\distsubcpx.3), we see
it is a distinguished subcomplex.
\bigskip

\sss{\distcpxZMU}{\it  Distinguished subcomplexes of $\bts_{i\in [1, r]}\Z(M_i, \cU_i)$.}
\quad 
 As in \S 1, let $M$ be a smooth variety and $\cU$ be an open covering of an open set of $M$;  
we then have the complex  $\Z(M, \cU)$. 
 We can repeat all of (\distcpxZM) 
for the complex $\Z(M, \cU)$. 

Since $\Z(M, \cU)=\bop_I\Z(U_I)$ as a module, where $I$ varies over subsets (including the empty set)
 of the indexing set for $\cU$, an element $u\in \Z(M, \cU)$ can be uniquely written as  $\sum_I u_I$ with
$u_I\in \Z(U_I)$. 
\bigskip 

(\distcpxZMU.1)\quad Since $\Z(U_I)$ is $\ZZ$-free on the set $\set\Z(U_I)$, 
$\Z(M, \cU)$ is free on  the set 
$$\set\Z(M, \cU):= \amalg_I \set\Z(U_I)\,\,.\eqno{(\distcpxZMU.a)}$$
An element $u\in \Z(M,\cU)$ can be uniquely expanded as
 $\sum c_\nu \al_\nu$, with non-zero $c_\nu\in \ZZ$
 and $\al_\nu\in \set\Z(M, \cU)$.
\bigskip

(\distcpxZMU.2)\quad
Let $M_i$ be a sequence of fiberings indexed by $[1, n]$, and let $\cU_i$ 
be a finite covering of an open set $U_i\subset M_i$. 
We will call such data a sequence of fiberings {\it equipped 
with open coverings}, and denote it by $(M,Y,  \cU)$ on $[1, n]$. 
\smallskip 

(1)  Let $A$ be a subset of $[1, n]$, and 
$(\al_i)_{i\in A}\in \prod_{i\in A}\set\Z(M_i, \cU_i)$
be an element. It is said to be properly intersecting with 
respect to $(M, Y, \cU)$ on $[1, n]$ if, assuming 
$\al_i\in \set\Z((U_i)_{I_i})$ with $I_i$ a subset of the index set 
of $\cU_i$, 
$\{p_i^*\al_i\,,\text{faces}\}_{i\in A}$ is properly intersecting in 
the set 
$$\bigcap_{i\in A}  p_i^{-1}(\, (U_i)_{I_i}\,)
\subset M_{[1, n]}\,.$$

(2) Generalizing this as in (\distcpxZM.2), an element 
$(u_i)_{i\in A}\in \prod_{A} \Z(M_i, \cU_i)$
 is defined to be  properly intersecting
if the following condition holds:
Let $A'=\{i\in A\mid u_i\neq 0\}$. 
For any choice of elements $\al_i\in \set\Z(M_i, \cU_i)$ for $i\in A'$, 
that appear in the basis expansion of $u_i$, the element 
$(\al_i)_{i\in A'}\in \prod_{A'} \set\Z(M_i, \cU_i)$
is properly intersecting.
\bigskip

More generally, let $\Tee_1, \cdots, \Tee_c$ be a partition of 
$[1, n]$, and assume given 
a sequence of fiberings $(M, Y)$ on each $\Tee_\lambda$; 
we will refer to such data as
$(M, Y, \cU)$ on $\{\Tee_\lambda\}$ (the open coverings $\cU$ were given from the beginning).
As in (\distcpxZM.2), we can define the condition of proper intersection 
with respect to $(M, Y, \cU)$ on $\{\Tee_\lambda\}$ for elements 
$(\al_i)_{i\in A}\in \prod_{i\in A}\set\Z(M_i, \cU_i)$. 
\bigskip 

(\distcpxZMU.3){\bf Definition.}\quad 
Let $(M, \cU)$ be an ordered sequence of smooth varieties
 equipped with
open coverings on $[1, r]$. 

(1) A {\it single condition of constraint} on the $(M, \cU)$ on  $[1, r]$
is data
$$\cC=([1, r]\injto \BT=\cup \BT_\lambda; (M,Y, \cU)\,\,\mbox{on $\{\BT_\lambda\}$}\, \, ; P\subset [1, r]; \{f_k\in 
\Z(M_k, \cU_k) \}_{k\in \BT'})\,.
\eqno{(\distcpxZMU.b)}$$
consisting of:
\smallskip 

(a) An embedding of finite ordered sets $[1, r]\injto \BT$, and a partition $\{\BT_1, \cdots, \BT_c\}$ of $\BT$;

(b) A sequence of fiberings $(M, Y)$ on each $\BT_\lambda$;

(c) A subset $P$ of $[1, r]$;

(d) A set of elements 
$f_k\in \Z(M_k, \cU_k)$ for $k\in \BT'=\BT-[1, r]$, such that $\{f_k\}_{k\in \BT'}$ is properly intersecting
with respect to $(M, Y, \cU)$ on $\{\BT_\lambda\}$.  
\smallskip 

(2) A {\it  condition of constraint} $\cC$ on $(M, \cU)$ 
is a finite set of single conditions of constraint:  $\cC=\{\cC_j\}$.
\bigskip 

Given a single condition of constraint (\distcpxZMU.b)
let 
$[\Z(M_1, \cU_1)\ts\cdots\ts\Z(M_r, \cU_r)]_\cC$ 
be the (multiple) subcomplex generated by elements
$\al_1\ts\cdots\ts \al_r$ with $\al_i\in \set\Z(M_i, \cU_i)$, satisfying the condition that
$$\mbox{the element}\quad(\,\al_i\,\, (i\in P),\quad f_k\,\, (k\in T)\, \,)\in \prod_{P\cup \BT'} \Z(M_i, \cU_i)\qquad \mbox{
be properly intersecting.}$$
For a general constraint $\cC=\{\cC_j\}$ the 
corresponding distinguished subcomplex is defined, as in 
(\distcpxZM.2), 
by taking the intersection of the distinguished subcomplexes 
with respect to $\cC_j$.   
\smallskip 

Although this is not a distinguished subcomplex in the sense of 
(\distsubcpx), we extend the usage of the name since the proposition below shows it is a quasi-isomorphic subcomplex.
As a special case, given a sequence of fiberings 
equipped with open coverings $(M, Y, \cU)$ on $[1, r]$, and  
a subset $P\subset [1, r]$, one has the subcomplex 
$C_P$ of  $\Z(M_1,\cU_1)\ts\cdots\ts\Z(M_r, \cU_r)$
generated by $\al_1\ts\cdots\ts\al_r$ with
 $\al_i\in \set\Z(M_i, \cU_i)$, where $\{\al_i\,(i\in P)\}$ is properly intersecting. 
\bigskip 

(\distcpxZMU.4)\Prop\quad{\it The subcomplex 
$[\ts_{i\in [1, r]}\Z(M_i, \cU_i)]_\cC$ 
is a quasi-isomorphic subcomplex of $\ts_{i\in [1, r]}
\Z(M_i, \cU_i)$. 
}\smallskip 

{\it Proof.}\quad
The complex $\Z(M, \cU)=\bop_I\Z(U_I)$ has a decreasing 
filtration by subcomplexes 
$F^k= \bop_{|I|\ge k} \Z(U_I)$, and a successive quotient 
$\Gr^k_F$ is the sum of the complexes $\Z(U_I)$ with $|I|=k$.  
Thus 
$$\bigotimes_{i\in [1, r]}\Z(M_i, \cU_i)=\bop_{(I_1, \cdots, I_n)} \Z(U_{I_1})\ts\cdots\ts\Z(U_{I_n})$$
also has a decreasing
filtration by subcomplexes 
$$F^{k}= \bop 
\Z(U_{I_1})\ts\cdots\ts\Z(U_{I_n})\,,$$
the sum over $(I_1, \cdots, I_r)$ with $\sum |I_i|\ge k$, 
and a successive quotient
is the sum of the complexes 
$\Z(U_{I_1})\ts\cdots\ts\Z(U_{I_n})$ with 
$\sum |I_i|=k$. 

Note that from the definitions we have
$$[\bigotimes_{i\in [1, r]}\Z(M_i, \cU_i)]_\cC=
\bop [\bigotimes_{i\in [1, r]}\Z(\cU_{I_i})]_\cC\,.$$
So there is also the filtration $F^{k}$ on 
the subcomplex $[\ts_{i\in [1, r]}\Z(M_i, \cU_i)]_\cC$, and 
the inclusion $[\ts_{i\in [1, r]}\Z(M_i, \cU_i)]_\cC \injto
\ts_{i\in [1, r]}\Z(M_i, \cU_i)$ respects the filtrations.
In successive quotients, 
we have induced inclusions 
$$[\bigotimes_{i\in [1, r]}\Z(U_{I_i})]_\cC\injto 
\bigotimes_{i\in [1, r]}\Z(U_{I_i})\,.$$
By (\distcpxZM.4) applied to $(U_{I_i})_{i\in [1, r]}$, these inclusions are quasi-isomorphisms. 
\bigskip

\sss{\almstdisj}\quad In order to transcribe the previous
subsections in the setting of \S 2, we need some preliminaries. 
\bigskip

(\almstdisj.1){\bf Definition.}\quad 
Let $I$ be a finite ordered set, and 
let $\{I_\al\}$ be a collection of sub-intervals (of cardinality $\ge 2$)
 of $I$,  indexed by a finite set $\cA$.  We say $\{I_\al\}$ is {\it almost disjoint} if 
for  distinct elements $\al$, $\be$ of $\cA$, 
$I_\al\cap I_\be$ consists of at most one element. 
Then there is a total ordering $<$ on the set $\cA$ 
given by $\al < \be$ if and
only if $\init(I_\al)<\init(I_\be)$. 
When $\al < \be$, we have either $I_\al$ and $I_\be$ disjoint, or 
$\term(I_\al)=\init(I_\be)$.

In case the indexed set is $\cA=[1, r]$, the ordering on $\cA$ 
thus defined
coincides with the natural one if and only if $\term(I_j)\le \init(I_{j+1})$ for each $j=1, \cdots, r-1$. We will then say that
$\{I_1, \cdots, I_r\}$ is {\it ordered}. 
\bigskip 

(\almstdisj.2)\quad
Let $I_1, \cdots, I_r$ be an ordered almost disjoint set of sub-intervals of $I$. \smallskip 

(1) The union of $I_j$ is $I$ 
 if and only if
$\term (I_j)= \init(I_{j+1})$ or $\term (I_j)+1= \init(I_{j+1})$
for $j=1, \cdots, r-1$, and $\init(I_1)=\init(I)$, $\term(I_r)=\term(I)$.  
We then say $\{I_j\}$ is a {\it pseudo-segmentation}. 
It is a segmentation of $I$ if and only if 
$\term (I_j)= \init(I_{j+1})$ for $j=1, \cdots, r-1$, and $\init(I_1)=\init(I)$, $\term(I_r)=\term(I)$.  
\smallskip 

(2) 
We have a partition $\Tee_1, \cdots, \Tee_c$ of $[1, r]$ 
such that each $\Tee_\lambda$ is a sub-interval maximal with 
the property that if $j, j+1\in \Tee_\lambda$, then $I_j \cap I_{j+1}\neq \emptyset$. We refer to $\{T_\lambda\}$ as the partition 
associated to the almost disjoint set of sub-intervals $\{I_j\}$.
[Example: If $r=5$, and $\term(I_1)=\init(I_2), \term(I_2)<\init(I_3), 
\term(I_3)<\init(I_4), \term(I_4)=\init(I_5)$, then $[1, 5]$ is partitioned into the sub-intervals $[1, 2], \{3\}$, and  $[4, 5]$.]

If we set 
$$I^\lambda =\bigcup_{j\in \Tee_\lambda} I_j, \qquad \lambda=1, \cdots, s,$$
then $I^1, \cdots, I^s$ is a disjoint set of sub-intervals
(it is a partition of $I$ in case $\{I_j\}$ is a pseudo-segmentation), each $I^\lambda$ segmented into sub-intervals ${I^\lambda}_1, \cdots, {I^\lambda}_{r_\lambda}$, so that the totality of the sub-intervals 
$\{{I^\lambda}_k\}$ coincides with $\{I_j\}$.
(Conversely an ordered,  disjoint  set of sub-intervals 
$I^1, \cdots, I^s$, plus a segmentation of each $I^\lambda$, 
gives an almost disjoint set of sub-intervals of $I$.)
\bigskip 


\sss{\distcpxcFI}{\it  Distinguished subcomplexes of $\cF(I_1)\ts\cdots\ts \cF(I_r)$.}
\quad 
Assume given a sequence of varieties $X_i/S$ on a finite ordered
set $I$\,
(see \S 2).
\smallskip 

(\distcpxcFI.1)\quad Let $\cJ$ be a subset of $\ctop{I}$. 
One has the associated smooth variety $M=X_I^{\cJ}$ and its cycle complex $\cF(I, \cJ)=\Z(X_I^{\cJ})$. 
The complex $\cF(I, \cJ)$ is free on the set
$$\set\cF(I, \cJ):=\set\cZ(X_I^{\cJ})\,.\eqno{(\distcpxcFI.a)}$$
\bigskip

(1) (See \S 2.) Assume given 
a segmentation of $I$ into sub-intervals of $I$,  
 $\{I_1, \cdots, I_r\}$, and a set of subsets $\cJ_j \subset \ctop{I_j}$. 
Then there corresponds a sequence of fiberings on $T=[1, r]$ given by 
$j\mapsto M_j=X^{\cJ_j}_{I_j}$, 
and there are open coverings $\cU_j=\cU(\cJ_j)$ on $M_j$. 


(2) More generally, assume given an almost disjoint set of
sub-intervals 
$\{I_1, \cdots, I_r\}$ of $I$,
and subsets $\cJ_j \subset \ctop{I_j}$;
such data will be denoted by $(X/S; I; \{I_j\supset \cJ_j\})$.
Then one has a partition $\Tee_1, \cdots, \Tee_c$ of $[1, r]$ as in (\almstdisj.2), and on each $\Tee_\lambda$ we have a sequence of fiberings $(M, Y, \cU)$; of course this depends only on 
the union $\cup I_j$. 
\bigskip 

%
%
%

(\distcpxcFI.2)\quad Let $(X/S; I; \{I_j\supset \cJ_j\})$ be as in (2) above, and 
$(M,Y, \cU)$ on $\{\Tee_\lambda\}$ be the associated sequence of fiberings. 
Let $A$ be a subset of $[1, r]$.  
By (\distcpxZMU.2), we have the condition of proper intersection for elements 
$$(u_j)\in \prod_A \Z(M_j, \cU_j)=\prod_A \cF(I_j, \cJ_j)\,.$$
\bigskip

(\distcpxcFI.3){\bf Definition.}\quad 
Let $(X/S; I; \{I_j\supset \cJ_j\})$ be as in (2) above.  
 
(1) A single condition of constraint on $(X/S; I; \{I_j\supset \cJ_j\})$ with respect to $\cF(I_j, \cJ_j)$ 
is data 
$$\cC=(I\injto \BI; \mbox{$X$  on $\BI$}; P; \{f_k\in \cF(J_k, \cJ'_k)\})\eqno{(\distcpxcFI.b)}$$ 
consisting of:
\smallskip 

(a)  An embedding of finite ordered sets $I\injto \BI$ such that each $I_j$ remains a sub-interval of $\BI$; 

(b)  A sequence of varieties $X$ (over $S$) on $\BI$, extending the given $X$ on $I$; 

(c) A subset $P$ of $[1, r]$;

(d) A collection of sub-intervals $\{J_k\}_{k=1, \cdots, t}$  
of $\BI$  such that 
$\{I_j\}\cup\{ J_k\}$ is almost disjoint in $\BI$, 
subsets $\cJ'_k\subset\ctop{J_k}$, 
and elements $f_k\in \cF(J_k, \cJ'_k)$ that are properly 
intersecting (namely $(f_k)$ is properly intersecting).  
\smallskip 

(2) By a condition of constraint on $(X/S; I; \{I_j\supset \cJ_j\})$ 
with respect to $\cF(I_j, \cJ_j)$
we mean a finite set of single conditions of constraints 
$\cC=\{\cC_j\}$. 
\bigskip

Recall from (\distcpxcFI.1), (2) that associated to $(X/S; I; \{I_j\supset \cJ_j\})$
is $(M, Y, \cU)$ on $\{T_\lambda\}$, in particular 
a sequence of varieties with open coverings $(M, \cU)$ on $[1, n]$. 
A condition of constraint on $(X/S; I; \{I_j\supset \cJ_j\})$
in the above sense gives
a condition of constraint on the corresponding  
$(M, \cU)$ on $[1, r]$
in the sense of (\distcpxZMU.3), as follows. 
For a single condition of constraint (\distcpxcFI.b):
\smallskip 

(a)  Let $\BT=[1, r]\amalg [1, t]$, the index set for the intervals $\{I_j\}\cup \{J_k\}$, 
ordered as in (\almstdisj.1); let $\{\BT_\lambda\}$ be the 
partition of $\BT$ associated to $\{I_j\}\cup \{J_k\}$
(see (2) of (\almstdisj.2)\,). 
There is an inclusion $[1, r]\injto \BT$. 
[Example.  Let $r=5$, $t=3$, and $\term(I_1)=\init(I_2)$, 
$\term(I_2)=\init(I_1)$
$\term(J_1)=\init(J_2)$, $\term(J_2)=\init(I_3)$, 
$\term(I_3)<\init(J_3)$, $\term(J_3)=\init(I_4)$, 
$\term(I_4)=\init(I_5)$. 
The figure below illustrate this. 
We have $\BT=\{1, 2, 1', 2', 3, 3', 4, 5\}$ in this order, and 
$\BT_1=\{1, 2, 1', 2', 3\}$, $\BT_2=\{3', 4, 5\}$.]

\vspace*{0.3cm}
\hspace*{1cm}
\unitlength 0.1in
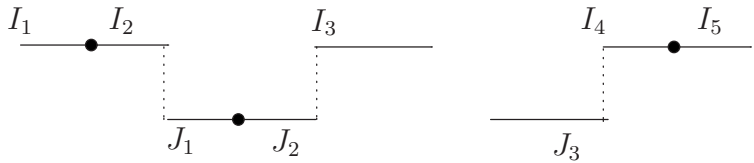
\begin{picture}( 38.9000,  6.6000)(  7.4000,-12.8000)
%
{\color[named]{Black}{%
\special{pn 8}%
\special{pa 810 810}%
\special{pa 1580 810}%
\special{fp}%
\special{pa 1580 810}%
\special{pa 1580 810}%
\special{fp}%
}}%
%
{\color[named]{Black}{%
\special{pn 8}%
\special{pa 1580 1200}%
\special{pa 2350 1200}%
\special{fp}%
\special{pa 2350 1200}%
\special{pa 2350 1200}%
\special{fp}%
}}%
%
{\color[named]{Black}{%
\special{pn 8}%
\special{pa 2350 820}%
\special{pa 2960 820}%
\special{fp}%
\special{pa 2960 820}%
\special{pa 2960 820}%
\special{fp}%
}}%
%
{\color[named]{Black}{%
\special{pn 8}%
\special{pa 3860 820}%
\special{pa 4630 820}%
\special{fp}%
\special{pa 4630 820}%
\special{pa 4630 820}%
\special{fp}%
}}%
%
{\color[named]{Black}{%
\special{pn 8}%
\special{pa 3270 1200}%
\special{pa 3880 1200}%
\special{fp}%
\special{pa 3880 1200}%
\special{pa 3880 1200}%
\special{fp}%
}}%
%
{\color[named]{Black}{%
\special{pn 4}%
\special{sh 1}%
\special{ar 1170 810 6 6 0  6.28318530717959E+0000}%
\special{sh 1}%
\special{ar 1170 810 6 6 0  6.28318530717959E+0000}%
\special{sh 1}%
\special{ar 1170 810 6 6 0  6.28318530717959E+0000}%
}}%
%
{\color[named]{Black}{%
\special{pn 0}%
\special{sh 1.000}%
\special{ia 1180 810 30 30  0.0000000 6.2831853}%
}}%
{\color[named]{Black}{%
\special{pn 8}%
\special{ar 1180 810 30 30  0.0000000 6.2831853}%
}}%
%
{\color[named]{Black}{%
\special{pn 0}%
\special{sh 1.000}%
\special{ia 1950 1200 30 30  0.0000000 6.2831853}%
}}%
{\color[named]{Black}{%
\special{pn 8}%
\special{ar 1950 1200 30 30  0.0000000 6.2831853}%
}}%
%
{\color[named]{Black}{%
\special{pn 0}%
\special{sh 1.000}%
\special{ia 4230 820 30 30  0.0000000 6.2831853}%
}}%
{\color[named]{Black}{%
\special{pn 8}%
\special{ar 4230 820 30 30  0.0000000 6.2831853}%
}}%
\put(7.4000,-7.5000){\makebox(0,0)[lb]{$I_1$}}%
\put(12.7000,-7.5000){\makebox(0,0)[lb]{$I_2$}}%
\put(23.3000,-7.6000){\makebox(0,0)[lb]{$I_3$}}%
\put(37.4000,-7.5000){\makebox(0,0)[lb]{$I_4$}}%
\put(43.5000,-7.5000){\makebox(0,0)[lb]{$I_5$}}%
\put(15.7000,-13.8000){\makebox(0,0)[lb]{$J_1$}}%
\put(21.2000,-13.9000){\makebox(0,0)[lb]{$J_2$}}%
\put(35.9000,-14.1000){\makebox(0,0)[lb]{$J_3$}}%
%
{\color[named]{Black}{%
\special{pn 8}%
\special{pa 1560 820}%
\special{pa 1560 1200}%
\special{dt 0.045}%
\special{pa 1560 1200}%
\special{pa 1560 1190}%
\special{dt 0.045}%
}}%
%
{\color[named]{Black}{%
\special{pn 8}%
\special{pa 2360 820}%
\special{pa 2360 1200}%
\special{dt 0.045}%
\special{pa 2360 1200}%
\special{pa 2360 1190}%
\special{dt 0.045}%
}}%
%
{\color[named]{Black}{%
\special{pn 8}%
\special{pa 3860 830}%
\special{pa 3860 1210}%
\special{dt 0.045}%
\special{pa 3860 1210}%
\special{pa 3860 1200}%
\special{dt 0.045}%
}}%
\end{picture}%
\vspace{0.5cm}

(b) Applying (\distcpxcFI.1), (2)  to 
$(X/S; \BI; \{I_j\supset \cJ_j\}\cup
\{J_k\supset \cJ'_k\})$, one obtains a 
a sequence of fiberings $(M, Y, \cU)$ on $\{\BT_\lambda\}$. 
Concretely, $(M, Y, \cU)$ consists of the varieties $M_j=X_{I_j}$ and 
$M_k=X_{J_k}$ with superscripts $\cJ_j$ and $\cJ'_k$, together with open coverings 
$\cU(\cJ_j)$, $\cU(\cJ_k')$, respectively. 
The corresponding cycle complexes are $\cF(I_j, \cJ_j)$ and 
$\cF(J_k, \cJ'_k)$.  

(c) Take the given set $P\subset [1, r]$. 

(d) For $k\in [1, t]$, we are given elements 
$f_k\in \Z(M_k, \cU_k)$, that are properly intersecting. 
\smallskip 

\noindent To a general constraint  $\cC$ on $(X/S; I; \{I_j\supset \cJ_j\})$, there corresponds a general 
constraint on $(M, \cU)$. 

According to (\distcpxZMU) one has the corresponding distinguished subcomplex 
$$[\cF(I_1, \cJ_1)\ts\cdots\ts\cF(I_r, \cJ_r)]_\cC\subset \cF(I_1, \cJ_1)\ts\cdots\ts\cF(I_r, \cJ_r)\,,$$ 
called the distinguished subcomplex with respect to $\cC$. 
Recall that it is generated by $\al_1\ts\cdots\ts\al_r$, 
$\al_j\in \set\cF(I_j, \cJ_j)$, such that 
$$(\al_j\,\, (j\in P), \quad f_k\,\,(k\in [1, t])\,)\in \prod_P\cF(I_j, \cJ_j)\times \prod_{[1, t]} \cF(J_k, \cJ'_k)$$
is properly intersecting with respect to $(M, Y, \cU)$ on 
$\{\BT_\lambda\}$. 

As a special case of a distinguished subcomplex, 
we have $\cF(I_1, \cJ_1)\hts\cdots\hts\cF(I_r, \cJ_r)$ as 
considered in \S 2.
\smallskip 

{\it Remark.}\quad In (d) of (\distcpxcFI.3), 
one may give elements 
$f_k\in \cF(J_k)$ that are properly 
intersecting (in the sense to be defined in (\distcpxcFI.4).)
This gives a general condition of constraint as follows. 
Write 
$f_k=\sum f_k(\cJ'_k)$ with $f_k(\cJ'_k)\in \cF(J_k, \cJ'_k)$. 
For any choice of subsets $\{\cJ'_k\}_{k=1, \dots, t}$, we have a single constraint in the sense of 
(\distcpxcFI.3). 
The collection of these, as we take all choices $\{\cJ'_k\}$, 
gives us a constraint.  
\bigskip 

(\distcpxcFI.4){\it Variant for $\cF(I)$.}\quad 
We can repeat the discussion (\distcpxcFI.1)-(\distcpxcFI.3)
with 
$\cF(I, \cJ)$ replaced by $\cF(I)$. 

The complex $\cF(I)$ is the direct sum of submodules $\cF(I, \cJ)$ for $\cJ\subset \ctop{I}$. 
Thus it is free on the set 
$$\set\cF(I):=\amalg \,\set\cF(I, \cJ)\,, \eqno{(\distcpxcFI.c)}$$
the disjoint union of the sets of generators for $\cF(I, \cJ)$. 

Let $\{I_j\}_{j=1, \cdots, r}$ be an almost disjoint set of 
sub-intervals of $I$. 
Let $A$ be a subset of $[1, r]$.  
An element 
$(u_j)\in \prod_{j\in A} \cF(I_j)$ is said to be properly intersecting with respect to $(X/S; I; \{I_j\})$ 
if, 
writing $u_j=\sum u_j(\cJ_j)$ with $u_j(\cJ_j)\in \cF(I_j, \cJ_j)$, the element 
$$(u_j(\cJ_j))_{j\in A}\in \prod_A\cF(I_j, \cJ_j)$$
is properly intersecting with respect to $(X/S; I; \{I_j\supset \cJ_j\})$
in the sense of (\distcpxcFI.2), for 
any choice of subsets $\cJ_j$. 
Note that the data $(X/S; I; \{I_j\})$ alone does not specify 
a sequence of fiberings. 

A single condition of constraint on $(X/S; I; \{I_j\})$
with respect to $\cF(I_j)$ is data
$$\cC=(I\injto \BI; \mbox{$X$  on $\BI$}; P; \{f_k\in \cF(J_k)
\}_{k=1, \cdots, t})\eqno{(\distcpxcFI.d)}$$ 
satisfying the same conditions as in (\distcpxcFI.3) except 
in (d) one must give properly intersecting elements 
$f_k\in \cF(J_k)$ instead of elements in $\cF(J_k, \cJ'_k)$. 
A condition of constraint is a finite set of single constraints. 

A condition of constraint $\cC$ specifies
a subcomplex of $\cF(I_1)\ts\cdots\ts\cF(I_r)$ called 
the distinguished subcomplex with respect to $\cC$. 
For a single constraint as (\distcpxcFI.d), the corresponding subcomplex 
$[\ts_{j\in [1, r]} \cF(I_j)]_\cC$ is defined to be  
the submodule generated by 
 $\al_1\ts\cdots\ts \al_r$ with $\al_j\in \set\cF(I_j)$
satisfying the condition that 
$$(\, \al_j\,\,(j\in P), \quad f_k\,\,(k\in [1, t]\,)\,)\in 
\prod_P\cF(I_j)\times \prod_{[1, t]} \cF(J_k)$$
be properly intersecting
 with respect to $(X/S; \BI; \{I_j\}\cup\{J_k\})$. 
From this definition one has equality
$$[\cF(I_1)\ts\cdots\ts\cF(I_r)]_\cC=\bop
[\cF(I_1, \cJ_1)\ts\cdots\ts\cF(I_r, \cJ_r)]_\cC\,,\eqno{(\distcpxcFI.e)}$$
the sum over $(\cJ_1, \cdots, \cJ_r)$, 
of the distinguished subcomplexes with respect to the same $\cC$
defined in (\distcpxcFI.3)
(see also Remark at the end of (\distcpxcFI.3)).
The inclusion of a distinguished subcomplex is a quasi-isomorphism. 

As a special case of a distinguished subcomplex, 
we have $\cF(I_1)\hts\cdots\hts\cF(I_r)$ (see \S 2).
\bigskip 

(\distcpxcFI.5){\it Variant for $\cF(I|\Sigma)$.}\quad
We give a further variant, where $\cF(I)$ is replaced with 
$\cF(I|\Sigma)$. 

If $\Sigma$ is a subset of $\ctop{I}$ and $I_1, \cdots, I_c$
is the corresponding segmentation, recall that   
$\cF(I|\Sigma)=\cF(I_1)\hts\cdots\hts\cF(I_c)$ is the subcomplex of $\cF(I_1)\ts\cdots\ts
\cF(I_c)$ with the set of generators
$$\set\cF(I|\Sigma):=\{(\al_i)_{i=1, \cdots, c}\in
\prod \set\cF(I_i)|\,\,  (\al_i)\quad\mbox{is properly intersecting.}\,\}.
\eqno{(\distcpxcFI.f)}$$ 

Assume now $I_1, \cdots, I_r$ is an almost disjoint set of 
sub-intervals of $I$, and 
$\Sigma_j$ are subsets of $\ctop{I_j}$; we denote the data
by $(X/S; I; \{I_j|\Sigma_j\})$. 
Each $I_j$ is segmented by $\Sigma_j$ to $I_{j,1}, \cdots, 
I_{j, s_j}$. 
Then an element ($A$ is a subset of $[1, r]$)
$$(\al_j)_{j\in A}\in \prod_{A} \set\cF(I_j|\Sigma_j)\quad\mbox{with}\quad
\al_j=(\al_{j k}) \in \prod_{k\in [1,, s_j]} \set\cF(I_{j k})
$$
is said to be properly intersecting if the element 
$$(\al_{j k})\in\prod_{j\in A}\prod_k \set\cF(I_{jk})\,$$
is  properly intersecting 
with respect to $(X/S; I; \{I_{jk}\})$.
Also, as in (\distcpxZM), (2), we can define the condition 
of proper intersection for an element 
$$(u_j)_{j\in A}\in \prod_{A} \cF(I_j|\Sigma_j)\,.$$

A single condition of constraint on the $(X/S; I; \{I_j|\Sigma_j \})$  is data 
$$\cC=(I\injto \BI; \mbox{$X$  on $\BI$}; P; \{f_k\in \cF(J_k|\Sigma'_k)\})\eqno{(\distcpxcFI.g)}$$ 
satisfying the same conditions as in (\distcpxcFI.3) except 
in (d) one must give sub-intervals $J_k$ such that $\{I_j\}\cup
\{J_k\}$ is almost disjoint,  subsets $\Sigma'_k\subset
\ctop{J_k}$, and
properly intersecting elements 
$f_k\in \cF(J_k|\Sigma'_k)$.  A condition of constraint
is a finite set of single constraints. 

Just as in (\distcpxcFI.3) or (\distcpxcFI.4), a condition of constraint specifies the corresponding distinguished subcomplex 
$$[\cF(I_1|\Sigma_1)\ts\cdots\ts\cF(I_r|\Sigma_r)]_\cC\subset \cF(I_1|\Sigma_1)\ts\cdots\ts\cF(I_r|\Sigma_r)\,.$$
\bigskip

\sss{\distcpxFI}{\it  Distinguished subcomplexes of $F(I_1)\ts\cdots\ts F(I_r)$.}
\quad We keep the assumption from the previous subsection
that $X/S$ is a sequence of varieties on $I$. 
We give a variant of (\distcpxcFI) for $F(I)$. 
\smallskip 

(\distcpxFI.1)\quad 
We have $F(I)=\bop \cF(I|\Sigma)$, the sum over subsets 
$\Sigma$ of $\ctop{I}$. 
So it is $\ZZ$-free over the set 
$$\set F(I):=\amalg\, \,\set\cF(I|\Sigma)\,.\eqno{(\distcpxFI.a)}
$$

(\distcpxFI.2)\quad 
If $\{I_j\}_{j=1, \cdots, r}$ is an almost disjoint set of sub-intervals of  
$I$, and $A\subset [1, r]$, 
we have the notion of proper intersection 
with respect to $(X/S; I; \{I_j\})$ for an 
element 
$$(u_j)\in \prod_{j\in A} F(I_j)\,$$
 by naturally extending that for an element of  
$\prod\cF(I_j|\Sigma_j)$ in (\distcpxcFI.5).  
\bigskip 

(\distcpxFI.3){\bf Definition. }\quad 
A single condition of constraint on $(X/S; I; \{I_j\})$
with respect to the complexes $F(I_j)$ can be 
defined as in (\distcpxcFI.3), except in (d) one must
replace $f_k\in \cF(J_k, \cJ_k')$ with properly intersecting elements $f_k\in F(J_k)$. 
Thus it is of the form 
$$\cC=(I\injto \BI; \mbox{$X$  on $\BI$}; P; \{f_k\in F(J_k)\}_{k=1, \cdots, t})\,.$$
A condition of constraint is a finite set of single constraints. 
\bigskip 

For a condition of constraint $\cC$ one has the associated 
distinguished subcomplex
$$[F(I_1)\ts\cdots\ts F(I_r)]_\cC\injto F(I_1)\ts\cdots\ts F(I_r)\,.$$
For a single constraint, it is the subcomplex generated by 
 $\al_1\ts\cdots\ts \al_r$ with $\al_j\in \set F(I_j)$
satisfying the condition that 
$$(\, \al_j\,\,(j\in P), \quad f_k\,\,(k\in [1, t]\,)\,)\in 
\prod_P F(I_j)\times \prod_{[1, t]} F(J_k)$$
be properly intersecting
 with respect to $(X/S; \BI; \{I_j\}\cup\{J_k\})$. 
\smallskip 

As a special case, let $\{I_j\}$ be a segmentation of $I$, 
and let $F(I_1)\hts\cdots\hts F(I_r)$
be the subcomplex generated by elements 
$\al_1\ts\cdots\ts\al_r$ with
$\al_j\in \set F(I_j)$ where $(\al_j)_{j\in [1, r]}$ is properly 
intersecting. 
If $S$ corresponds to the segmentation $\{I_j\}$, then we have
$$F(I_1)\hts\cdots\hts F(I_r)=F(I|S)\,.\eqno{(\distcpxFI.b)}$$
(where the complex $F(I|S)$ is defined in \S 2 as the sum 
of $\cF(I|\Sigma)$ for $\Sigma\supset S$). 
Indeed the left hand side is generated by properly 
intersecting elements 
$\al_1\ts\cdots\ts\al_r$ with 
$\al_j \in \set\cF(I_j |\Sigma_j)$, 
namely by elements $\al\in \set\cF(I|\Sigma)$ for $\Sigma\supset S$. 
\bigskip


\section{The diagonal cycle and the diagonal extension}

We keep the notation of \S 2.  
In this section $I$ denotes a finite ordered set of 
cardinality $\ge 2$, and we often take $I=[1, n]$ 
for simplicity.
\bigskip

\sss{\diagcycle} {\it The diagonal cycles $\bDelta(I)$.}
\quad 
Let $X$ be a smooth variety, projective over $S$, and $X_i=X$ be a 
constant sequence of varieties on $I=[1, n]$, $n\ge 2$
(or on a finite ordered set of cardinality $\ge 2$).  
There is the diagonal embedding $\Delta: 
X\to X\times_S\cdots\times_S X$; denote the image of 
the fundamental class of $X$ by $\Delta([1, n])\in \Z(X
\times_S\cdots\times_S X)$. 
Recall that there is a natural quasi-isomorphism
$$\iota:
\Z(X\times_S\cdots\times_S X)\to \Z(X^\emptyset_{[1, n]}, \{U^\emptyset_{[1, n]}\})=\cF([1, n], 
\emptyset)\,\,.$$
 We  use the same $\Delta([1, n])$ to denote 
its image under this map. It thus consists of $\Delta([1, n])$
 in $\Z(X^\emptyset_{[1, n]})$, 
and the zero element in $\Z(U^\emptyset_{[1, n]})$. 
Similarly for a finite ordered set $I$ of cardinality $\ge 2$, 
 we have an element
$\Delta(I)\in \cF(I, \emptyset)^0$ which is a cocycle;
viewed as an element of $\cF(I)$, it is of degree 1 (of bi-degree $(1, 0)$), and is 
a cocycle:
$\partial(\Delta(I)\,)=0$. 

For a subset $\Sigma\subset\ctop{I}$, letting  $I_1, \cdots, I_c$ be 
 the segmentation of $I$ given by $\Sigma$, 
one verifies immediately that the tensor product 
$$\Delta(I|\Sigma):=\Delta(I_1)\ts\Delta(I_2)\ts\cdots\ts\Delta(I_c)
\in \cF(I_1, \emptyset)\ts\cdots\ts\cF(I_c, \emptyset)$$
 is in the subcomplex $\cF(I, \emptyset|\Sigma)$, 
 namely that $\{\Delta(I_1), \cdots, \Delta(I_c)\}$ is properly 
 intersecting. 
It is a cocycle of degree $c$: $\partial (\Delta(I|\Sigma)\,)=0$. 
As an element of the complex $\cF(I|\Sigma)^{\shift}$ (see (\cpxFI)\,), its degree is $0$. 
It is easily verified that   
the elements $\Delta(I)$ are subject to the following properties:
\smallskip 

(1) For $k\in \Sigma$, $\rho_k(\Delta(I|\Sigma)\, )=r_k
 (\Delta(I|\Sigma-\{k\})\, )
$ in $\cF(I, \{k\}|\Sigma-\{k\})$. 
Here the map $r_{\cJ, \cJ'}: \cF(I, \cJ|\Sigma)
\to \cF(I, \cJ'|\Sigma)$ is written $r_k$ if $\cJ'=\cJ\cup\{k\}$.

(2) For $K\subset\ctop{I}- \Sigma$, $\pi_K(\Delta(I|\Sigma)\,)
=\Delta(I-K|\Sigma)$ in $\cF(I- K|\Sigma)$.
\smallskip 

From the $\partial$-closedness and (1), it follows that the collection 
$$\bDelta(I):=(\Delta(I|\Sigma)\, )_{\Sigma}\in 
 \bop \cF(I, \emptyset|\Sigma)\subset F(I)\eqno{(\diagcycle.a)}$$
is a cocycle of degree 0 in the complex $F(I)$.
Indeed, 
$$\begin{array}{cl} 
 \bar{d}(\Delta(I|\Sigma)\,)  &=-d(\Delta(I|\Sigma)\,)\\
    & = - \sum_{k\in \Sigma} r_k(\Delta(I|\Sigma)\,)
\end{array}$$
(the first identity by the definition of $\bar{d}$ in (\cpxFI), and 
the second identity by $\partial (\Delta(I|\Sigma)\,)=0$), and 
$$\bar{\rho}_k (\Delta(I|\Sigma)\,)={\rho}_k (\Delta(I|\Sigma)\,)$$
by the definition of $\bar{\rho}$. 
From the definition (\diagcycle.a) and (2) above, we have:
\bigskip 

\sss{\diagcycleProp}
\Prop{\it The $0$-cocycle $\Delta(I)\in F(I)$ satisfies:

(1) If $|I|=2$, then $\bDelta(I)=\Delta(I)\in F(I)$. 

(2) If $S\subset \ctop{I}$, and $I_1, \cdots, I_c$ the corresponding 
segmentation, one has 
$$\tau_S(\bDelta(I)\, )=\bDelta(I_1)\ts\cdots\ts\bDelta(I_c)$$
in $F(I\tbar S)=F(I_1)\ts\cdots\ts F(I_c)$. 
(Recall from (\cpxFI) that 
$\tau_S: F(I)\to F(I\tbar S)$ is the composition of 
$\sigma_{\emptyset, S}: F(I)\to F(I|S)$ and $\iota_S: F(I|S)\to F(I\tbar S)$.)

(3) For $K\subset\ctop{I}$, 
$\vphi_K(\bDelta(I)\,)=\bDelta(I-K)$.}
\bigskip 

Since $\Delta(I)$ and $\bDelta(I)$ depends on $X$, we will write $\Delta_X(I)$ for 
$\Delta(I)$ and $\bDelta_X(I)$ for $\bDelta(I)$. If $|I|=2$, 
$\Delta_X(I)$ is the usual diagonal $\Delta_X$. 
\bigskip 

\sss{\deltadiag}
{\it The diagonal embedding $\delta_*$.}\quad 
Let $X$ be a sequence of varieties on $I=[1, n]$.
Given an element $k\in I$ (we allow $k=1$ or $k=n$)
and an integer $m\ge 1$, let 
$I'$ be the ordered set 
$$\{1, \cdots, k-1, k_1, \cdots, k_m, k+1,
\cdots, n\}
\,\,,$$
where $k$ is repeated $m$ times. 
There is a natural surjection $\lambda: I' \to I$
 which sends $k_j$ to $k$ and is the identity on $I'-\{k_j\}$,
so there is an induced sequence of varieties on $I'$. 
One has the induced sequence of varieties $\lambda^* X$ on $I'$, 
which maps $i\in I'-\{k_j\}$ to $X_i$ and $k_j$ to $X_k$. 
For a subset $\cJ$ of $\ctop{I}$, recall that we have the variety 
$$X_I^\cJ= \prod_{i\in I} X'_i, \quad X'_i=
\left\{
\begin{array}{cl}
\bar{X}_i &\mbox{if $i\not\in \cJ$}, \\
{X}_i &\mbox{if $i\in \cJ$}. 
\end{array}
\right.$$
Similarly for a subset $\cJ'$ of $\ctop{I'}$ one has the associated variety 
$(\lambda^*X)_{I'}^{\cJ'}$ which is abbreviated to $X_{I'}^{\cJ'}$. 
Explicitly,
$$X_{I'}^{\cJ'}=X'_1\times\cdots\times X'_{k-1}\times (X'_{k_1}\times\cdots\times X'_{k_m})
\times X'_{k+1}\times\cdots\times X'_n\,$$
where if $i\not\in\{k_j\}$, 
$$X'_i=
\left\{
\begin{array}{cl}
\bar{X}_i &\mbox{if $i\not\in \cJ'$}, \\
{X}_i &\mbox{if $i\in \cJ'$}. 
\end{array}
\right.\, , $$
and 
$X'_{k_j}$ equals $\bar{X}_k$ if $k_j\not\in \cJ'$, and equals $X_k$ if $k_j\in \cJ'$. 

Assume that the image of $\cJ'$ by $\lambda$ is contained in $\cJ$, and 
that the induced map $\cJ'\to \cJ$ is bijective. 
Then for an element $i\in I'-\{k_1, \cdots, k_m\}\cong I-\{k\}$, 
the variety $X'_i$ is the same for $X_I^\cJ$ and for $X_{I'}^{\cJ'}$, and 
for each $k_j$ there is an open immersion 
$X'_k\injto X'_{k_j}$. Indeed if $k\not\in \cJ$ and $k_j\not\in\cJ'$, then $X'_k=X'_{k_j}=\bar{X}_k$;
if $k\in \cJ$ and $k_j\in\cJ'$, $X'_k=X'_{k_j}={X}_k$; 
if $k\in\cJ$ and $k_j\not\in \cJ'$, then $X'_k=X_k$ and $X'_{k_j}=
\bar{X}_k$. Taking the product of them, one obtains a closed immersion $\delta: X_I^{\cJ} \injto X_{I'}^{\cJ'}$. 

 Recall to a subset $J\subset \ctop{I}$ there corresponds
a closed set $A^\cJ_I\subset X^\cJ_{I}$, and 
$U^\cJ_I$ is its complement. 
\bigskip 

(\deltadiag.1)\Prop{\it 
(1) Let $\cJ\subset \ctop{I}$, $\cJ'\subset \ctop{I'}$ be subsets 
such that $\cJ'\isoto \cJ$, and $J\subset I$, $J'\subset I'$ be
subsets such that $\lambda$ induces a surjection $J'\to J$. 
Then the following square is Cartesian:
$$\begin{array}{ccc}
A^\cJ_I&\injto &X^\cJ_{I}\\
\mapd{}& &\mapdr{\delta} \\
A^{\cJ'}_{I'}&\injto &\phantom{\,\,.}X^{\cJ'}_{I'}\,\,.
\end{array} 
$$
Thus $\delta^{-1}(U^{\cJ'}_{I'})= U^{\cJ}_I$. 

(2) Let $\cJ\subset \ctop{I}$, $\cJ'\subset \ctop{I'}$ be subsets 
such that $\cJ'\isoto \cJ$. Then the closed immersion $\delta$ induces a map of complexes
$$\delta_*: \cF(I, \cJ)\to \cF(I', \cJ')\,.$$
}\smallskip 

{\it Proof.}\quad (1) Obvious from the definitions. 

(2) Let $\{J^0, \cdots, J^r\}$ be the segmentation of 
$I$ by $\cJ$. 
Similarly let $\{J'^0, \cdots, J'^{r}\}$ be the segmentation of 
$I'$ by $\cJ'$.
Then one has surjections ${J'}^{i}\to J^i$, so by (1), 
$\delta^{-1}(U^{\cJ'}_{J'^i})=U^\cJ_{J^i}$. 
Thus $\delta^{-1}\cU(\cJ')=\cU(\cJ)$. 
By (\ZMUpush) there is the induced map of complexes 
$$\delta_*: \Z(X_{I}^\cJ, \cU(\cJ)\, )\to \Z(X_{I'}^{\cJ'}, \cU(\cJ')\, )\,\,.$$
\bigskip 

\sss{\Deltadiag}
{\it The maps $\delta_*$ and $\Delta(\Sigma, \Sigma')$.}\quad
We keep the notation for $I$ and $I'$ from the previous subsection. 
Assume that $\cJ$, $\Sigma$ are disjoint subsets of $\ctop{I}$, and
$\cJ'$, $\Sigma'$ are disjoint subsets of $\ctop{I'}$, satisfying 
the following conditions:
\smallskip 

$\bullet$\quad $\lambda$ induces a bijection $\cJ'\isoto \cJ$.

$\bullet$\quad $\lambda$ induces a bijection $\Sigma'-\{k_1, \cdots, k_m\}\isoto \Sigma-\{k\}$. 
If $1<k<n$, it is also assumed that $\lambda$ induces a bijection 
$\Sigma'\to \Sigma$, in other words, 
if $k\in \Sigma$, then $\Sigma'\cap\{k_1, \cdots, k_m\}\neq \emptyset$, and 
if $k\not\in\Sigma$, then $\Sigma'\cap\{k_1, \cdots, k_m\}= \emptyset$. 
(When $k=1$ or $n$, $\Sigma'\cap\{k_1, \cdots, k_m\}$ may 
be empty or not.)
\smallskip 

Then we will define a map of complexes
$$\cF(I, \cJ|\Sigma)\to \cF(I', \cJ'|\Sigma')\,.$$
 According to cases (there will be Type (0), (I), and (II)), we will give it 
the name $\delta_*$ or $\Delta(\Sigma, \Sigma')$.
\bigskip

(Type 0) Case $\Sigma'\cap\{k_1, \cdots, k_m\}=\emptyset$
($k=1, n$ allowed)\quad  
We have the map
$\delta_*: \cF(I, \cJ|\Sigma)\to  \cF(I', \cJ'|\Sigma')$
given as follows.

First assume $1<k<n$. 
If $\Sigma=\Sigma'=\emptyset$, this is the map defined in (\Deltadiag.2). 
In general, let $I_1, \cdots, I_r$ be the segmentation of $I$ by $\Sigma$; then 
$k$ is contained in some $\ctop{I_i}$.
Let $I'_1, \cdots, I'_r$ be the segmentation of $I'$ by $\Sigma'$. 
Then $I'_j=I_j$ for $j\neq i$ and $I'_i\to I_i$ is a surjection. 
For an element 
$u_1\ts\cdots\ts u_r\in \cF(I_1, \cJ_1)\hts\cdots\hts\cF(I_r, \cJ_r)
=\cF(I, \cJ|\Sigma)$, we let 
$$\delta_*(u_1\ts\cdots\ts u_r)=u_1\ts\cdots\ts\delta_*(u_i)\ts\cdots\ts u_r
\in \cF(I'_1, \cJ'_1)\hts\cdots\hts\cF(I'_r, \cJ'_r)$$
where $\delta_*$ on the right hand side is the map $\cF(I_i, \cJ_i)\to \cF(I'_i, \cJ'_i)$ defined in (\deltadiag.1).
If $k=1, n$ and $\Sigma'\cap\{k_1, \cdots, k_m\}=\emptyset$, one defines $\delta_*$ in a similar manner. 

Note that there are two subcases:
\smallskip 

\quad (0-a) Case $k\not\in \cJ$ ($k=1, n$ allowed). 
Then $\cJ'$ as above is uniquely 
determined.

\quad (0-b) Case $k\neq 1, n$ and $k\in\cJ$.  Then  $\cJ'=(\cJ-\{k\})\cup\{k_j\}$ for $
j=1, \cdots, m$. 
So we write $(\delta_j)_*$ for $\delta_*$. 
\bigskip

(\Deltadiag.1) \Prop{\it  \quad  For the map $\delta_*:\cF(I, \cJ|\Sigma)\to  \cF(I', \cJ'|\Sigma')$, we have:

(1) In cases (a) and (b), $\delta_*$ commutes with $r_{\ell}$
if $\ell\not\in \{k\}\cup \cJ\cup\Sigma$.

(2) If $k \not\in \cJ\cup\Sigma$, then the following diagram commutes:
$$\begin{array}{ccc}
\cF(I,\cJ|\Sigma)&\mapr{\delta_*} &\cF(I', \cJ'|\Sigma') \\
\mapd{r_k}& &\mapdr{r_{k_j}} \\
\cF(I,\cJ\cup\{k\}|\Sigma)&\mapr{(\delta_j)_*} &\cF(I', \cJ'\cup\{k_j\}|\Sigma')\,\,.
\end{array}
$$

(3) If $k\in \cJ$ so that $\cJ=\cJ_0\cup\{k\}$, 
 the following commutes for any $j < j'$:
$$\begin{array}{ccc}
\cF(I,\cJ|\Sigma)&\mapr{(\delta_j)_*} &\cF(I', \cJ_0\cup\{k_j\}|\Sigma')
 \\
\mapd{(\delta_{j'})_*}& &\mapdr{r_{k_{j'}}} \\
\cF(I',\cJ_0\cup\{k_{j'}\}|\Sigma')&\mapr{r_{k_{j}}} &\cF(I', \cJ_0
\cup\{k_j, k_{j'}\}|\Sigma')\,\,.
\end{array}
$$

(4) For $\ell\in \Sigma$, the following diagram commutes:
$$\begin{array}{ccc}
\cF(I,\cJ|\Sigma)&\mapr{\delta_*} &\cF(I', \cJ'|\Sigma') \\
\mapd{\rho_\ell}& &\mapdr{\rho_\ell} \\
\cF(I,\cJ\cup\{\ell\}|\Sigma-\{\ell\})&\mapr{\delta_*} &\cF(I', \cJ'\cup\{\ell\}|\Sigma'-\{\ell\})\,\,.
\end{array}
$$
Similarly $(\delta_j)_*$ commutes with $\rho_\ell$. 
}
\bigskip 

{\it Proof.}\quad  (1), (2) and (4) are obvious from the definitions.

(3) For simplicity of notation consider the case $\Sigma=\emptyset$ and 
$\cJ=\{k\}$ (the general case is parallel). 
Write $k'=k_j$ and $k''=k_{j'}$. 
The map $r_{k''} (\delta_j)_*$ is induced by the map of coverings 
$$\cU(\cJ)=\delta^{-1}(\cU(\{k'\})\, ) \to \delta^{-1}(\cU(\{k', k''\})\,)\,.$$
But $\cU(\cJ)=\{U^\cJ_{[1, k]}, U^{\cJ}_{[k, n]}\}$, and 
$$\begin{array}{cl} 
  \delta^{-1}(\cU(\{k', k''\})) 
&=\delta^{-1}\{U^{\cJ_0\cup\{k', k''\}}_{[1, k']}, 
U^{\cJ_0\cup\{k', k''\}}_{[k', k'']}, 
U^{\cJ_0\cup\{k', k''\}}_{[k'', n]}\} \\
   &   =\{U^\cJ_{[1, k]}, \emptyset, U^{\cJ}_{[k, n]}  \}
\end{array}$$
since one clearly has 
$\delta^{-1}(U^{\cJ_0\cup\{k', k''\}}_{[k', k'']})=\emptyset$.

Let $U_0=U^\cJ_{[1, k]}$, $U_1=\emptyset$, and $U_2=U^{\cJ}_{[k, n]}$. 
Then $\cF(I', \cJ_0\cup \{k', k''\})=\bop_J \Z(U_J)$ where 
$J$ varies over subsets of $\{0, 1, 2\}$. 
Then the map in question, composed with the projection to the 
sum $\bop_{J\not\ni 1} \Z(U_J)$, 
$$\cF(I,\cJ|\{k\})\mapr{r_{k''} (\delta_j)_*}
\cF(I', \cJ_0\cup\{k_j, k_{j'}\}|\Sigma')=\bop_J \Z(U_J)
\mapr{proj}\bop_{J\not\ni 1} \Z(U_J)$$
is the identity map.
The projection to the sum $\bop_{J\ni 1} \Z(U_J)$ is zero by 
$\delta^{-1}(U^{\cJ_0\cup\{k', k''\}}_{[k', k'']})=\emptyset$.
The other map $ (\delta_{j'})_*r_{k'}$ in the diagram has the same description, so the two maps coincide. 
\bigskip

(Type I) Case $k\in \Sigma$ and 
$|\Sigma'|=|\Sigma|$. \quad  Note that $1<k<n$.  One will define a map
$\Delta(\Sigma, \Sigma')
: \cF(I, \cJ|\Sigma)\to  \cF(I', \cJ'|\Sigma')$.
For simplicity assume $\Sigma=\{k\}$, and 
let $I_1, I_2$ be the segmentation of $I$ by $k$.
Let $k':=k_j$ be the element in $\Sigma'$, and 
 $I'_1$, $I'_2$ be the segmentation of 
$I'$ by $\ell$, and 
$${\delta_1}: X_{I_1}^{\cJ_1}\injto X_{I'_1}^{\cJ'_1},\qquad
{\delta_2}: X_{I_2}^{\cJ_2}\injto X_{I'_2}^{\cJ'_2}
$$
be the 
embeddings corresponding to the surjections $I'_i\to I_i$.
Then the map 
$\Delta(\Sigma, \Sigma'): \cF(I, \cJ|\Sigma)\to  
\cF(I', \cJ'|\Sigma')$
is defined  
by $\Delta(\Sigma, \Sigma')(u_1\ts u_2)={\delta_1}_*(u_1)\ts {\delta_2}_*(u_2)$.
That the right hand side is indeed in $\cF(I', \cJ'|\Sigma')$ 
is the first assertion of the following lemma, which is immediately
verified. 
\bigskip

(\Deltadiag.2){\bf Lemma.} \quad{\it  
Let $u_i$ be elements in $\cF(I_i, \cJ_i)$ for $i=1, 2$,
such that $\{u_1, u_2\}$ is properly intersecting 
with respect to $(X/S; I; \{I_i\supset \cJ_i\})$ in the 
sense of (\distcpxcFI.2).
Then $(\delta_i)_*(u_i)\in \cF(I_i, \cJ_i)$, $i=1, 2$, 
are properly intersecting with respect to 
$(X/S; I'; \{I'_i\supset \cJ'_i\})$, and one has
$$\delta_*(u_1\scirc u_2)={\delta_1}_*(u_1)\scirc
 {\delta_2}_*(u_2)$$
in $\cF(I,\cJ)$. }
\bigskip 

The following proposition follows from (\Deltadiag.1) and (\Deltadiag.2).
\bigskip 

(\Deltadiag.3)\Prop\quad{\it Assume we are in case (I); let $\Sigma'=(\Sigma-\{k\})
\cup\{k_j\}$. 

(1) $\Delta(\Sigma, \Sigma')$ commutes with $r_{\ell}$ if $\ell\in \ctop{I}-(\cJ\cup\Sigma)$. 

(2) $\Delta(\Sigma, \Sigma')$ commutes with $\rho_{\ell}$ if $\ell\neq k$. 

(3) The following square commutes: 
$$\begin{array}{ccc}
\cF(I, \cJ|\Sigma)&\mapr{\Delta(\Sigma, \Sigma')} &
\cF(I', \cJ'|\Sigma') \\
\mapd{\rho_k}& &\mapdr{\rho_{k_j}} \\
\cF(I, \cJ\cup\{k\}|\Sigma-\{k\})
&\mapr{(\delta_j)_*} &\cF(I', \cJ'\cup\{k_j\}| \Sigma'-\{k_j\})\,\,.
\end{array}
$$}
\bigskip

(Type II) Case $k\in \Sigma$ and
$|\Sigma'|>|\Sigma|$ ($k=1, n$ allowed).\quad We will define the map 
$$\Delta(\Sigma,\Sigma'): \cF(I, \cJ|\Sigma)\to \cF(I', \cJ'|\Sigma')
$$
as follows. 

First consider the case $1<k<n$. 
For simplicity assume $\Sigma=\{k\}$, the general case being 
similar. Let $I_1, I_2$ be the segmentation of $I$ by $k$, 
and $I'_1,\cdots, I'_{c}$ the segmentation of $I'$ by 
$\Sigma'$. One has 
$\cF(I, \cJ|\Sigma)=\cF(I_1, \cJ_1)\hts \cF(I_2, \cJ_2)$, and 
$$\cF(I', \cJ'|\Sigma')=\cF(I'_1, \cJ'_1)\hts
\cF(I'_2, \emptyset)\hts\cdots\hts\cF(I'_{c-1}, \emptyset)
\hts\cF(I'_{c}, \cJ'_{c})\,\,.$$
Note $I'_2, \cdots, I'_{c-1}$ correspond to constant sequences on $X_k$. 
The map $\Delta(\Sigma, \Sigma')$ is defined by 
$$u_1\ts u_2\mapsto {\delta_1}_*(u_1)\ts\Delta(I'_2)\ts\cdots\ts\Delta(I'_{c-1})\ts
{\delta_2}_*(u_2)$$
where ${\delta_1}_* :\cF(I_1, \cJ_1)\to \cF(I'_1, \cJ'_1)$ is the map 
associated to the surjection $I'_1\to I_1$, and similarly for the 
map ${\delta_2}_*$.  We have used the following lemma. 
\bigskip 

(\Deltadiag.4)\Lem\quad{\it  
Let $u_i$ be elements in $\cF(I_i, \cJ_i)$ for $i=1, 2$,
such that $\{u_1, u_2\}$ is properly intersecting 
with respect to $(X/S; I; \{I_i\supset \cJ_i\})$.
Then 
$$\{{\delta_1}_*(u_1), \Delta(I'_2), \cdots, \Delta(I'_{c-1}), 
{\delta_2}_*(u_2)\, \}$$
is properly intersecting with respect to $(X/S; I'; \{I'_i\supset \cJ'_i\})$
One has 
$${\delta_1}_*(u_1)\scirc \Delta(I'_2)=
\bar{{\delta_1}}_*(u_1)\,\,$$
where $\bar{{\delta_1}}$ is associated 
to the surjection $I'_1\cup I'_2\to I_1$; similarly for 
$ \Delta(I'_{c-1})\scirc {\delta_2}_*(u_2)$. }
\bigskip 

Next assume $k=1$ (the case $k=n$ being similar), and assume 
for simplicity $\Sigma=\emptyset$. 
Let $I'_1, \cdots, I'_{c-1}, I'_c$ be the segmentation of 
$I'$ by $\Sigma'$. 
The surjection $I'_c\to I$ corresponds to a closed 
immersion $\delta': X_I^\cJ\injto X_{I'_c}^{\cJ'}$. 
The map $\Delta(\Sigma, \Sigma')$ is given by 
$$u\mapsto \Delta(I'_1)\ts\cdots\ts\Delta(I'_{c-1})\ts 
\delta'_*(u)\in \cF(I'_1, \emptyset)\hts\cdots\hts\cF(I'_{c-1}, \emptyset)
\hts \cF(I'_c, \cJ')=\cF(I', \cJ'|\Sigma')\,.$$

For $1\le k\le n$, from the definitions and the above lemma, we have:
\bigskip

(\Deltadiag.5) \Prop{\it  \quad  Assume we are in case (II). 

(1)  $\Delta(\Sigma, \Sigma')$ commutes with $r_{\ell}$.
if $\ell\neq k$ for $\ell\neq k$. 

(2) $\Delta(\Sigma, \Sigma')$ commutes with $\rho_{\ell}$
if $\ell\neq k$. 

(3) If $k'=k_j\in \Sigma'$, the following commutes: 
$$\begin{array}{ccc}
\cF(I, \cJ|\Sigma)&\mapr{\Delta(\Sigma, \Sigma')} 
&\cF(I',\cJ'|\Sigma') \\
\mapd{\Delta(\Sigma, \Sigma'-\{k'\})}& &\mapdr{\rho_{k'}} \\
\cF(I',\cJ'|\Sigma'-\{k'\})
 &\mapr{r_{k'}} &\cF(I',\cJ'\cup\{k'\}|\Sigma'-\{k'\})\,\,.
\end{array}
$$}
\bigskip 

(\Deltadiag.6)\quad We conclude this subsection by noting  
of the maps $\delta_*$ and $\Delta(\Sigma, \Sigma')$. 
Let $k'\in \{k_1, \cdots, k_m\}$ and $\lambda': I''\to I'$ be a map such that
$\lambda'$ is surjective and bijective over 
$I'-\{k'\}$. 
Let $\cJ'', \Sigma''$ be subsets of $\ctop{I'}$ satisfying the condition at the beginning of (\Deltadiag)
with respect to $(\cJ', \Sigma')$. 
We thus have the map $\delta_*$ or $\Delta(\Sigma', \Sigma'')$, $\cF(I', \cJ'|\Sigma')\to \cF(I'', \cJ''|\Sigma'')$. 
Then the composition of the maps 
$$\cF(I, \cJ |\Sigma)\mapr{\delta_*}\cF(I', \cJ' |\Sigma')\mapr{\delta_*} \cF(I'', \cJ''|\Sigma'')$$
is shown to coincide with the map $\delta_*$ associated to $\lambda\lambda': I''\to I$. The same for the composition of the maps 
$\Delta(\Sigma, \Sigma')$ and $\Delta(\Sigma', \Sigma'')$. 
The verifications are immediate from the definitions. 
\bigskip

\sss{\defdiag}\quad
We shall define a degree preserving map 
$$\diag: F(I)\to F(I')\eqno{(\defdiag.a)}$$
as the sum of the maps $\delta_*$ or $\Delta(\Sigma, \Sigma')$,
$\cF(I, \cJ|\Sigma)\to \cF(I', \cJ'|\Sigma')$ for $(\cJ, \Sigma)$ and 
$(\cJ', \Sigma')$ satisfying the condition at the beginning of (\Deltadiag) 
(the other components are set to be zero).
We also use the notation $\diag(I, I')$ or $\lambda^*$. 
\bigskip 

(\defdiag.1){\bf Proposition.}\quad{\it 
The map  $\diag$ is a map of complexes.}
\smallskip 

{\it Proof.}\quad  
One must show $d_F\diag=\diag d_F$.
This follows from the propositions in 
(\Deltadiag)
together with examination of the signs, as 
we will elaborate below. 

Recall that $d_F=\bar{d}+ \bar{\rho}$, and
$\bar{d}=\sum \pm 1\ts\cdots\ts (r+(-1)^a\partial )\ts\cdots\ts 1$. 
So $d_F$ is a signed sum of $1\ts\cdots\ts r\ts\cdots\ts 1$, 
$1\ts\cdots\ts \partial\ts\cdots\ts 1$, and $\rho$. 
Of these, $1\ts\cdots\ts \partial\ts\cdots\ts 1$ commutes with the map $\diag$. 
Thus we have only to examine the commutativity of $\diag $ and the signed sum of the maps 
$r$ and $\rho$. 

The composition $d_F\diag$ is the signed sum of the compositions 
$$\cF(I, \cJ| \Sigma) \mapr{f} \cF(I', \cJ' | \Sigma')\mapr{g} \cF(I', \cJ''|\Sigma'')$$
where $f=\delta_*$ or $\Delta(\Sigma, \Sigma')$, and $g=r_k$, $\rho_k$, or $\partial$. 
Specifically, according to the type of $f$, they consist of the following
(we omit the case $g=\partial$): 
\smallskip 

(Type 0-a) For $f=\delta_*$ of Type (0-a), 
\newline \qquad\quad
 $r_\ell \delta_*$\,  ($\ell\not\in \cJ'\cup\Sigma'\cup\{k_1, \cdots, k_m\}$), \quad
 $r_{k_j} \delta_*$\, ($j=1, \cdots, m$), 
\quad $\rho_\ell \delta_*$\,  ($\ell\in \Sigma'$). 
\smallskip 

(Type 0-b) For $f=(\delta_j)_*$ of Type (0-b),
\newline \qquad\quad 
$r_\ell (\delta_j)_*$  ($\ell\not\in \cJ'\cup\Sigma'\cup\{k_1, \cdots, k_m\}$), 
\quad  $r_{k_{j'} } (\delta_j)_*$\, ($j'\neq j$), 
\quad $\rho_\ell  (\delta_j)_*$\, ($\ell\in \Sigma'$). 
\smallskip 

(Type I) For $f=\Delta(\Sigma, \Sigma')$ of Type (I), with
$\Sigma'\cap\{k_1, \cdots, k_m\}=\{k_j\}$, 
\newline\qquad\quad
 $r_\ell \Delta(\Sigma, \Sigma')$  \, ($\ell\not\in \cJ'\cup\Sigma'\cup\{k_1, \cdots, k_m\}$),
\quad $r_{k_i} \Delta(\Sigma, \Sigma')$\,   ($i\neq j$), 
\newline\qquad\quad $\rho_\ell \Delta(\Sigma, \Sigma')$ \, ($\ell\in\Sigma'-\{k_j\}$),
\quad $\rho_{k_j} \Delta(\Sigma, \Sigma')$.
\smallskip 

(Type II) For $f=\Delta(\Sigma, \Sigma')$ of Type (II),
\newline\qquad\quad 
 $r_\ell \Delta(\Sigma, \Sigma')$  \, ($\ell\not\in \cJ'\cup\Sigma'\cup\{k_1, \cdots, k_m\}$),
\quad $r_{k_i} \Delta(\Sigma, \Sigma')$\,   ($k_i\not\in\Sigma'$), 
\newline\qquad\quad $\rho_\ell \Delta(\Sigma, \Sigma')$ \, ($\ell\in\Sigma'-\{k_1, \cdots, k_m\}$),
\quad $\rho_{k_j} \Delta(\Sigma, \Sigma')$ \, ($k_j\in \Sigma'$).
\smallskip 

\noindent 
On the other hand, the composition $\diag d_F$ is the signed sum of the compositions 
$$\cF(I, \cJ| \Sigma) \mapr{g} \cF(I, \cJ_1| \Sigma_1)\mapr{f} \cF(I', \cJ'_1|\Sigma'_1)\,.$$ They consist of (omitting the 
case $g=\partial$):
\smallskip

(Type 0-a) For $f=\delta_*$ of Type (0-a), 
\newline \qquad\quad 
 $ \delta_* r_\ell$\,  ($\ell\not\in\cJ\cup\Sigma$), \quad
 $ \delta_* \rho_\ell$\, ($\ell\in\Sigma$).
\smallskip 

(Type 0-b) For $f=(\delta_j)_*$ of Type (0-b),
\newline \qquad\quad 
$ (\delta_j)_* r_\ell$  ($\ell\not\in\cJ\cup\Sigma$), 
\quad  $  (\delta_j)_* r_k$\, ($j=1, \cdots, m$), 
\newline \qquad\quad
 $ (\delta_j)_*\rho_\ell $\, ($\ell\in \Sigma$),
\quad  $ (\delta_j)_*\rho_k $. \phantom{$j=1, \cdots, m$}
\smallskip 

(Type I) For $f$ of Type (I), 
\newline \qquad\quad 
$\Delta(\Sigma, \Sigma')r_\ell$ \, ($\ell\not\in \cJ\cup\Sigma$),
\quad 
$\Delta(\Sigma-\{\ell\}, \Sigma')\rho_\ell$ \, ($\ell\in\Sigma'$).
\smallskip 

(Type II) For $f$ of Type (II), 
\newline \qquad\quad 
$\Delta(\Sigma, \Sigma')r_\ell$ \, ($\ell\not\in \cJ\cup\Sigma$),
\quad 
$\Delta(\Sigma-\{\ell\}, \Sigma')\rho_\ell$ \, ($\ell\in\Sigma'$).
\smallskip

(i) Among these, one has identities 
$$r_\ell f= f r_\ell\qquad 
(\mbox{$f=\delta_*$ or $\Delta(\Sigma, \Sigma')$}\,)$$ 
and 
$$\rho_\ell f= f \rho_\ell\qquad 
(\mbox{$f=\delta_*$ or $\Delta(\Sigma, \Sigma')$}\,)$$ 
by (\Deltadiag.1), (1), (4) and (\Deltadiag.3), (1), (3);
 one also verifies that the signs assigned to $r_\ell f$ 
and $f r_\ell$ (resp.  $\rho_\ell f$ and $f \rho_\ell$) are
equal.

(ii) One has ${r_{k_j}}\delta_*=(\delta_j)_* r_k$ by 
(\Deltadiag.1), (2). To show that the signs assigned to the two terms are also 
equal, assume, say, $\Sigma=\emptyset$
(the general case being similar), and let $I_1, I_2$ be the
segmentation of $I$ by $k$, $\cJ_i=\cJ\cap I_i$. 
Then ${r_{k_j}}\delta_*$ appears in $d_F\diag$ with the sign 
$-\{|\cJ_2|\}$; 
indeed the sign $\{|\cJ'_{>k_j}|\}=\{|\cJ_2|\}$ is assigned to $r_{k_j}$ in defining $r$ (cf. (\cFI)), and minus sign 
is further assigned in defining $\bar{d}$ (cf. (\cpxFI.b))\,.
On the other hand,   
$(\delta_j)_* r_k$ appears in $\diag d_F$ with the same sign, since 
$\{|\cJ_{>k}|\}=\{|\cJ_2|\}$ is assigned to $r_{k}$ for $r$, and minus sign 
is further assigned for $\bar{d}$. 
(Recall our notation $\{c\}:=(-1)^c$ for an integer $c$.)

(iii) One has, for $j<j'$, $r_{k_{j'} } (\delta_j)_*=r_{k_{j} } (\delta_{j'})_*$ by (\Deltadiag.1), (3). 
We show that they appear in $d_F\diag$ with {\it opposite}
signs, so they cancel each other. 
Indeed, when $\Sigma=\emptyset$ and $\cJ=\cJ_0\cup\{k\}$, let  
$I_1, I_2$ be as above and $\cJ_i=\cJ\cap\ctop{I_i}$. 
The map $r_{k_{j'} } (\delta_j)_*$ appears with sign $-\{|(\cJ_0\cup\{k_j\})_{>k_{j'}}|\}=-\{|\cJ_2|\}$
while $r_{k_{j} } (\delta_{j'})_*$ appears with 
$-\{|(\cJ_0\cup\{k_{j'}\})_{>k_j}|\}=-\{|\cJ_2|+1\}$. 

(iv) By (\Deltadiag.3), (2) one has, for $\Delta(\Sigma, \Sigma')$ of type (I), 
$\rho_{k_j} \Delta(\Sigma, \Sigma')=(\delta_j)_*\rho_k$.
We show that they appear with the same signs. 
Assume $\Sigma=\{k\}$, and let $I_1, I_2$, $\cJ_i$ be as above. 
For $u=u_1\ts u_2\in \cF(I_1, \cJ_1)\hts \cF(I_2, \cJ_2)$, 
with $u_1\in \cF(I_1, \cJ_1)^p, a=|\cJ_1|+1$ 
(namely of bi-degree $(a, p)$) and 
$u_2\in \cF(I_2, \cJ_2)^q, b=|\cJ_2|+1$, 
recall that $\bar{\rho}_k(u)=(-1)^{aq}u_1\scirc u_2$
by (\cFISigma.e) and (\cpxFI.c). 
The sign for $\bar{\rho}_{k_j}\Delta(\Sigma, \Sigma')(u)$ is the 
same, since $\Delta(\Sigma, \Sigma')$ preserves the bi-degrees 
$(a, p)$ and $(b, q)$. 

(v) For each map $\Delta(\Sigma, \Sigma'): 
\cF(I, \cJ| \Sigma) \to \cF(I', \cJ' | \Sigma')$ of type 
(I) or (II), and $k_j\not\in\Sigma'$, one has by (\Deltadiag.5), (3)
$$\rho_{k_j} \Delta(\Sigma, \Sigma'\cup\{k_j\})
= r_{k_j}\Delta(\Sigma, \Sigma')\,.$$
We claim that the two maps appear with {\it opposite} signs in 
$d_F\diag$, canceling each other. 

We will consider three typical cases (the general case is
parallel to one of these cases.)
First assume $\Sigma=\{k\}$, $\Sigma'=\{k_{j'}\}$, and $j<j'$. 
Let $J=[k_j, k_{j'}]$, and $I_1, I_2$, $\cJ_i$ be as before. 
Consider an element $u=u_1\ts u_2\in \cF(I_1, \cJ_1)\hts \cF(I_2, \cJ_2)$, with $u_1$ of bi-degree $(a, p)$, and 
$u_2$ of bi-degree $(b, q)$. 
Then 
$$\Delta(\Sigma, \{k_j, k_{j'}\})(u)=
(\delta_1)_*u_1\ts\Delta(J)\ts (\delta_2)_* u_2\,,$$
using the notation introduced where we defined the map
$\Delta(\Sigma, \Sigma')$. 
Then  $\rho_{k_j}\Delta(\Sigma, \{k_j, k_{j'}\})(u)$ appears with the sign 
$$(-1)^{a\cdot 0}\{\epsilon(\Delta(J))+\epsilon(u_2)\}=\{b+q-1\}$$
by (\cFISigma.e) and (\cpxFI.c), 
while $r_{k_j}\Delta(\Sigma, \{k_{j'}\})(u)$ appears with
$-\{\epsilon(u_2)\}=-\{b+q-1\}$. 

Second assume $\Sigma=\{k\}$ and $\Sigma'=\{k_{j'}\}$ with $j'<j$. 
Let $J=[k_{j'}, k_{j}]$. 
For $u=u_1\ts u_2$ as above, 
$\Delta(\Sigma, \{k_j', k_{j}\})(u)=
(\delta_1)_*u_1\ts\Delta(J)\ts (\delta_2)_* u_2$.
In this case $\rho_{k_j}\Delta(\Sigma, \{k_{j'}, k_{j}\})(u)$ appears with sign $\{q\}\{b+q-1\}$ (since $\Delta(J)$ is of
bi-degree $(1, 0)$, the sign $(-1)^{aq}$ for $\rho$ is $(-1)^q$), 
and $r_{k_j}\Delta(\Sigma, \{k_{j'}\})(u)$ appears with 
$-\{|\cJ_2|\}=\{b\}$. 

The third case is where $\Sigma=\{k\}$ and $\Sigma'=\{k_{j'},
k_{j''}\}$ and $j'<j<j''$. 
Then 
$$\Delta(\Sigma, \Sigma'\cup\{ k_{j}\})(u)=
(\delta_1)_*u_1\ts\Delta(J_1)\ts\Delta(J_2)\ts (\delta_2)_* u_2\,,$$ 
with $J_1=[k_{j'}, k_j]$, $J_2=[k_{j}, k_{j''}]$.
Then $\rho_{k_j}\Delta(\Sigma, \Sigma'\cup\{ k_{j}\})(u)$ 
appears with sign $\{b+q-1\}$, while 
$r_{k_j}\Delta(\Sigma, \Sigma'\})(u)$ appears with 
$-\{b+q-1\}$. 
This completes the proof of the proposition.
\bigskip 

(\defdiag.2){\it Generalization.}\quad 
Let $\lambda: I'\to I$ be a surjective map of finite 
ordered sets of cardinality $\ge 2$. 
Given a sequence of varieties $X$ on $I$, let $\lambda^*X$ be 
the induced sequence of varieties on $I'$.   
Generalizing the map (\defdiag.a), we define 
the map of complexes 
$\lambda^*=\diag=\diag(I, I'): F(I)\to F(I')$
by writing $\lambda=\lambda_1\lambda_2\cdots\lambda_r$, with 
each $\lambda_i$ satisfying the condition in (\deltadiag), and 
setting $\lambda^*=\lambda^*_r\cdots\lambda^*_1$. 
The well-definedness follows from (\Deltadiag.6).
Note that we have functoriality: If $\lambda':I''\to I'$ is another surjective map, then we have
${\lambda'}^*\lambda^*=(\lambda\lambda')^*: F(I)\to F(I'')$. 
\bigskip

\sss{\secdiag.6}
\Prop{\it  \quad Let $\lambda: I'\to I$ be a surjective map of finite ordered sets. 
The map  $\diag:F(I)\to F(I')$ 
is compatible with $\vphi$ and $\tau $ in the following sense. 

 (1) Let $\ell\in \ctop{I'}$. 
If $\sharp \lambda^{-1}(\lambda(\ell))=1$, then 
$\vphi_{\lambda(\ell)}\diag(I, I')=\diag(I-\{\ell\}, I'-\{\ell\}
)\vphi_{\ell}$, namely the following square commutes:
$$\begin{array}{ccc}
F(I)&\mapr{\diag(I, I')} &F(I') \\
\mapd{\vphi_{\lambda(\ell)}}& &\mapdr{\vphi_{\ell}} \\
F(I-\{k\})&\mapr{\diag(I-\{\ell\}, I'-\{\ell\})} &F(I'
-\{\ell\})\,\,.
\end{array}
$$
If $\sharp \lambda^{-1}(\lambda(\ell))>1$, then 
$\vphi_\ell\diag(I, I')=\diag(I, I'-\{\ell\})$, namely the diagram 
$$\begin{array}{ccc}
 F(I)  &\mapr{\diag(I, I')} &F(I') \\
\phantom{aaaaaaaa}\scriptstyle{\diag(I, I'-\{\ell\}) }   &\!\!\!\!\!\searrow & \mapdr{\vphi_\ell} \\
   &  &F(I'-\{\ell\})
 \end{array}$$
commutes. 

(2) Let $\ell\in \ctop{I'}$. 
If $\lambda (\ell)\neq 1, n$, let 
 $I_1, I_2$ be the segmentation of $I$ by $\lambda(\ell)$, and 
 $I'_1, I'_2$ be the segmentation of $ I'$ by $\ell$. 
One then has a commutative diagram: 
 $$\begin{array}{ccc}
F(I)&\mapr{\diag(I, I')} &F(I') \\
\mapd{\tau_{\lambda(\ell)}}& &\mapdr{\tau_\ell} \\
F(I_1)\ts F(I_2)&\mapr{} &F(I'_1)\ts  F(I'_2)\,\,, 
\end{array}
$$
where the lower horizontal arrow is $\diag(I_1, I'_1)
\ts \diag (I_2, I'_2)$. 

If $\lambda(\ell)=1$ (so that $\lambda^{-1}(1)$ has cardinality 
$\ge 2$), let $I'_1, I'_2$ be the 
segmentation of $I'$ by $\ell$. Then the following 
diagram commutes:
$$\begin{array}{ccc}
F(I)&\mapr{\diag(I, I')} &F(I') \\
\mapd{\diag(I, I'_2)}& &\mapdr{\tau_\ell} \\
F(I'_2)&\mapr{} &F(I'_1)\ts F(I'_2)\,\,.
\end{array}
$$
The lower horizontal map is $u\mapsto u\ts \bDelta(I'_2)$. 
(Note $I'_2$ parametrizes a constant sequence of varieties, so one 
has the element $\bDelta(I'')\in F(I'')$.) 
Similarly in case $\lambda(\ell)=n$. }
\smallskip 

{\it Proof.}\quad 
(1) For $\ell$ with $\sharp \lambda^{-1}(\lambda(\ell))=1$, and for the map 
$\delta_*$ in (\deltadiag), it follows immediately from the definitions that the following square commutes:
$$\begin{array}{ccc}
\cF(I,\cJ|\Sigma)&\mapr{\delta_*} &\cF(I', \cJ'|\Sigma') \\
\mapd{\pi_\ell}& &\mapdr{\pi_\ell} \\
\cF(I-\{\ell\},\cJ|\Sigma)&\mapr{\delta_*} &\cF(I'-\{\ell\}, \cJ'|\Sigma')\,\,.
\end{array}
$$
The same holds for the map $\Delta(\Sigma, \Sigma')$ in (\Deltadiag). 
If $\sharp \lambda^{-1}(\lambda(\ell))>1$, the following diagram commutes:
$$\begin{array}{rcc}
\cF(I,\cJ|\Sigma)&\mapr{\delta_*} &\cF(I', \cJ'|\Sigma') \\ 
\scriptstyle{\delta_* } &\!\!\!\!\!\searrow & \mapdr{\pi_\ell} \\
   &  &\cF(I'-\{\ell\}, \cJ'|\Sigma')\,.
 \end{array}$$
The same for the map $\Delta(\Sigma, \Sigma')$. 
The assertion follows from these. 

(2) If $\lambda(\ell)\in\Sigma$ and $\ell\in\Sigma'$, then 
for the map $\Delta(\Sigma, \Sigma')$, 
the following square commutes:
$$\begin{array}{ccc}
\cF(I, \cJ|\Sigma)&\mapr{\Delta(\Sigma, \Sigma')} 
&\cF(I', \cJ'|\Sigma')\\
\mapd{}& &\mapdr{} \\
\cF(I_1, \cJ_1|\Sigma_1)\ts\cF(I_2, \cJ_2|\Sigma_2)
&\mapr{} &\cF(I'_1, \cJ'_1|\Sigma'_1)\ts\cF(I'_2, \cJ'_2|\Sigma'_2)
\,\,.
\end{array}
$$
Here $\cJ_i=\cJ\cap \ctop{I}_i$, $\Sigma_i=\Sigma\cap \ctop{I}_i$, 
and similarly for $\cJ'_i$ and $\Sigma'_i$. The vertical inclusions are the 
canonical ones, and the lower horizontal arrow is 
$\Delta(\Sigma_1, \Sigma'_1)\ts \Delta(\Sigma_2, \Sigma'_2)$. 
Taking the sum over $\Delta(\Sigma, \Sigma')$ we obtain the 
claim. 
The proof for the case $\lambda(\ell)=1, n$ is similar. 
\bigskip

{\bf Acknowledgments.} We would like to thank S. Bloch, B. Kahn,
 and P. May for helpful discussions. 
\bigskip

{\bf References.}
\smallskip




\RefBlone Bloch, S. : Algebraic cycles and higher $K$-theory, 
Adv. in Math. 61 (1986), 267 - 304. 

\RefBltwo  ---  : The moving lemma for higher Chow groups, J. Alg. Geom. 
3 (1994), 537--568. 

\RefBlthree --- :  
Some notes on elementary properties of higher chow groups, including functoriality properties and cubical chow groups, preprint on Bloch's home page.

\RefCH Corti, A. and Hanamura, M. : Motivic decomposition and intersection Chow
groups I, Duke Math. J. 103 (2000), 459-522.


\RefFu Fulton, W. : Intersection Theory, Springer-Verlag, Berlin, New York, 1984.



\RefHaoneone  Hanamura, M.  :  
Mixed motives and algebraic cycles I, Math. Res.
Letters 2(1995), 811-821.

\RefHaonetwo  --- :
Mixed motives and algebraic cycles II, 
Invent. Math. 158(2004), 105-179.

\RefHaonethree  --- :
Mixed motives and algebraic cycles III, 
Invent. Math. 158(2004), 105-179.
Math. Res. Letters  6(1999), 61-82.

\RefHatwo   --- :  Homological and cohomological motives of algebraic varieties,
Invent. Math. 142 (2000), 319-349. 

\RefHathree  --- : Quasi DG categories and mixed motivic sheaves, preprint. 


\bigskip 

{}
\bigskip

{\it Department of Mathematics, 

Tohoku University, 

Aramaki Aoba-Ku, 980-8587, Sendai, Japan}

\end{document}